\documentclass{article}%
\usepackage{amssymb}
\usepackage{amsfonts}
\usepackage{amsmath}
\usepackage{graphicx}%
\providecommand{\U}[1]{\protect\rule{.1in}{.1in}}

\begin{document}

\title{A Survey of the Development of Geometry up to 1870\thanks{This monograph was
written up in 2008-2009, as a preparation to the further study of the early
geometrical works of\ Sophus Lie and Felix Klein at the beginning of their
career around 1870. The author apologizes for possible historiographic
shortcomings, errors, and perhaps lack of updated information on certain
topics from the history of mathematics. Comments or corrections from the
reader are most welcome, and may contribute to an improved later edition. }}
\author{Eldar Straume\\Department of mathematical sciences\\Norwegian University of Science and Technology (NTNU)\\N-9471 Trondheim, Norway}
\maketitle

\begin{abstract}
This is an expository treatise on the development of the classical geometries,
starting from the origins of Euclidean geometry a few centuries BC up to
around 1870. At this time classical differential geometry came to an end, and
the Riemannian geometric approach started to be developed. Moreover, the
discovery of non-Euclidean geometry, about 40 years earlier, had just been
demonstrated to be a "true" geometry on the same footing as Euclidean
geometry. These were radically new ideas, but henceforth the importance of the
topic became gradually realized. As a consequence, the conventional attitude
to the basic geometric questions, including the possible geometric structure
of the physical space, was challenged, \ and foundational problems became an
important issue during the following decades.

Such a basic understanding of the status of geometry around 1870 enables one
to study the geometric works of Sophus Lie and Felix Klein at the beginning of
their career in the appropriate historical perspective.

\end{abstract}
\tableofcontents

\section{Euclidean geometry, the source of all geometries}

At the end of the 18th century, the notion of geometry was largely synonymous
with Euclidean geometry, namely the classical Greek geometry which had
prevailed for more than 2000 years. In the 17th century, Kepler, Galileo and
Newton were leading figures in the Copernican revolution which had paved the
way for the birth of modern science, and moreover, which finally abandoned the
long lasting doctrines and the supreme authority of Aristotle among the
scholars. Certainly, the 17th century was also a century of great advances in
mathematics, such as the rise of calculus, but still the basic role of
Euclidean geometry was not challenged. Indeed, in his famous work Principia
Mathematica (1687) Newton was careful to recast his demonstrations in
geometric terms, and analytical calculations are almost completely missing. In
particular, Newton formulated his basic laws for the universe in the framework
of Euclidean geometry.

However, during the 18th century another dominating authority had established
himself in the intellectual world, namely the German philosopher Immanuel Kant
(1724--1804). He maintained that Euclidean geometry was the only absolute
geometry, known to be true \emph{a priori }in our mind as an inevitable
necessity of thought, and no other geometry was thinkable. Moreover, Kant
regarded the Newtonian universe as the true model of the physical space,
supported by our endowed intuition about space and time, and independent of
experience. In reality, the Kantian doctrine on the nature of space and
geometry hampered the development of science, until the outburst of the
inevitable geometric "revolution" in the 19th century led to the discovery of
non-Euclidean geometries and radically changed geometry as a mathematical science.

\subsection{Early geometry and the role of the real numbers}

Geometry is in fact encountered in the first written records of mankind. But
what are the origins of the Euclidean geometry? In his truly remarkable work,
organized into 13 books usually referred to as \emph{Euclid's} \emph{Elements,
}Euclid (325--265 BC) presented a large part of the geometric and algebraic
knowledge of Babylonian, Egyptian and Greek scholars at his time. He did this
in a deductive style, which has become known as the axiomatic method of
mathematics. Geometry as developed in the \emph{Elements} is usually referred
to as \emph{synthetic} geometry. The basic undefined geometric objects are
points, lines, planes, whereas notions such as line segments, angles and
circles can be defined, and moreover, the undefined notion of congruence
expresses "equality" among them. Starting from five basic postulates and five
common notions (or rules of logic), a logical chain consisting of 465
propositions were deduced and presented in the \emph{Elements}. For the
convenience of the reader, let us be more specific and state the postulates
--- somewhat modernized--- as follows:

\begin{description}
\item $E1$ : There is a unique line passing through any two distinct points.

\item $E2$ : Any segment on a line may be extended by any given segment.

\item $E3$ : For every point $O$ and every other point $P$, there is circle
with center $O$ and radius $OP.$

\item $E4$ : All right angles are congruent to each other.

\item $E5$ : If a line intersects two other lines so that the two interior
angles on one side of it are together less than two right angles, then the two
lines will meet at a point somewhere on that side.
\end{description}

The ancient geometry flourished for about 1000 years, primarily in Greece and
Alexandria in Egypt. Besides Euclid some of the prominent men from this epoch
are Thales of Miletus (ca. 600 BC), Pythagoras (585--501 BC), Plato (429--348
BC), Eudoxus of Cnidus (408--355? BC), Aristotle (384--322 BC), Archimedes of
Syracuse (287--212 BC), Appolonius of Perga (262--190 BC), and Pappus of
Alexandria (290--350 AD). The Pythagorean school played a crucial role,
developing number theory which they also linked to geometry, number mystics
and music theory.\ In fact, Pythagoras himself is said to have coined the
words "philosophy" and "mathematics". The Platonic Academy in Athens became
the mathematical center of the world and, for example, Aristotle and Eudoxus
had been students at this academy.

Many ot the ancient mathematical texts have been lost, including advanced
works by Euclid and Appolonius. For example, Appolonius introduced the names
\emph{ellipse}, \emph{parabola} and \emph{hyperbola} in his famous treatise
Conic Sections, but some of the volumes are irretrievably lost. Moreover,
Euclid's three-volume work called the \emph{Porisms}, apparently on advanced
geometry, is also lost. However, the \emph{Elements} of Euclid have survived
throught the centuries by manuscript copies which have transmitted the
geometric truths to new generations, and with minor changes. Through Euclids's
authority, the \emph{Elements} teach geometry with a definite approach, but
with no motivation even of the most sophisticated terms. Moreover, previous or
alternative approaches, which could have filled some of the gaps, are ignored.
The first printed edition was in 1482, and over thousand editions have
appeared since then.\ In effect, the \emph{Elements} have remained in use as a
textbook practically unchanged for more than 2000 years, and up to modern
times the first six books have served as the student's usual introduction to geometry.

It is somewhat surprising that Euclid did not use numbers in his geometry.
Namely, the geometric objects such as line segments, angles, areas etc. are
not measured by numbers, but they are related to each other in terms of
congruence, similarity and proportion (or ratio). The idea behind this is,
loosely speaking, that congruent objects have the same "shape and size", that
is, they are \emph{similar} and have the same \emph{magnitude.} In the
simplest case of line segments, denoted by letters $a,b,c,..$, congruence
means "same magnitude".
\ \ \ \ \ \ \ \ \ \ \ \ \ \ \ \ \ \ \ \ \ \ \ \ \ \ \ \ \ \ \ \ \ \ \ \ \ \ \ \ \ \ \ \ \ \ \ \ \ \ \ \ \ \ \ \ \ \ \ \ \ \ \ \ \ \ \ \ \ \ \ \ \ \ \ \ \ \ \ \ \ \ \ \ \ \ \ \ \ \ \ \ \ \ \ \ \ \ \ \ \ \ \ \ \ \ \ \ \ \ \ \ \ \ \ \ \ \ \ \ \ \ \ \ \ \ \ \ \ \ \ \ \ \ \ \ \ \ \ \ \ \ \ \ \ \ \ \ \ \ \ \ \ \ \ \ \ \ \ \ \ \ \ \ \ \ \ \ \ \ \ \ \ \ \ \ \ \ \ \ \ \ \ \ \ \ \ \ \ \ \ \ \ \ \ \ \ \ \ \ \ \ \ \ \ \ \ \ \ \ \ \ \ \ \ \ \ \ \ \ \ \ \ \ \ \ \ \ \ \ \ \ \ \ \ \ \ \ \ \ \ \ \ \ \ \ \ \ \ \ \ \ \ \ \ \ \ \ \ \ \ \ \ \ \ \ \ \ \ \ \ \ \ \ \ \ \ \ \ \ \ \ \ \ \ \ \ \ \ \ \ \ \ \ \ \ \ \ \ \ \ \ \ \ \ \ \ \ \ \ \ \ \ \ \ \ \ \ \ \ \ \ \ \ \ \ \ \ \ \ \ \ \ \ \ \ \ \ \ \ \ \ \ \ \ \ \ \ \ \ \ \ \ \ \ \ \ \ \ \ \ \ \ \ \ \ \ \ \ \ \ \ \ \ \ \ \ \ \ \ \ \ \ \ \ \ \ \ \ \ \ \ \ \ \ \ \ \ \ \ \ \ \ \ \ \ \ \ \ \ \ \ \ \ \ \ \ \ \ \ \ \ \ \ \ \ \ \ \ \ \ \ \ \ \ \ \ \ \ \ \ \ \ \ \ \ \ \ \ \ \ \ \ \ \ \ \ \ \ \ \ \ \ \ \ \ \ \ \ \ \ \ \ \ \ \ \ \ \ \ \ \ \ \ \ \ \ \ \ \ \ \ \ \ \ \ \ \ \ \ \ \ \ \ \ \ \ \ \ \ \ \ \ \ \ \ \ \ \ \ \ \ \ \ \ \ \ \ \ \ \ \ \ \ \ \ \ \ \ \ \ \ \ \ \ \ \ \ \ \ \ \ \ \ \ \ \ \ \ \ \ \ \ \ \ \ \ \ \ \ \ \ \ \ \ \ \ \ \ \ \ \ \ \ \ \ \ \ \ \ \ \ \ \ \ \ \ \ \ \ \ \ \ \ \ \ \ \ \ \ \ \ \ \ \ \ \ \ \ \ \ \ \ \ \ \ \ \ \ \ \ \ \ \ \ \ \ \ \ \ \ \ \ \ \ \ \ \ \ \ \ \ \ \ \ \ \ \ \ \ \ \ \ \ \ \ \ \ \ \ \ \ \ \ \ \ \ \ \ \ \ \ \ \ \ \ \ \ \ \ \ \ \ \ \ \ \ \ \ \ \ \ \ \ \ \ \ \ \ \ \ \ \ \ \ \ \ \ \ \ \ \ \ \ \ \ \ \ \ \ \ \ \ \ \ \ \ \ \ \ \ \ \ \ \ \ \ \ \ \ \ \ \ \ \ \ \ \ \ \ \ \ \ \ \ \ \ \ \ \ \ \ \ \ \ \ \ \ \ \ \ \ \ \ \ \ \ \ \ \ \ \ \ \ \ \ \ \ \ \ \ \ \ \ \ \ \ \ \ \ \ \ \ \ \ \ \ \ \ \ \ \ \ \ \ \ \ \ \ \ \ \ \ \ \ \ \ \ \ \ \ \ \ \ \ \ \ \ \ \ \ \ \ \ \ \ \ \ \ \ \ \ \ \ \ \ \ \ \ \ \ \ \ \ \ \ \ \ \ \ \ \ \ \ \ \ \ \ \ \ \ \ \ \ \ \ \ \ \ \ \ \ \ \ \ \ \ \ \ \ \ \ \ \ \ \ \ \ \ \ \ \ \ \ \ \ \ \ \ \ \ \ \ \ \ \ \ \ \ \ \ \ \ \ \ \ \ \ \ \ \ \ \ \ \ \ \ \ \ \ \ \ \ \ \ \ \ \ \ \ \ \ \ \ \ \ \ \ \ \ \ \ \ \ \ \ \ \ \ \ \ \ \ \ \ \ \ \ \ \ \ \ \ \ \ \ \ \ \ \ \ \ \ \ \ \ \ \ \ \ \ \ \ \ \ \ \ \ \ \ \ \ \ \ \ \ \ \ \ \ \ \ \ \ \ \ \ \ \ \ \ \ \ \ \ \ \ \ \ \ \ \ \ \ \ \ \ \ \ \ \ \ \ \ \ \ \ \ \ \ \ \ \ \ \ \ \ \ \ \ \ \ \ \ \ \ \ \ \ \ \ \ \ \ \ \ \ \ \ \ \ \ \ \ \ \ \ \ \ \ \ \ \ \ \ \ \ \ \ \ \ \ \ \ \ \ \ \ \ \ \ \ \ \ \ \ \ \ \ \ \ \ \ \ \ \ \ \ \ \ \ \ \ \ \ \ \ \ \ \ \ \ \ \ \ \ \ \ \ \ \ \ \ \ \ \ \ \ \ \ \ \ \ \ \ \ \ \ \ \ \ \ \ \ \ \ \ \ \ \ \ \ \ \ \ \ \ \ \ \ \ \ \ \ \ \ \ \ \ \ \ \ \ \ \ \ \ \ \ \ \ \ \ \ \ \ \ \ \ \ \ \ \ \ \ \ \ \ \ \ \ \ \ \ \ \ \ \ \ \ \ \ \ \ \ \ \ \ \ \ \ \ \ \ \ \ \ \ \ \ \ \ \ \ \ \ \ \ \ \ \ \ \ \ \ \ \ \ \ \ \ \ \ \ \ \ \ \ \ \ \ \ \ \ \ \ \ \ \ \ \ \ \ \ \ \ \ \ \ \ \ \ \ \ \ \ \ \ \ \ \ \ \ \ \ \ \ \ \ \ \ \ \ \ \ \ \ \ \ \ \ \ \ \ \ \ \ \ \ \ \ \ \ \ \ \ \ \ \ \ \ \ \ \ \ \ \ \ \ \ \ \ \ \ \ \ \ \ \ \ \ \ \ \ \ \ \ \ \ \ \ \ \ \ \ \ \ \ \ \ \ The
magnitude of $a$ relative to $b$ is represented by the proportion $a:b$, but
let us write it as a ratio $a/b$. The idea of proportion is the clue to
proving many theorems, since by subdivision of the geometric figure similitude
can be reduced to congruence. The ultimate procedure for defining and handling
the ratios $a/b$ in Greek mathematics, as presented in books V --VI of the
\emph{Elements}, is generally attributed to Eudoxus. Previously, the
Pythagoreans had failed at this point, having stumbled into the existence of
incommensurable ratios and hence the discovery of irrationality, which
influenced the further development of Greek geometry in a fundamental way.

In more detail, the Pythagorean approach was to associate a (rational) number
$m/n$ to each ratio $a/b$ of magnitudes, based on their belief that $a $ and
$b$ are always commensurable, that is, they are both integral multiples of
some suitably small $c$, say $a=mc$, $b=nc$. However, taking $a$ and $b$ to be
the side and diagonal of a square, and assuming $n/m$ is reduced so that $n$
and $m$ are relatively prime, this example yields the contradictory identity
$2n^{2}=m^{2}$. Thus they had, in fact, encountered a pair of incommensurable
line segments. Hippasus (ca. 500 BC) is credited with the discovery, but most
likely he first discovered the incommensurability of the side and diagonal of
a regular pentagon.

Later, Aristotle pointed out how rational numbers can approximate any ratio by
a "take away from each other" procedure, which is closely related to the
continued fraction construction in modern mathematics. Namely, if $a_{0}$ and
$a_{1}$ are the magnitudes to be compared and $a_{0}>a_{1}$ say, one finds
successively unique positive integers $m_{1},m_{2},...,$ so that
$a_{p-1}=m_{p}a_{p}+a_{p+1}$ with $a_{p}>a_{p+1}$, for $p=1,2$ etc. Let us
write
\begin{equation}
a_{0}/a_{1}\longleftrightarrow\left[  m_{1},m_{2},m_{3},..\right]
=m_{1}+\frac{1}{m_{2}+\frac{1}{m_{3}+...}} \label{approx}%
\end{equation}
which also indicates that the procedure may never stop, in which case
$a_{0}/a_{1}$ is incommensurable. Otherwise, $a_{p+1}=0$ for some $p$, and the
procedure is known as the \emph{Euclidean algorithm}. Then $a_{0}=ma_{p}%
,a_{1}=na_{p}$ for some $m$ and $n,$ and the above finite continued fraction
ending with $m_{p}$ is just the rational number $m/n$ measuring the
commensurable ratio $a_{0}/a_{1}$.

In a more elegant way, however, Eudoxus's approach to \emph{ratios }is the
principal source to the modern view of real numbers, irrational or not.
Regarding \emph{ratio }as an undefined relation between magnitudes, he
declares that $a/b=$ $c/d$, if for any two positive integers $n$ and $m$, one
of the following three relations holds for the pair $(a,b)$, namely
\begin{equation}
i)\text{ }ma=nb\text{ \ \ or }\ \ ii)\text{ }ma>nb\text{ \ \ or \ }iii)\text{
}ma<nb\text{,} \label{equality}%
\end{equation}
if and only if the corresponding relation holds for the pair $(c,d)$.
Furthermore, he declares that $a/b$ is less than $c/d$ if for some $m$\ and
$n$
\[
na<mb\text{ \ \ and \ \ }nc>md
\]

For real numbers $a,b,c,d$ in the modern sense, the definition (\ref{equality}%
) of equality between ratios certainly yields $a/b=$ $c/d$ as numbers.
Furthermore, Aristotle demonstrates the density of the rationals, locating an
incommensurable ratio $a/b$ by comparing it with a special sequence of
approximating commensurable ratios, similar to the modern definition of real
numbers by decimal fractions. On the other hand, Eudoxus compares $a/b$ with
all commensurable ratios and thus anticipates the modern view that a real
number is determined by its order relations with respect to all rationals.

However, in the \emph{Elements} the existence of the ratio $a/b$ of two given
magnitudes $a,b$ is not questioned. What is needed is the existence of
integers $m,n$ so that $ma>b$ and $nb>a$, which was taken for granted. In the
19th century it became clear that this property, referred to as the
\emph{Archimedian axiom}, has to be postulated, directly or indirectly. This
is one of many examples illustrating the Greek philosophy of mathematics
which, despite the use of the deductive method, is very different from the
modern ideal of an axiomatic system, where all conclusions are strictly
deduced from a few fundamental axioms. For Plato or Aristotle, geometric
objects are real and knowable and their basic properties are, anyhow, settled
and may not be explicitly announced.\ 

\ Another major contribution of Eudoxus is his method of \emph{exhaustion},
which is elaborated in Book XII of the \emph{Elements}. Namely, by using a
concept close to the idea of integration in modern calculus, he shows how to
subdivide a known magnitude into decreasingly small pieces whose totality
approaches that of an unknown magnitude. He uses the method to show, for
example, that the volume of a pyramid is one third of the volume of the prism
with the same base and height. The method can be used to compute areas and
volumes bounded by specific curves and surfaces. In 1906 an unknown treatise
of Archimedes called \emph{The method} was discovered, in which he did not
repudiate "infinitesimal" methods. It is likely that he made progress beyond
Eudoxus in this direction, but the improved tools must have failed to meet the
rigour of the \emph{Elements} and hence eliminated by Euclid. The ideas of
Eudoxus and Archimedes were, indeed, anticipating the integral calculus
initiated by Newton and Leibniz in the 17th century.

Skillful techniques\ were developed to apply the theory of proportions, and
the ancient use of the ratios came very close to the segment arithmetic
introduced by Hilbert (1899). Namely, by choosing a "unit" segment $e$, so
that $a$ (or rather its magnitude) can be identified with the ratio $a/e$, the
algebraic\ operations addition and multiplication could be effectuated by
compass and ruler constructions. A well known theorem, due to Menelaus of
Alexandria (ca. 70--130 AD), is the following formula involving the product of
three ratios:\ \ \ \ \
\begin{equation}
\frac{AP}{PB}\frac{BO}{QC}\frac{CR}{RA}=-1 \label{Menelaus}%
\end{equation}
where $P,Q,R$ \ are points on the (possibly extended) edges $AB,BC,CA$ of a
triangle $ABC$, respectively. The theorem says that the three points are
collinear if and only if the identity (\ref{Menelaus}) holds.

\subsubsection{Geometric algebra, constructivism, and the real numbers}

The ancient scholars certainly developed some mathematical rigor and logical
analysis through the problem of incommensurability.The area formula for a
rectangle proved by the Pythagoreans is valid only when the sides are
commensurable, but following Eudoxus the formula holds in general. The
traditional opinion, say up to the beginning of the 20th century, was that the
Greeks created \emph{geometric algebra} by translating algebraic relations
into geometry. For example, the Pythagoreans solved quadratic equations by a
geometric procedure involving the notion of area. Namely, the equations%
\begin{equation}
A=(a+x)x\text{ \ \ and \ \ }A=(a-x)x \label{conics}%
\end{equation}
express the area of two rectangles arising from a rectangle with sides $a$ and
$x$, by adding the \emph{excess} (hyperbole) or subtracting the \emph{defect}
(ellipse) represented by the square of side $x$. However, for the last 50
years no historian would phrase the origins of geometric algebra in this way.

On the other hand, quadratic equations correspond to the geometric problems
which the Greeks could study by means of the ruler and compass constructions,
still an indispensable tool in the modern education of plane geometry. Today
we can say that the ratios $a/b$ obtained in this way, starting from a unit
element $e$, will only yield the so-called \emph{constructible} numbers, that
is, obtained from the rationals by the successive addition of square roots.
Some modern writers (e.g. Hartshorne [2000: 42]) claim the above geometric
approach to algebra prevented Euclid and other ancient scholars from
conceiving of real numbers beyond the constructible ones, for example,
$\sqrt[3]{2}$ or transcendental numbers such as $\ $ $\pi$.

The numbers $\sqrt[3]{2}$ and $\pi$ are the solutions of two famous geometric
problems of antiquity which remained unsolved, namely the "dublication of the
cube", and the "quadrature of the circle" which Anaxagoras (499-428 BC) first
attempted to solve. Eudoxus would have approximated such irrational "numbers"
by rationals, but we must remind ourselves that in Greek mathematics the
irrational numbers did, in fact, not have the status of being "numbers". After
all, the irrationality of $\pi$ was first proved by Lambert in the 18th century.

The problem Eudoxus solved with his theory of proportions was not really
understood until the late 19th century, more than 2000 years after the
\emph{Elements} were written down. During this long period a clear conception
of the nature of the real numbers seemed to be missing, and even doubts about
the soundness of irrationals were expressed by some scholars. For example, the
German algebraist Michael Stifel (1487--1567), who discovered logarithms and
was the first to use the term "exponent", \ argued in 1544 that "just as an
infinite number is not a number, so an irrational number is not a true
number". \ Finally, in 1871 Richard Dedekind reexamined the ancient problems
on incommensurables, and with his epoch-making essay \emph{Continuity and
irrational numbers} (1872), he established the theory of real numbers on a
logical foundation and without the extraneous influence of geometry. It should
be mentioned, however, that the Dedekind cut construction is essentially the
same idea as Eudoxus used.

It is interesting to observe that the ancient Euclidean constructivism has
survived up to modern times, manifesting itself in the belief that mathematics
should deal only with constructible numbers and with a finite number of
operations. The famous Berlin professor Kronecker, who Klein and Lie met
during the fall 1869, was the first to doubt the significance of
non-constructive existence proofs. He is well known for his remark that "God
created the integers, all else is the work of man", and he opposed the use of
irrational numbers. As late as in 1886, when Klein's previous student
Lindemann lectured on his proof from 1882 that $\pi$ is transcendental,
Kronecker complimented Lindemann on a beautiful proof, but as he added, it
proved nothing since transcendental numbers did not exist.

In their study of equations like (\ref{conics}) the ancient scholars also came
close to the development of coordinate (or analytic) geometry. For example, by
putting $A=y^{2}$ in (\ref{conics}) the equations describe a hyperbola and an
ellipse in the coordinates $x$ and $y$. However, with the Greek passion for
geometry, Menaechmus (390--320 BC), who was a student of Eudoxus, was led to
the observation that these curves are plane sections of a cone. His study was
continued by Archimedes, and Appolonius finished the project with his
celebrated treatise \emph{Conic Sections}, where he also set forth the
principal properties of \emph{conjugate }diameters. In a way coordinates were
used, but always in an rather awkward language dictated by the geometry.
Therefore, today we regard analytic geometry as originating from Descartes and
Fermat in the early 17th century.

\subsubsection{The downfall of the ancient geometry}

The Second Alexandrian School (around 300 AD), with the mathematicians
Diophantus and Pappus, continued the tradition dating back to Pythagoras and
brought again fame to Alexandria which lasted for another century or two.
Pappus was actually on the track of a new type of geometric truths, naturally
belonging to projective geometry. These are statements involving only points,
lines and their incidence relations. According to the celebrated
\emph{Pappus's theorem}, if $A,B,C$ and $A^{\prime},B^{\prime},C^{\prime}$ are
points on different lines $l$ and $l^{\prime\text{ }}$respectively, and
$AB^{\prime}$ denotes the line joining $A$ and $B^{\prime}$ etc., then the
three points of intersection
\begin{equation}
P=AB^{\prime}\cap A^{\prime}B\text{, \ \ }Q=BC^{\prime}\cap B^{\prime}C\text{,
\ \ }R=CA^{\prime}\cap C^{\prime}A\text{\ \ } \label{Pappus}%
\end{equation}
lie on the same line. Pappus was the last of the great geometers from
classical Greek mathematics, and he stated this result as an exercise in one
of his books. Most likely Pappus knew about Menelaus's theorem (cf.
(\ref{Menelaus})), which can be used to solve his exercise.

But the downfall of Greek geometry was unrelenting; it came with the decadence
of the classical Greek culture around the 5th century AD. At this time the
flow of written records and oral traditions carrying the unofficial
mathematical knowledge was suddenly broken, probably due to political events
or pressure of Roman culture. In addition, the surviving literature,
influenced by the selective role played by Euclid's \emph{Elements}, was
unable to inspire further mathematical creativity.

Stepping forward to the early Middle Ages and the Islamic golden age, one
finds that progress in geometry beyond Euclid's \emph{Elements} and the
Alexandrian School was still rather modest. Certainly, there were some
extensions and refinements of the Euclidean postulates, as well as attempts to
prove the parallel postulate and further studies of conic sections. Thus,
there seemed to be no remaining challenges from the ancient geometry. Stepping
forward another century or two, the works of Kepler, Desargues and Pascal, in
the spirit of Pappus, also remind us that the ancient geometry was not capable
of much further extension.

\subsubsection{The ancient geometry: Its failures and its final
algebraization}

The idea of transformations is absent in Euclid's \emph{Elements}, and
therefore the \emph{Elements} study the properties of individual triangles
from a static viewpoint only. This must have limited the scope of geometric
thought and its interaction with related sciences such as kinematics and
mechanics. This is all the more astounding since transformations were, in
fact, known and used long before Euclid and, for example, the modern notion of
\emph{symmetry} played an important role in the geometry of Thales in the 6th
century BC. Even stone age decorations witness that it was known very early to
mankind. But geometric transformations amount to changing of figures according
to specific rules, so they were eliminated from Euclid's work as they seemed
to belong to mechanics rather than geometry. This is certainly far from the
modern viewpoint where motions and deformations of triangles are ideal study
object of kinematic geometry.

Despite its fame as the outstanding example of a deductive theory, the claim
that all propositions are deduced logically from the definitions, axioms or
postulates, is unfounded and modern criticism have shown essential gaps, as
pointed out earlier. Without the basic idea of a transformation the
\emph{Elements} use artificial methods, which cannot be justified from the
basic postulates. For example, to circumvent or demonstrate congruency between
triangles, Euclid applies a \emph{superposition method} --- an unstated
congruence axiom--- which allows him to move a triangle from one place to
another. Tacitly assuming this he establishes the properties usually referred
to as the congruence propositions SSS (side-side-side) and SAS
(side-angle-side). Unfortunately, this "magic" method is still found in school
geometry textbooks nowadays.

So, for many reasons, during the 19th century there was a growing conviction
that the classical geometry needed a thorough upgrading of its logical
fundament. Based on the joint efforts of many previous geometers, David
Hilbert (1862-1943) finally presented in his \emph{Foundations of geometry}
(1899) and subsequent works the modern axiomatization of the classical
geometries (Euclidean, hyperbolic, or projective). In particular, the
Euclidean space was given a sound basis using about 15 axioms, in the spirit
of Euclid, with congruence via its SAS property postulated as an axiom.

In analytic geometry (see below), the real Cartesian plane $\mathbb{R}^{2}$ is
used as a model of the Euclidean plane. Thus the geometry study in the plane
is reduced to algebra, and geometric properties depend on properties of the
real numbers, justified by the Cantor--Dedekind axiom dating back to 1872. In
fact, in their study of axiomatic systems and dependency relations, Hilbert
and others introduced the Cartesian plane $\mathbb{F}^{2}$ over various
(ordered) fields $\mathbb{F}$, even skew fields. Points and lines are
introduced as in $\mathbb{R}^{2}$, but the geometric properties of these
planes reflect the algebraic properties of $\mathbb{F}$. For example, Pappus's
theorem (see (\ref{Pappus})) holds in $\mathbb{R}^{2}$, and it is an
interesting result due to Hilbert that the field $\mathbb{F}$ is commutative
if and only if the theorem holds in $\mathbb{F}^{2}$. \ \ \ \ \ \ \ \ \ \ \ \ \ \ \ \ 

Finally, it should be noted that critical questions with regard to Euclid's
\emph{Elements} are associated with two major events in the development of
geometry in the 19th century, in which Felix Klein was largely involved around
1870. The first event is concerned with the fifth postulate $E5$, the
so-called \emph{parallel postulate}, and the discovery of non-Euclidean
geometry.\ Whereas the first four postulates are local in nature and easy to
accept, the parallel postulate is less intuitive since one cannot "see" what
happens indefinitely far out in the plane. It was precisely the doubts about
the nature of this postulate which spurred the discovery, and early works of
Klein in 1871--72 contributed to the understanding of the new geometry which
he called \emph{hyperbolic} geometry (see Chapter 4).

The other event was the appearance of the notion of symmetry group in
geometry, which heavily uses transformations and thus provides an alternative
approach to the classical geometry. This new viewpoint is exemplified in
Klein's Erlanger Programm (1872). For example, the congruences in Euclidean
geometry are the \emph{rigid motions} and they constitute the symmetry group
of the geometry. Conversely, in Klein's approach a geometry is largely
characterized by its group. Namely, by focusing attention on the group itself,
geometric properties are precisely those which are invariant under the group.
But certainly the idea of transformation groups and the study of their
possible structures extend far beyond the known classical geometries at that
time. Sophus Lie's general study of continuous groups arose naturally from his
geometric experiences in the early 1870's, in those years when Klein and Lie
were in close contact and mutually stimulated each other.

\subsection{The decline of pure geometry and rise of analytic geometry}

Recall that since antiquity the language of algebra was largely provided by
geometry, and maybe pure algebra made little progress since the ancients
believed the theory of ratios could not be granted algebraically. So, it
seemed that only persons unaware of the Greek scruples would be in the
position to resume and develop genuine algebra. But a mathematical symbolic
language was still missing, and this may have hampered the progress during the
Middle Ages. The creator of such a language, and usually regarded as the
father of modern algebraic notation, is the French lawyer and mathematician
Francois Vi\`{e}te (1540--1603). Some of the symbols he introduced are still
in use today.

Following Vi\`{e}te and greatly influenced by Pappus, there was a revival of
mathematics in France in the 17th century, with leading figures such as Girard
Desargues (1591-1661), Ren\'{e} Descartes (1596--1650), Pierre Fermat
(1601--1665), and Blaise Pascal (1623--1662). In the hands of Descartes and
Fermat, the symbolic analysis led to an entirely new way of investigating
mathematical problems, namely the analytic or coordinate approach to geometry.
Actually, Menaechmus (390--320 BC) and Appolonius came close to such an
approach in their study of conic sections. But, whereas Appolonius had failed
by attaching the coordinates to the conic itself, Descartes and Fermat
initiated the new approach which attaches a coordinate system to the
underlying plane, thus enabling them to study the different figures in their
mutual relations. Therefore, the beginning of the analytic geometric approach
is generally attributed to Descartes and Fermat.

Actually, both Euclid and Descartes unified algebra and geometry, but their
approaches were converse to each other. With his basic mathematical work
\emph{La Geometrie }(1637), Descartes initiated the algebraization of geometry
by associating to each point in the Euclidean plane a pair $(x,y)$ of real
numbers called coordinates, namely the (signed) distances from the point to
two fixed perpendicular axes. In fact, Descartes considered skew coordinate
systems as well, but we shall not do so. The pair $(x,y)$ became known as
\emph{Cartesian }(or \emph{rectangular}) coordinates, and the totality of
pairs as the \emph{Cartesian plane}, usually denoted by $\mathbb{R}^{2}$ for
set theoretic reasons. Henceforth, geometric figures and their properties
could be expressed and studied in terms of coordinates; for example, Descartes
and Fermat represented a curve by an equation $f(x,y)=0$ in the variables
$x,y$. Thus, various geometric locii could now be described in the Cartesian
plane in terms of algebraic or transcendental functions and equations.
\ Notice, however, that the Cartesian distinction between geometrical and
mechanical curves corresponds to the terms algebraic and transcendental in
modern terminology. To begin with, orthogonality of the coordinate axes was
not necessarily assumed, Originally, Descartes the coordinate axes could be
more general, were not necessarily perpendicular, but We remark that
perpendicularity of the coordinate axes was not originally assumed,

This was the birth of \emph{ccordinate geometry}, usually referred to
as\emph{\ analytic geometry} after Lacroix introduced the term for the first
time in his famous two-volume textbook \emph{Trait\'{e} de calcul ..
}(1797--98). But it was also the birth of \emph{algebraic geometry}, in the
sense that the equations involve the coordinates in a purely algebraic way.
Equations of type $ax+by+c=0$ represent the straight lines, whereas those of
degee 2%
\begin{equation}
ax^{2}+bxy+cy^{2}+dx+ey+f=0 \label{quadratic}%
\end{equation}
represent curves which were found to be the familiar conic sections. Descartes
himself initiated a classification of algebraic curves according to the degree
of the equation.\ In the Cartesian plane, problems of higher degree, even
transcendental or with 3 variables, could also be handled in a similar way.

As a consequence of all this, however, the ancient pure geometry fell into
oblivion, at least temporarily, because the new approach turned out to be a
much more powerful method of proof and discovery. Therefore, the algebraic and
analytic methods continued to dominate geometry almost to the exclusion of
synthetic methods. In effect, the Euclidean plane or space was virtually
replaced by the Cartesian plane $\mathbb{R}^{2}$ or 3-space $\mathbb{R}^{3}$,
respectively. In reality, however, this identification rests on a kind of
"continuity" postulate, tacitly accepted since the days of Descartes, saying
that to each point of the line there corresponds a real number, and
conversely. This kind of subtlety, however, was not fully understood until
late in the 19th century.

In the 18th century, Euler was a leading mathematician, writing on essentially
everything, so he combined algebra, analysis and geometry to solve many types
of problems. In solid analytic geometry, his formulae for translation and
rotation of the axes, in terms of the so-called Euler angles, are still well
known and used today. About 100 years after Descartes and Fermat attempted a
unified treatment of binary quadratics (\ref{quadratic}) and conics, Euler
gave for the first time a unified treatment of the general quadratic
equations
\begin{equation}
ax^{2}+by^{2}+cz^{2}+2dxy+2eyz+2fzx+gx+hy+iz+k=0 \label{quadric}%
\end{equation}
which involve three variables and have up to ten terms. His work indicated
that the equation can be reduced by transformation to the five canonical forms
of quadrics, namely an ellipsoid, two types of hyperboloids and two types of
paraboloids, but he did not list all the degenerate types.

Another leading 18th century mathematician, born three decades after Euler,
was Lagrange, whose work resembles Euler's in its elegance and generality, but
neither of the two were typical geometers. At some occations the analyst
Lagrange even boasted of his omission of diagrams or figures. In 1773 Lagrange
turned his attention to the basic problem in solid geometry, namely the
geometry of four points. They span a tetrahedron, and with one point chosen to
be the origin Lagrange was seeking analytic formulae for its various geometric
invariants, such as area, volume, center of gravity, and centers and radii of
the inscribed and circumscribed spheres. His "tetranometry" presented maybe
one of the earliest associations of linear algebra with analytic geometry, but
as Lagrange himself put it, the importance of the work lay more in the point
of view than in the substance (cf. Boyer[1956], Chap. III).

On the other hand, the invention of infinitesimal calculus by Newton and
Leibniz in the later half of the 17th century did not only divert attention
from pure geometry, but even geometry as a whole. Instead, by exploring the
consequences of the fundamental theorem of calculus, a new type of equations
arose, namely differential and integral equations, which have ever since
played an important role in the modern developments of mathematics, analytic
mechanics, and the natural sciences. In reality, however, one important aspect
of geometry profited largely from the new tools in analysis, namely the metric
study of Euclidean geometry, with all its embedded curves and curved surfaces.
So, here we are witnessing the birth of \emph{differential geometry}, to which
Euler and Monge gave many basic contributions.

But Lagrange's interests tended towards physics, and with his famous
\emph{Mechanique analytique} (1788) he developed the whole subject of
mechanics, from the time of Newton, into a branch of mathematical analysis,
starting from a few basic principles and using the theory of differential
equtions. Mechanics may be regarded as the geometry of a (3+1)-dimensional
space, Lagrange remarked, but he did not initiate such a multidimensional
coordinate geometry.

Among synthesists as well as analysts, the greatest geometer of the 18th
century was Monge. We shall encounter him many times later, notably for his
basic contributions to differential geometry and projective geometry. Lagrange
must have envied him for his fertility of imagination and geometrical
innovation, but they both seemed to have realized, more than anyone before,
that analysis and geometry can be combined into a very useful alliance. So,
geometry was rapidly approaching a new stage, towards the turn to the 19th
century, when Paris became the major scene for many historical events,
political as well as scientific.

\subsection{The advent of the new geometries}

During the 19th century geometry developed roughly in three major directions

\qquad\qquad\qquad$(i)$ differential geometry,

$\qquad\qquad\qquad(ii)$\ projective geometry,

$\qquad\qquad\qquad(iii)$\ non-Euclidean geometry, \newline and in the
following chapters we shall discuss them separately in some detail, up to
around 1870, say. It is natural to use this year as a limit for our review,
and for many reasons. At this time projective geometry had established itself
as the central topic of geometry, in fact, as the new geometry of the century.
Since the beginning of the century, the methodology of analytic geometry had
been successfully developed during the decades prior to 1870, to the extent of
being the "golden age" of analytic geometry (cf. also Boyer[1956]). Classical
differential geometry also came to an end around 1870, as the Riemannian
geometric approach became known among geometers and was facing a rapid
development. Moreover, the discovery of non-Euclidean geometry, about 40 years
earlier, had just been demonstrated to be a "true" geometry on the same
footing as Euclidean geometry. These were radically new ideas, but henceforth
the importance of the topic became gradually realized. This was the status of
geometry anno 1870, when Felix Klein og Sophus Lie came to Paris in the
spring, at the beginning of their mathematical career.

On the other hand, with the 1870's a new era started for mathematics as a
whole, where the distinction between the various fields such as algebra and
geometry became increasingly more difficult. Geometry was clearly dominated by
the analytic approach, but there were also geometers who opposed the idea of
reducing geometric thinking to analytic geometry and perhaps relying too much
on physical intuition or experience. The discovery of non-Euclidean geometry,
which became known as \emph{hyperbolic} geometry, also challenged the
conventional attitude to the basic geometric questions, including the possible
geometric structure of the physical space. Therefore, from the synthesist's
viewpoint, the foundations of projective geometry as well as Euclidean and
hyperbolic geometry, maybe in contrast to analysis, still needed a thorough
revision of its basic concepts and postulates. Consequently, the foundational
problems became an important issue during these decades, and for the sake of
completeness we have also included a brief account on these topics\ in Chapter 5.

\ \ \ \ \ \ \ 

Now, let us return to the year 1787, when Lagrange came to Paris and joined
the group of leading mathematicians, such as Monge, Laplace, and Legendre. A
great flowering of French mathematics was about to start, in the spirit of the
French revolution and the new Republic (1789--99), followed by the interlude
of the emperor Napoleon. The \'{E}cole Normale and \'{E}cole Polytechnique
were founded in 1794--95, and these schools were given a leading role with
regard to the higher education.

Among the faculty of the new schools we find for example Lagrange, Monge,
Carnot and Lacroix. In fact, as a favorite of Napoleon, Monge also served as
the director of \'{E}cole Polytechnique, responsible for the education of
military engineers and officers. Geometry played an important role in the
mathematics curriculum, and with his lectures during the years 1795--1809,
Monge's enthusiasm inspired pupils like Brianchon, Dupin, Servois, Biot,
Gergonne, Poncelet and others. They contributed to the further development of
Monge's ideas in differential geometry, but also to the rediscovery and
revival of projective geometry, which had somehow fallen into oblivion since
the days of Desargues and Pascal 150 years earlier. \ \ \ \ \ \ \ \ \ \ \ \ \ \ \ \ \ \ \ \ \ \ \ \ \ \ \ \ \ \ \ \ \ \ \ \ \ \ \ \ \ \ \ \ \ \ \ \ \ \ \ \ 

The new generation of mathematicians also included prominent figures like
Fourier, Cauchy, Poisson and Chasles, but perhaps the innovations in geometry
were, after all, less than in some other areas. An important step was taken by
the establishment of a periodical at \'{E}cole Polytechnique, which certainly
made publication easier than before, since the traditional publication
channels were only the journals at the science academies. Even better, in 1810
an earlier student of Monge, J. D. Gergonne (1771--1859), had just retired as
an artillery officer, and turning to mathematics he founded and edited his own
Gergonne's Annales, the first periodical devoted entirely to mathematics. The
journal terminated when he finally retired in 1832, but fortunately, in 1826
the first purely mathematical journal in Germany, Crelle's Journal, had been
established in Berlin.

Truly, Paris was the center from which the new spirit of analytic geometry
spread to the rest of the world. Former students of Monge, S.F. Lacroix
(1765--1843) and J.B. Biot (1774--1862), are well known for their popular
textbooks on the topic, which also influenced writers in other countries.
Lacroix was said to be the most prolific writer in "modern times", and when
Pl\"{u}cker reported that his introduction to analytic geometry took place in
1825, he referred to the 6th edition of Biot's textbook (cf. Boyer [2004:
245]). But numerous textbooks similar to these appeared in many countries, see
for example Salmon's treatise [1865a] which first appeared in 1862, and its
popular German translation [1865b] which Klein used.

Differential geometry in the tradition of Euler, Monge, Dupin and others
continued to flourish in France in the 19th century. But after Monge, the next
great step was actually made by Gauss in Germany, a master of many
mathematical disciplines. His major work [1828] on curved surfaces also
established him in the forefront of geometry. On the French side, however,
they did not keep an eye on what the Germans were doing until around 1850, but
during the remaining two decades of the pre-Riemannian era differential
geometry progressed steadily with participants also from England, Italy and
other countries. (For a review, see Reich [1973].)

On the other hand, the French geometers also missed the discovery of
non-Euclidean geometry in the 1820's. This landmark in the history of geometry
is generally attributed to Bolyai in Hungary, Lobachevsky in Russia, as well
as Gauss. We refer to Chapter 4\ for more information, including the long
prehistory of the topic.

Finally, turning to projective geometry, let us start with the remark that the
history of geometry is full of discoveries and rediscoveries, as well as
rivalry, priority conflicts, and controversies. This is amply exemplified by
the revival of projective geometry at the turn of the century, initiated by
Monge's pupil L.N.M. Carnot (1753--1823) with his major work
\emph{G\'{e}om\'{e}trie de position} (1803). Although he includes a brief
section on coordinate geometry, in fact the most general view of coordinate
systems since Newton's time, his main purpose is to "free geometry from the
hieroglyphics of analysis". Monge himself abided by the joint use of analysis
and pure geometry, but gradually several outstanding mathematicians held the
opinion that synthetic geometry had been unfairly and unwisely neglected in
the past, and now they would make an effort to revive and extend that
approach. The champion of the synthetic method was Poncelet, a previous
officer of Napoleon. With his classic \emph{Trait\'{e} des propri\'{e}t\'{e}s
projectives des figures} (1822)\emph{\ }he is said to have introduced
projective geometry as a new discipline. But he had strong opponents like
Gergonne and Chasles, who headed the analytic trend and its use of algebra,
and they were also joined by the foremost analysts in Germany, namely
M\"{o}bius and Pl\"{u}cker.

Although the reception of Poncelet's work was rather poor in France, his ideas
were followed up in Germany. In fact, his strict synthetic approach was taken
over by Steiner, the first of a German school who favored strict geometric
methods to the extent of even detesting analysis. Steiner's noble goal was to
develop projective geometry as a unification of the classical geometry,
whereas his compatriot von Staudt wanted to establish projective geometry
independent of Euclidean geometry and its metric concepts. Von Staudt almost
succeeded in the 1860's, but around 1870 a flaw was discovered and traced back
to the implicit usage of the Euclidean parallel postulate. However, this was
later remedied by others, and remaining fundamental questions about projective
spaces were settled during the following decades by German and Italian geometers.

The controversy between the proponents of the two geometric approaches lasted
for many decades. In retrospect, there were in fact good reasons for this
since the methods of analysis were incomplete and even logically unsound. The
pure geometer rightly questioned the validity of the analytic proofs and would
credit them merely with suggesting the results. The analysist, on the other
hand, could retort only that the geometric proofs were clumsy and not so
elegant. (cf. Kline [1972: Chap.35, \S 3]).

In reality, analytic geometry is based upon the Cartesian geometry, and the
subtle distinction between the Euclidean plane (appropriately defined) and the
real Cartesian plane $\mathbb{R}^{2}$ critically depends on the properties of
the real numbers $\mathbb{R}$. However, the Cantor--Dedekind axiom (1872)
finally placed the arithmetization of analytic geometry upon a solid logical
foundation. Therefore, we may well regard 1872 as the terminal year of
classical analytic geometry, 200 years after Descartes and Fermat (cf. also
Boyer[1956], Ch.IX).

\section{ Classical differential geometry}

Classical differential geometry is a term used since 1920 for pre-Riemannian
differential geometry in the tradition of Euler, Clairaut, Monge, Gauss,
Dupin, say, to the end of the 1860's. Curves and surfaces in 3-space can be
rather complicated geometric object, but the development of calculus has
provided infinitesimal techniques to analyze these objects and to distinguish
between and classify the abundance of various geometric forms. Some of the
simplest ones are lines, circles, planes, and spheres, and these are also used
for comparison reason or for approximations of the more general ones. The
basic notion of curvature, essentially a measure of the deviation from
linearity, was introduced for these purposes. For modern references, cf. e.g.
Spivak[1979], Vol. II, III, and Rosenfeld[1988].

\subsection{Curves and their geometric invariants}

Let us start with a brief review of curves and their curvature theory. The
ancients investigated many aspects of the conic sections, which are curves of
the second degree. In his study of plane curves Euler continued with a
classification of cubic curves, and moreover, he also listed 146 different
types of quartics. But Clairaut (1831) was the first who published a treatise
on space curves; for him a space curve was the intersection of two surfaces.

In the older literature we encounter many well known geometric curves, often
constructed as the solution of a mechanical problem, and usually they are
expressed by transcendental functions. One of the constructions yields a
\emph{roulette,} namely a curve generated by a curve rolling on another curve.
For example, the focus of a parabola rolling on a straight line traces out a
\emph{catenary}, and a fixed point on a circle rolling on a straight line
traces out a \emph{cycloid.} The catenary has the shape of a flexible chain
suspended by its ends and acted on by gravity. Its equation is
\begin{equation}
y=a\cosh(\frac{x}{a}) \label{catenary}%
\end{equation}
and was obtained by Leibniz. Let us also recall the Newton--Leibniz
controversy which provoked a great rivalry between the British and the
continentals. Johann Bernoulli (1696) proposed to Leibniz, l'H\^{o}pital, and
others, the celebrated \textquotedblleft brachistocrone\textquotedblright%
\ problem, to find the path between two fixed points which minimizes the time
of fall of a point mass acted on by gravity. The problem reached Wallis and
Gregory in Oxford, who could not solve the problem, but when the challenge
reached Newton, he found out, in a few hours, that the brachistocrone is a cycloid.

The \emph{tractrix} is another distinguished type of curve in the history of
mathematic. It is characterized by the property that the length of its tangent
between a line ($y$-axis) and the point of tangency is the constant $a$, which
leads to the equation
\begin{equation}
\pm y=a\operatorname{arcsech}(\frac{x}{a})-\sqrt{a^{2}-x^{2}}.
\label{tractrix}%
\end{equation}
The problem was posed by Leibniz, but it was first studied by Huygens in 1692,
who also gave the curve its name.

Various other curves are associated with a given curve. For example, when
light rays are reflected off a fixed curve, the curve arising as the envelope
of the reflected rays is known as a \emph{caustic} curve. Such curves were
also studied by Huygens (1678), and the Bernoullis, l'H\^{o}pital, and
Lagrange also gave their contributions. Recall that a family of curves
$F(x,y,a)=0$ depending on a parameter $a$ has an enveloping curve which at
each of its points is tangent to some member of the family. Its equation is
found by elimination of $a$ from the equations%
\begin{equation}
\frac{\partial F}{\partial a}(x,y,a)=F(x,y,a)=0. \label{envelope}%
\end{equation}
The \emph{evolute} of a curve can be constructed in this way, namely as the
curve enveloped by the family of normal lines of the given curve. Huygens
(1673) is credited for being the first to study the evolute, but the geometric
idea dates back, in fact, to the work of Appolonius on conics. The reverse of
the evolute is an involute, which is not unique, since different curves may
have the same evolute. We mention that the evolute of the tractrix is the
catenary, but the cycloid is its own evolute (with a shift). Similarly, the
logarithmic spiral, described in polar coordinates $(r,\theta)$ by the
equation $r=ae^{b\theta}$, is congruent to its own evolute. The study of
envelope in calculus dates back to Leibniz, and it is implicit in the early
works on the evolute. As Lagrange also pointed out, a singular solution of a
differential equation is generally an envelope of the integral curves.

The only intrinsic geometric invariant of a curve is its arc-length, and the
arc-length $s$ measured from a starting point provides a natural
parametrization of the curve. However, what is interesting with curves, as
exemplified by the curves described above, is their extrinsic geometry, namely
their metric properties and position relative to the ambient Euclidean plane
or space. For a plane curve parametrized by $s,$ its \emph{curvature }measures
the rate of change of its tangent direction at a given point, that is,%
\begin{equation}
\kappa\ =\frac{\mathrm{d}\alpha}{\mathrm{d}s} \label{curv1}%
\end{equation}
where $\alpha$ is the oriented angle (in radians) between a fixed direction in
the plane and the tangent direction at $p$.

Note that $\alpha=\alpha(s)$ is also the arc-length traced out on a circle by
the unit tangent vector. The sign of $\kappa$ depends on the direction of the
curve and the orientation of the plane. The \emph{radius} \emph{of curvature,}
$\rho=\pm1/\kappa$, is the radius of the osculating circle, which is the
limiting circle passing through three points on the curve and tending to $p$.
In fact, these concepts were already known to Leibniz, who also had similar
thoughts about the osculating sphere and the curvature of a surface (cf.
Rosenfeld [1988: 280]).

For a curve in the plane, the center of the osculating circle is a focal
point, situated at the distance $\rho$ from $p$ on the normal line, and as the
point $p$ varies the focal points trace out the \emph{focal }curve. The latter
curve is, in fact, the evolute of the given curve, but now arising from a
different viewpoint.

Euler took up the subject of space curves in 1775. He extended the definition
(\ref{curv1}) of curvature by taking $\alpha\geq0$ as the length of the
spherical curve traced out by the unit tangent vector. Let $x,y,$and $z$ be
rectangular coordinates in $\mathbb{R}^{3}$ so that $\sqrt{x^{2}+y^{2}+z^{2}}$
measures the distance from the origin. In terms of these coordinates the
arc-length element $\mathrm{d}s$ is subject to the constraint
\begin{equation}
\mathrm{d}s^{2}=\mathrm{d}x^{2}+\mathrm{d}y^{2}+\mathrm{d}z^{2}
\label{Emetric}%
\end{equation}
and the curvature expresses as
\[
\kappa=\frac{\mathrm{d}\alpha}{\mathrm{d}s}=\sqrt{\frac{\mathrm{d}^{2}%
x}{\mathrm{d}s^{2}}+\frac{\mathrm{d}^{2}y}{\mathrm{d}s^{2}}+\frac
{\mathrm{d}^{2}z}{\mathrm{d}s^{2}}}%
\]

Before 1850 the geometric properties of space curves were laboriously
investigated, and many types of curves have been described in the literature.
They were also called curves of \emph{double curvature}, because the curve is
characterized by two scalar functions $\kappa_{i}(s)$, $i=1,2$, namely the
above curvature $\kappa(s)>0$ and the \emph{torsion} $\tau(s)$, classically
known as\emph{\ }the\emph{\ }first and second curvature. A huge step forward
was taken around 1850, with the advent of the Frenet--Serret formulae, which
in modern vector algebra notation read%
\begin{equation}
\frac{\mathrm{d}\mathbf{t}}{\mathrm{d}s}=\kappa\mathbf{n}\text{, \ }%
\frac{\mathrm{d}\mathbf{n}}{\mathrm{d}s}=-\kappa\mathbf{t}+\tau\mathbf{b}%
\text{, \ \ }\frac{\mathrm{d}\mathbf{b}}{\text{ }\mathrm{d}s}=-\tau\mathbf{n}
\label{F-S}%
\end{equation}
Here $(\mathbf{t,n,b)}$ is an orthonormal frame along the the given curve
$\gamma(s)$, actually explained by the above equations, where the velocity
vector $\mathbf{t=}$ $\mathrm{d}\gamma/\mathrm{d}s$ is the unit tangent,
$\mathbf{n}$ is called the \emph{principal normal}$\mathbf{,}$ and
$\mathbf{b=t\times n}$ is the \emph{binormal.} The equations were obtained
independently by J.A. Serret in 1851 and J.F. Frenet (1816--1900) in 1852, but
partly also in Frenet's thesis (1847). Camille Jordan discovered in 1874 their
generalization for curves with curvatures $\kappa_{1},..,\kappa_{n-1} $ in $n$-space.

Geometrically, the space curve is completely determined by the system
(\ref{F-S}). Thus, for given functions $\kappa(s)>0$ and $\tau(s)$ on the
interval $[0,L]$ and a given initial position and direction, the system yields
by integration a unique space curve of length $L$. The plane spanned by
$\mathbf{t}$ and $\mathbf{n}$ is the \emph{osculating plane}, and hence the
torsion is measuring the rate at which the plane is changing. The
\emph{osculating circle,} lying in this plane, has radius equal to the radius
of curvature $\rho=1/\kappa$, as before. The \emph{normal plane} is the plane
spanned by $\mathbf{n}$ and $\mathbf{b}$, and hence the principal normal
$\mathbf{n}$ spans the line common to the osculating plane and the normal plane.

\subsection{Surfaces and their curvature invariants}

Henceforth, let us turn to surfaces in 3-space and those special curves
confined to them. The curvature theory of surfaces is considerably more
complicated that that of curves, mainly because of the interlocking
relationship between their intrinsic and extrinsic geometry. Surfaces are
typically presented and studied via parametrizations, or as level surfaces of
functions. Euler also introduced differential equations to define surfaces,
such as solutions of the equation
\begin{equation}
\text{\ \ }P\mathrm{d}x+Q\mathrm{d}y+R\mathrm{d}z=0. \label{diff2}%
\end{equation}
But a systematic study of surfaces appeared for the first time in Monge's book
[1807], \emph{Applications of analysis to geometry}. Monge introduced notions
like \textquotedblleft families of surfaces\textquotedblright\ and he used his
theory of surfaces to elucidate the solutions of partial differential
equations. His theory of surfaces and his geometric approach to the study of
differential equations inspired many geometers in the 19th century.

Euler worked out the foundations of the analytic theory of curvature of
surfaces. In 1767 he characterized the curvature of a surface at a point $p$
by describing all \emph{normal curvatures} $\kappa_{\theta}$, where the number
$\kappa_{\theta}$ is the signed curvature (see (\ref{curv1})) of the curve,
called the \emph{normal section}, cut out by a plane $P_{\theta}$
perpendicular to the surface. The plane $P_{\theta}$ contains the line $l$
perpendicular to the surface at $p$, so it is determined by the angle $\theta$
which specifies the tangential direction of the curve at $p$.

What did Euler find out about the numbers $\kappa_{\theta}$? If $\kappa
_{\theta}$ is not a constant for all $\theta$, then $\kappa_{\theta}$ takes
its minimum (resp. maximum) value $\kappa_{1}$ (resp.$\ \kappa_{2})$ for
directions $\theta_{1}$ and $\theta_{2}$ which differ by $90^{\circ}$, and by
choosing the zero angle $\theta=0$ so that $\kappa_{0}$ is smallest, Euler's
formula reads%
\begin{equation}
\kappa_{\theta}=\kappa_{1}\cos^{2}\theta+\kappa_{2}\sin^{2}\theta.
\label{curv2}%
\end{equation}
Actually, Euler expressed his formula in terms of curvature radii, and the
above formula is a modification by Dupin (1837). The numbers $\kappa_{1}$ and
$\kappa_{2}$ became known as the \emph{principal curvatures}, and the
corresponding directions $\theta_{i}$ are the \emph{principal directions}.
Dupin was a student of Monge at Ecole Polytechnique in Paris. His geometric
ideas were influencial for many decades, and below we shall occasionally
return to some of his achievements.

More generally, in 1776 J.B. Meusnier (1754--1793) extended Euler's result to
finding the curvature of the curve $C$ cut out by any plane $P$ through $p $.
In this case, let $C$ have tangent direction $\theta$ at $p$, let $P_{\theta}$
be the corresponding normal section, and let $\varphi<\pi/2$ be the (dihedral)
angle between the planes $P_{\theta}$ and $P$. Then the curvature
$\kappa=\tilde{\kappa}_{\varphi}$ of the curve $C$ at $p$ is determined by
Meusnier's identity (for a proof, see (\ref{curv3}))%
\begin{equation}
\tilde{\kappa}_{\varphi}\cos\varphi=\kappa_{\theta}. \label{Meus}%
\end{equation}

\subsubsection{The Gaussian approach}

We shall deduce the formulae (\ref{curv2}) and (\ref{Meus}), but in the more
general setting due to Gauss. With his paper [1828] on general investigations
of curved surfaces, Gauss unified the surface theories of Euler and Monge, but
he also went much further. This paper, briefly called the
\emph{Disquisitiones},\emph{\ } is maybe the single most important work in the
history of differential geometry. Here he also proved his celebrated results
on the intrinsic geometry of surfaces. Gauss related the curvature of a
surface $S$ with the variation of the tangent planes, or equivalently, the
normal directions. So, as an astronomer himself he borrowed from astronomy the
notion of spherical representation and introduced the so-called \emph{Gauss
map} \emph{\ }%
\begin{equation}
\mathbf{\eta}:S\rightarrow S^{2}=\left\{  (x,y,z)|\text{ }x^{2}+y^{2}%
+z^{2}=1\right\}  \label{Gauss}%
\end{equation}
which, for an orientable surface, specifies the \textquotedblleft
outward\textquotedblright\ normal direction by a continuously varying unit
vector $\mathbf{\eta}_{p}$ perpendicular to $S$ at each $p$.

In his early investigations of curvature, Gauss defined the \emph{total
curvature} of a subset $R$ of $S$ as the area of the image $\mathbf{\eta
(}R\mathbf{)}$, with negative sign if $\mathbf{\eta}$ reverses orientation.
Then, by comparison with the area of $R$ itself and letting $R$ decrease to a
point $p$, he defined the curvature of $S$ at $p$ as the limit (if it exists)%
\[
K=\lim_{R\rightarrow p}\frac{Area(\mathbf{\eta(}R\mathbf{))}}{Area(R)}.
\]
This definition, perhaps not so rigorous, is still useful in special cases;
for example, $K\ $must vanish on a region $R$ whose image $\mathbf{\eta}(R)$
has zero area. Obvious examples are cones or cylinders, where the image of
$\mathbf{\eta}$ lies on a curve of finite length.

Both before and after Gauss various definitions of curvature have been
proposed, say by Euler, Meusnier, Monge and Dupin, but they were never
established and fell into oblivion. Instead, the adopted definition well known
today is due to the more rigorous Gaussian approach leading to the definition
of $K$ as the product of Euler's principal curvatures $\kappa_{i} $, as
follows. Firstly, observe that the tangent plane $T_{p}$ of $S$ at $p$, viewed
as a linear subspace of $\mathbb{R}^{3}$, is a Euclidean plane with the inner
product $\mathbf{v\cdot w}$ inherited from $\mathbb{R}^{3}$. Secondly, by
noticing that $T_{p}$ naturally identifies with the tangent plane of the
sphere $S^{2}$ at $\mathbf{\eta}_{p}$, the differential of $\mathbf{\eta}$ at
$p$ becomes a linear operator today known as the \emph{shape} operator (or
Weingarten map)
\begin{equation}
\mathrm{d}\mathbf{\eta}:T_{p}\rightarrow T_{p}\mathbf{\ } \label{shape}%
\end{equation}
Therefore, associated with this operator, or rather its negative
$-d\mathbf{\eta}$, is the bilinear form%

\begin{equation}
\mathrm{II}_{p}(\mathbf{v},\mathbf{w})=-\mathrm{d}\mathbf{\eta(v)\cdot w}
\label{second1}%
\end{equation}
called the \emph{second fundamental form}. Below we shall see why this form is
actually symmetric, and consequently it has eigenvalues $\kappa_{1}\leq
\kappa_{2}$ and orthonormal eigenvectors $\mathbf{t}_{1},\mathbf{t}_{2}$ in
$T_{p}$ so that
\begin{equation}
\mathrm{II}_{p}(\mathbf{t}_{1},\mathbf{t}_{1})=\kappa_{1}\text{,
\ }\mathrm{II}_{p}(\mathbf{t}_{2},\mathbf{t}_{2})=\kappa_{2}\text{,
\ }\mathrm{II}_{p}(\mathbf{t}_{1},\mathbf{t}_{2})=0\text{\ \ } \label{second2}%
\end{equation}
Moreover, the numbers $\kappa_{i}$ are just Euler's principal curvatures, and
hence the vectors $\mathbf{t}_{i}$ point in the principal directions.

The \emph{Gaussian curvature} $K$ and the \emph{mean curvature} $H$ are
defined to be
\begin{equation}
K=\kappa_{1}\kappa_{2}\text{, \ \ }H=\frac{1}{2}(\kappa_{1}+\kappa_{2}).
\label{H-K}%
\end{equation}
With his \textquotedblleft Theorema Egregium\textquotedblright\ Gauss showed
that $K$ is an intrinsic invariant, that is, it depends only on the geometry
of the surface itself, so that $K$ is unchanged when the surface is bent or
isometrically deformed in any way. Contrary to this, however, $H$ \ reflects
the way the surface is embedded in the ambient space, so it is an extrinsic invariant.

The surface is nowadays called \emph{flat} (resp. \emph{minimal)} if $K=0$
(resp. $H=0$) holds at all points, and the reason for these terms will become
clear later. More generally, if there is a functional relation $W(\kappa
_{1},\kappa_{2})=0$, for example when $K$ or $H$ is constant, the surface is
referred to as a \emph{Weingarten surface} (or W-surface), after J. Weingarten
(1836--1910) who made important contributions, in the 1860's and onward, to
the theory of surfaces in the spirit of Gauss. Note, however, Klein and Lie
used the term \textquotedblleft W-surface\textquotedblright\ with a different
meaning (see letter 3.3.1870). At a single point $p$ the surface is said to be
\emph{hyperbolic, elliptic, parabolic,} or \emph{planar} if (i) $K<0$ or (ii)
$K>0$, or (iii) $K=0\neq H$ or (iv) $K=H=0$, respectively, at the point $p$.
The terms hyperbolic, elliptic, parabolic, with reference to the sign of the
curvature, is due to Klein (1871).

Let $(x,y)$ be the coordinate system of the tangent plane $T_{p}$ relative to
the principal frame $(\mathbf{t}_{1},\mathbf{t}_{2})$. The \emph{Dupin
indicatrix }at $p$ is the following curve
\begin{equation}
\mathfrak{D}_{p}:\kappa_{1}x^{2}+\kappa_{2}y^{2}=\pm1. \label{Dupin}%
\end{equation}
It consists of either two hyperbolas, an ellipse, two parabolas, two parallel
lines, or is empty, according to whether $p$ is hyperbolic, elliptic,
parabolic, or planar, respectively. The indicatrix describes the local
geometry around $p$, as follows. The two planes parallel to $T_{p}$ and at a
small distance $\delta$ cuts the surface in a set which projects orthogonally
to a set $C_{\delta}$ in $T_{p}.$ By scaling $C_{\delta}$ with the facor
$(2\delta)^{-1/2}$, the limiting set as $\delta\rightarrow0 $ is the curve
$\mathfrak{D}_{p}$. (For a proof, see Spivak [1975], Vol. 3, p. 68.) At a
hyperbolic point $p$ the two hyperbolas have asymptotic lines given by
\[
y=\pm(\sqrt{-\kappa_{1}/\kappa_{2}})x.
\]
Therefore, the vectors $\mathbf{t}_{\theta}=\cos\theta$ $\mathbf{t}_{1}%
+\sin\theta\mathbf{\ t}_{2},$ with $\tan\theta=\pm\sqrt{-\kappa_{1}/\kappa
_{2}}$, represent the directions of the asymptotes, and they yield vanishing
normal curvature%
\begin{equation}
\mathrm{II}_{p}(\mathbf{t}_{\theta},\mathbf{t}_{\theta})=\kappa_{1}\cos
^{2}\theta+\kappa_{2}\sin^{2}\theta=\kappa_{\theta}=0. \label{second3}%
\end{equation}
In the literature these directions are called the \emph{asymptotic directions}
at the point $p$, and we note that the directions are perpendicular if
$\kappa_{1}=-\kappa_{2}$, that is, when $H=0$. On the other hand, at a
parabolic point there is only one asymptotic direction, namely the principal
direction $\mathbf{t}_{i}$ corresponding to $\kappa_{i}=0$. Let us also say
that all directions are both principal and asymptotic at a planar point.

The indicatrix is useful in the\ geometric study of two of the most
interesting families of curves on surfaces, namely the\emph{\ \ lines of
curvature }and the \emph{asymptotic lines}.\emph{\ }These are the curves which
at each point are heading in a principal or asymptotic direction,
respectively. Euler (1760) was the first who investigated the lines of
curvature, and his work inspired Monge to develop his general theory of
curvature which he applied in 1795 to the central quadrics
\[
\lambda_{1}x^{2}+\lambda_{2}y^{2}+\lambda_{3}z^{2}=c.
\]
However, asymptotic lines were introduced by Dupin. In 1813 he published his
\emph{Developpements de geometrie}, with many contributions to differential
geometry such as the idea of asymptotic lines. The indicatrix was not invented
by him, but he showed how to make more effective use of this suggestive conic.
For example, a pair of conjugate diameters of this conic are referred to as
\emph{conjugate tangents} in Dupin's theory.

The classical approach to surfaces is, of course, simplified by the usage of
vector calculus. So, let $S$ be a given surface in Euclidean 3-space
parametrized by coordinates $u,v$, assumed (for simplicity) valid on all $S$.
Thus, the parametrization
\begin{equation}
\Phi:(u,v)\rightarrow\Phi(u,v)=(x(u,v),y(u,v),z(u,v)) \label{para}%
\end{equation}
is a one-to-one smooth map from a region $D$ in the $uv$-plane onto $S$, and
the coordinate vectors $\frac{\partial\Phi}{\partial u},\frac{\partial\Phi
}{\partial v}$ at the point $p$ span the tangent plane of $S$ at $p.$ Let us
define its positive orientation and associated normal field $\mathbf{\eta}$ on
$S$, in other words the Gauss map (\ref{Gauss}), by taking \emph{\ }
\begin{equation}
\ \mathbf{\eta}=\left\vert \frac{\partial\Phi}{\partial u}\times\frac
{\partial\Phi}{\partial v}\right\vert ^{-1}\left(  \frac{\partial\Phi
}{\partial u}\times\frac{\partial\Phi}{\partial v}\right)  \label{normal}%
\end{equation}

By\ expressing $\mathrm{d}x,\mathrm{d}y,\mathrm{d}z$ in terms of $\mathrm{d}u
$ and $\mathrm{d}v$ and substituting into (\ref{Emetric}), we obtain the
squared line element restricted to the surface, namely the quadratic form
\begin{equation}
\mathrm{d}s^{2}|_{S}=E\mathrm{d}u^{2}+2F\mathrm{d}u\mathrm{d}v+G\mathrm{d}%
v^{2} \label{metric1}%
\end{equation}
whose coefficients are the following inner products%
\begin{equation}
E=\left\vert \frac{\partial\Phi}{\partial u}\right\vert ^{2},\text{ \ }%
F=\frac{\partial\Phi}{\partial u}\cdot\frac{\partial\Phi}{\partial v},\text{
\ }G=\left\vert \frac{\partial\Phi}{\partial v}\right\vert ^{2} \label{E-F-G}%
\end{equation}
Euler considered, for example, the case of a surface which can be developed on
the plane with rectangular coordinates $u,v$. Then a small triangle on $S$
must be mapped to an isometric triangle on the plane, from which he deduced,
in effect, that the coordinate vectors $\frac{\partial\Phi}{\partial u}%
,\frac{\partial\Phi}{\partial v}$ are orthonormal, so the quadratic form
(\ref{metric1}) becomes the following simple form
\begin{equation}
\mathrm{d}s^{2}=\mathrm{d}u^{2}+\mathrm{d}v^{2} \label{plane}%
\end{equation}
which characterizes the geometry of the Euclidean plane.

The surface theory of Gauss involves in fact two quadratic forms in
$\mathrm{d}u,\mathrm{d}v$, referred to as the \emph{first} and \emph{second
fundamental form,}
\begin{align}
\mathrm{I}  &  =E\mathrm{d}u^{2}+2F\mathrm{d}u\mathrm{d}v+G\mathrm{d}%
v^{2}\label{first}\\
\mathrm{II}  &  =L\mathrm{d}u^{2}+2M\mathrm{d}u\mathrm{d}v+N\mathrm{d}v^{2}
\label{second}%
\end{align}
where \textrm{I} is the metric form (\ref{metric1}), and the coefficients of
$\mathrm{II}$ are the normal components of the second order derivatives of the
map (\ref{para}),
\begin{equation}
L=\mathbf{\eta}\cdot\frac{\partial^{2}\Phi}{\partial u^{2}}\text{,
\ \ }M=\text{\ }\mathbf{\eta}\cdot\frac{\partial^{2}\Phi}{\partial u\partial
v}\text{, \ \ }N=\mathbf{\eta}\cdot\frac{\partial^{2}\Phi}{\partial v^{2}%
},\ \text{\ \ }\ \text{ \ } \label{LMN}%
\end{equation}
Indeed, the extrinsic geometry of the surface is encoded in the second form
(\ref{second}). This approach is motivated by the study of\ parametrized
curves $t\rightarrow\gamma(t)=\Phi(u(t),v(t))$ on the surface, whose
acceleration is the vector%
\begin{equation}
\gamma^{\prime\prime}=\frac{\mathrm{d}^{2}\Phi}{\mathrm{d}t^{2}}=u^{\prime
2}\frac{\partial^{2}\Phi}{\partial u^{2}}+2u^{\prime}v^{\prime}\frac
{\partial^{2}\Phi}{\partial u\partial v}+v^{\prime2}\frac{\partial^{2}\Phi
}{\partial v^{2}}+\left(  u^{\prime\prime}\frac{\partial\Phi}{\partial
u}+v^{\prime\prime}\frac{\partial\Phi}{\partial v}\right)  \label{accel}%
\end{equation}
and consequently its normal component times $\mathrm{d}t^{2}$, namely the
expression $(\gamma^{\prime\prime}\cdot\mathbf{\eta})\mathrm{d}t^{2}$, is just
the above quadratic form $\mathrm{II}$. In particular, using the natural
parameter $t=s$, the Frenet-Serret formulas (\ref{F-S}) yield%
\begin{equation}
\mathrm{II}=(\gamma^{\prime\prime}\cdot\mathbf{\eta)}\mathrm{d}s^{2}%
=(\kappa\cos\varphi)\mathrm{d}s^{2} \label{pi}%
\end{equation}
where $\kappa$ is the curvature of the curve $\gamma$ and $\varphi$ is the
angle between its principal normal $\mathbf{n}$ and the surface normal
$\mathbf{\eta.}$

Finally, by combining (\ref{metric1}) with (\ref{pi}) there is the general
formula%
\begin{equation}
\kappa\cos\varphi=\frac{Lu^{\prime2}+2Mu^{\prime}v^{\prime}+Nv^{\prime2}%
}{Eu^{\prime2}+2Fu^{\prime}v^{\prime}+Gv^{\prime2}}=\mathrm{II}_{p}%
(\mathbf{t}_{\theta},\mathbf{t}_{\theta}) \label{curv3}%
\end{equation}
where $\mathbf{t}_{\theta}$ is the unit tangent of the curve at $p$. In
particular, this explains why the form (\ref{second1}) is symmetric. For the
curve cut out by the normal section $P_{\theta}$ we have $\kappa
=\kappa_{\theta}$ and $\cos\varphi=\pm1$, so the above formula (\ref{curv3})
yields, in fact, that Euler's \textquotedblleft normal\textquotedblright%
\ curvature $\kappa_{\theta}$ equals $\mathrm{II}_{p}(\mathbf{t}_{\theta
},\mathbf{t}_{\theta})$. Consequently, we also find that Meusnier's formula
(\ref{Meus}) is just the special case of (\ref{curv3}) for curves which are planar.

Recall from (\ref{H-K}) the definition of the curvature $K$ as the product of
the principal curvatures
\[
\kappa_{1}=\mathrm{II}_{p}(\mathbf{t}_{1},\mathbf{t}_{1}),\text{ \ }\kappa
_{2}=\mathrm{II}_{p}(\mathbf{t}_{2},\mathbf{t}_{2})
\]
which are the minimum and maximum values of $\mathrm{II}_{p}(\mathbf{t}%
,\mathbf{t})$ as $\mathbf{t}$ runs over all tangent vectors of unit length,
that is, $1=$ $\left\vert \mathbf{t}\right\vert ^{2}=\mathrm{I}_{p}%
(\mathbf{t},\mathbf{t})$. The metric form $\mathrm{I}_{p}$ has matrix
coefficients $E,F,G$ with respect to the basis\ $\left\{  \frac{\partial\Phi
}{\partial u},\frac{\partial\Phi}{\partial u}\right\}  $, so let us apply the
Gram--Schmidt algorithm and replace it by an orthonormal basis and calculate
the associated symmetric matrix $\Omega$ of the form \textrm{II}$_{p}$. Then
it will be clear that the eigenvalues of $\Omega$ are just the numbers
$\kappa_{i}$, and a simple calculation yields the following useful fomulas for
the Gaussian and the mean curvature
\begin{equation}
K=\kappa_{1}\kappa_{2}=\frac{LN-M^{2}}{EG-F^{2}}\text{, \ \ \ \ \ }H=\frac
{1}{2}(\kappa_{1}+\kappa_{2})=\frac{1}{2}\frac{EN-2FM+GL}{EG-F^{2}}.
\label{HK}%
\end{equation}

For later usage, let us calculate these quantities for the graph of a function
$z=f(x,y)$, the most familiar case in elementary vector calculus. The function
in (\ref{para}) becomes $\Phi(x,y)=(x,y,f(x,y))$, and writing $f_{x}%
,f_{y},f_{xx}$ etc. for the partial derivatives of $f$, the expressions in
(\ref{HK}) read
\begin{equation}
K=\frac{f_{xx}f_{yy}-f_{xy}^{2}}{(1+f_{x}^{2}+f_{y}^{2})^{2}}\text{,
\ \ \ \ }H=\frac{(1+f_{x}^{2})f_{yy}+(1+f_{y}^{2})f_{xx}-2f_{x}f_{y}f_{xy}%
}{(1+f_{x}^{2}+f_{y}^{2})^{3/2}}. \label{HK1}%
\end{equation}
Meusnier derived this formula for $H$ in 1776. To prove his \textquotedblleft
Theorema Egregium\textquotedblright\ Gauss expressed $K$ in terms of the three
functions $E,F,$and $G$ (\ref{E-F-G}). But this also implies there is a
general formula for $LN-M^{2},$ purely in terms of $E,F,$and $G$ without
invoking any parametrization function. Liouville and Brioschi also derived
such explicit formulas, see below.

Gauss himself introduced the kind of coordinates $(u,v)$ called
\emph{isotermal} coordinates, a term due to Lam\'{e} (1833), so that the
metric $\mathrm{d}s^{2}$ and its Gaussian curvature takes the simple form
\begin{equation}
\mathrm{d}s^{2}=e^{\varphi}(\mathrm{d}u^{2}+\mathrm{d}v^{2})\text{,
\ \ \ \ }K=-\frac{1}{2}e^{-\varphi}\left(  \frac{\partial^{2}\varphi}{\partial
u^{2}}+\frac{\partial^{2}\varphi}{\partial v^{2}}\right)  . \label{isoterm}%
\end{equation}
He proved only the existence of isotermal coordinates when the given functions
$E,F,$and $G$ are analytic, wheras in the differentiable case their existence
was not proved until the 20th century. In the special case of constant
curvature, the expression (\ref{isoterm}) for $K$ is often referred to as
Liouville's equation. Liouville is well-known, for example, for his studies of
conformal geometry, where metrics differing only by a function multiple, say
$\mathrm{d}s^{2}=e^{\varphi}\mathrm{d}s^{\prime2}$ , are regarded as
\textquotedblleft identical\textquotedblright\ and are said to be
\emph{conformally} equivalent.

\subsubsection{French and Italian response to the Gaussian approach}

Actually, the geometric works of Gauss, which culminated with his
\emph{Disquisitiones} (1828), were largely unnoticed or neglected for many
years to come, say up to around 1850 or so. Ten years after the publication of
the Disquisitiones, E.F.A Minding (1806--1885) was the first who continued the
work of Gauss, with his special study of surfaces of constant curvature (see
below). In France, the tradition after Monge and his pupils still dominated
the preferences and geometrical way of thinking. S.D. Poisson (1784--1840)
published in 1832 a memoir on the curvature of surfaces, including a
historical review, but without mentioning Gauss (see Reich [1973]). During the
following decade, essentially only the versatile engineer and applied
mathematician G. Lam\'{e}, who in 1832 had accepted a chair of physics at the
\'{E}cole Polytechnique, wrote papers on differential geometry. Among his
interests we find triply orthogonal systems of sufaces, and through his
studies of heat conduction he was also led to a general theory of curvilinear
coordinates. Up to the 1860's, the \textquotedblleft best\textquotedblright%
\ non-trivial example of a triply orthogonal system was the family of confocal
quadrics, and both Lam\'{e} [1839] and Jacobi [1839] introduced elliptic (or
ellipsoidal) coordinates using this system. Lam\'{e} used them to separate and
solve the Laplace equation $\Delta(f)=0$, whereas Jacobi calculated the
geodesics on the ellipsoid.

But in the early 1840's some younger French geometers, such as Bertrand and
Bonnet in their early twenties, made reference to Gauss and his Disquisitiones
in their first papers, and at the end of the decade an ebullient interest in
the geometric works of Gauss 20 years earlier burst out in France and Italy.
In particular, as Monge's celebrated paper [1807] was republished by Liouville
in 1850, he also gave a new proof of the Theorema Egregium which appeared in
the appendix, together with the original latin version of the Disquisitiones.

Liouville had an explicit formula for $K$ in 1851 which, in fact, the Italian
Beltrami made clever use of in his papers in the mid 1860's. In Italy, F.
Brioschi (1824--97) had greatly influenced the direction of higher education
and research in mathematics, and among his doctoral students we find Cremona
(1853) and Beltrami (1856). In 1863 he founded the Technical University in
Milan, where he served both as a director and professor of mathematics and
hydraulics. Maybe Beltrami's teacher had not yet worked out his formula by
1865, but it is Brioschi's formula for $K$ rather than Liouville's \ which is
most easily found in textbooks and online pages on differential geometry. In
orthogonal coordinates $u,v$, namely when $F=0$, the formula simplifies to
$\ $%
\begin{equation}
K=-\frac{1}{2\sqrt{EG}}\left(  \frac{\partial}{\partial u}(\frac{G_{u}}%
{\sqrt{EG}})+\frac{\partial}{\partial v}(\frac{E_{v}}{\sqrt{EG}})\right)  .
\label{K0}%
\end{equation}

Next, let us also recall some geometric invariants naturally arising in the
study of curves $\gamma(s)$ on a given oriented surface $S$, using their
arc-length $s$ as the natural parameter. Along the curve there is the velocity
vector $\gamma^{\prime}(s)=\mathbf{t}$, the acceleration $\gamma^{\prime
\prime}(s)=\kappa\mathbf{n}$, and the resulting Frenet-Serret frame
$(\mathbf{t,n,b})$, see (\ref{F-S}). But there is also the orthonormal
\emph{Darboux frame} $(\mathbf{t,u,\eta})$, with $\mathbf{u}$ in the tangent
plane, chosen so that $\mathbf{t\times u=}$ $\mathbf{\eta}$ is the positively
directed normal, as before. Again, in analogy with (\ref{F-S}),
differentiation with respect to $s$ yields a dynamical system with a skew
symmetric matrix, namely for suitable coefficients $\kappa_{g},\kappa_{\eta
},\tau_{g}$ we can write%

\begin{equation}
\mathbf{t}^{\prime}=\kappa_{g}\mathbf{u+}\kappa_{\eta}\mathbf{\eta}\text{,
\ }\mathbf{u}^{\prime}=-\text{\ }\kappa_{g}\mathbf{t}+\tau_{g}\mathbf{\eta
}\text{, \ }\mathbf{\eta}^{\prime}=-\kappa_{\eta}\mathbf{t-}\tau
_{g}\mathbf{u.}\text{\ } \label{Darboux}%
\end{equation}
The leftmost equation expresses the decomposition of the acceleration
(\ref{accel}) into its tangential and normal component%
\[
\gamma^{\prime\prime}=\mathbf{t}^{\prime}=\gamma_{g}^{\prime\prime}%
+\gamma_{\eta}^{\prime\prime}=\kappa_{g}\mathbf{u+}\kappa_{\eta}\mathbf{\eta}%
\]
called the \emph{geodesic} and \emph{normal} \emph{curvature} vector,
respectively. It follows that $\kappa_{\eta}$ is just the normal curvature
introduced by Euler, see (\ref{curv2}), namely $\kappa_{\eta}=\kappa_{\theta}$
where the angle $\theta$ gives the direction of $\mathbf{t}$, so that
\begin{equation}
\kappa^{2}=\kappa_{g}^{2}+\kappa_{\eta}^{2},\ \ \ \mathrm{II}_{\gamma
(s)}(\mathbf{t,t)=}\ \kappa_{\eta}. \label{curv}%
\end{equation}
Gauss referred to $\kappa_{g}$ as the \textquotedblleft
Seitenkr\"{u}mmung\textquotedblright, and the term tangential curvature was
also used. However, the term used today is \emph{geodesic curvature}, which
dates back to Bonnet (1848).

\subsubsection{Curves on a surface and the Gauss-Bonnet theorem}

The terms introduced above serve to characterize the following three classical
main types of curves $\gamma(s)$ confined to a surface $S$, as follows:

\begin{itemize}
\item \emph{Geodesic curves:} Their geodesic curvature $\kappa_{g}$ vanishes.
The curves are locally characterized by the geodesic equation $\kappa_{g}=0$,
which in local coordinates is a 2nd order nonlinear differential equation.
Therefore, a geodesic starting from a point $p_{0}$ is uniquely determined by
its initial direction $\mathbf{t}(0)$. Their crucial geometric property is the
\textquotedblleft shortest length\textquotedblright\ property, which generally
holds for small segments of the curve, but not necessarily for longer
segments. But surely, a shortest curve between two given points must be a
geodesic. Geodesics are the natural generalization of the \textquotedblleft
straight lines\textquotedblright\ in Euclidean or hyperbolic geometry. Being a
geodesic curve is an intrinsic property, that is, it remains a geodesic under
isometric deformations.

\item \emph{Asymptotic lines: } Their normal curvature $\kappa_{\eta}$
vanishes. Other equivalent conditions are (i) $\kappa_{g}=\pm\kappa$, or (ii)
the osculating plane (if defined) equals the tangent plane, or (iii) the
Frenet--Serret frame and the Darboux frame differ at most by signs, namely
$(\mathbf{t,n,b})=$ $(\mathbf{t,}\pm\mathbf{u,}\pm\mathbf{\eta}).$ Moreover,
the identity $\tau_{g}=\tau$ (if defined) also holds. The equation
$\kappa_{\eta}=0$ means the vanishing of the second fundamental form
(\ref{second}) along the curve, namely
\begin{equation}
L\mathrm{d}u^{2}+2M\mathrm{d}u\mathrm{d}v+N\mathrm{d}v^{2}=0 \label{asymp}%
\end{equation}
By regarding this relation as a quadratic equation this yields (in general)
two vector fields on the surface, whose integral curves are the asymptotic lines.

\item \emph{Lines of curvature:} Their geodesic torsion $\tau_{g}$ vanishes.
Equivalently, their velocity vector $\gamma^{\prime}(s)=\mathbf{t}$($s)$ is a
principal direction $\mathbf{t}_{i}(s)$, for all $s$. It is easy to express
the principal directions as the zero directions of the following quadratic
form
\begin{equation}
(EM-FL)\mathrm{d}u^{2}+(EN-GL)\mathrm{d}u\mathrm{d}v+(FN-GM)\mathrm{d}v^{2}=0,
\label{curvline}%
\end{equation}
and in analogy with the case (\ref{asymp}) this also yields two vector fields
on the surface, whose integral curves are the lines of curvature.
\end{itemize}

Euler (1732) was the first to work out a differential equation for the
shortest curve between two points on a surface which solves equation
(\ref{diff2}) (see Rosenfeld[1988: 282]). Dupin's theorem, which we shall
return to below, provides in specific cases a geometric construction of the
lines of curvature curves. Obvious examples of asymptotic curves, which are
geodesics as well, are straight lines lying on the surface. Their unit speed
motion along the line has a constant velocity vector $\mathbf{t,}$ so
$\mathbf{t}^{\prime}=0$ and consequently $\kappa_{g}=$ $\kappa_{\eta}=0$ by
(\ref{Darboux}). But $\kappa_{\eta}$ lies between $\kappa_{1}$ and $\kappa
_{2}$, so $K=\kappa_{1}\kappa_{2}\leq0$ must hold along the line. Moreover,
strict equality $K=0$ holds if and only if our line is also a line of
curvature, or equivalently, the tangent plane as the same normal direction and
hence is constant along the line. The reader may verify these statements, say
using (\ref{Darboux}) and (\ref{second1}).

In his differential geometric studies Monge had found envelopes of families of
surfaces. For a one-parameter family $F(x,y,z,a)=0$ enveloping a surface, the
procedure is similar to that of curves, see the equations (\ref{envelope}).
For each value of $a$ the equations yield a curve, called the
\emph{characteristic}, on the surface labelled by $a$. Then by varying $a$ the
characteristics sweep out the enveloping surface. For example, the Dupin
cyclides are surfaces which can be enveloped by a one-parameter family of
spheres, in fact, in two different ways.

A 2-parametric system of lines is usually called a ray systems or a line
congruence. They play an important role in line geometry, where their envelope
is referred to as the \emph{focal set}, or the focal surface
(Brennfl\"{a}che), which may degenerate to a curve called the evolute. A
typical focal surface $F$ has two components $F_{1}$ and $F_{2}$, and the ray
system consists of their common tangent lines. This also applies to the study
of surfaces, by considering the ray system of lines normal to a given surface
$S$.

On the other hand, using the notion of curvature there is another construction
of the focal sets $F_{i}$ of a given surface. Namely, by pointwise pushing the
surface in the normal direction a distance equal to its curvature radii, one
obtains the sets $F_{i}$ as the locus of points
\[
p+\frac{1}{\kappa_{i}(p)}\mathbf{\eta}(p),\text{ \ }p\in S\text{, \ }i=1,2.
\]
Letting $S$ be a torus obtained by rotating a circle, for example, the two
sets will be a circle and a straight line (axis if symmetry).

In his paper [1828] Gauss considered geodesic triangles $\Delta$ on a surface
$S$, and denoting the angles by $\alpha,\beta,$ and $\gamma$ he derived the
simple formula%
\begin{equation}%
{\displaystyle\iint\limits_{\Delta}}
K\mathrm{d}S=\alpha+\beta+\gamma-\pi\label{excess}%
\end{equation}
which expresses the \textquotedblleft angular excess\textquotedblright\ as the
total curvature of the triangle. Twenty years later, Bonnet 1848 considered
the general case of a triangle with piecewise smooth boundary edges $C_{i}$,
and he extended the formula by adding correction terms on the left side,
namely the line integral of the geodesic curvature $\kappa_{g}$ along the
three edges. This is the essense of the general Gauss--Bonnet formula
(\ref{GaussBonnet}).

The remaining part amounts to some combinatorial bookkeeping which arises when
$n$ triangles fill up a region $S_{n}$ of an oriented surface $S,$ and we
apply and add together the modified formula (\ref{excess}) for each triangle.
It is assumed that any two triangles are either disjoint or have just one
common vertex or edge. Then the number of vertices ($V)$, edges ($E$) and
triangles ($T$) are related by the Euler characteristic $\mathcal{X}%
(S_{n})=V-E+F$ of the region, and this is the only global topological
invariant needed here. Clearly, for a closed disk and a sphere the number is
$1$ and $2$, respectively. The boundary $\partial S_{n}$ of the region
consists of those edges belonging to a single triangle.

Next, let us make use of the induced orientation of the triangles and their
edges, noting that adjacent triangles induce opposite directions on their
common edge. Therefore, the line integrals away from the boundary cancel each
other. Now, let us replace the angles of a triangle by its oriented outer
angles $\alpha^{\prime}=\pi-\alpha,\beta^{\prime}=\pi-\beta$ etc., which
measure the \textquotedblleft jumps\textquotedblright\ of the tangential
direction at the vertices. In particular, at each corner of the boundary there
is a net angular \textquotedblleft jump\textquotedblright, and let us move
them to the left side of our equation. Then the remaining \textquotedblleft
jumps\textquotedblright\ on the right side, in fact, add up to a multiple of
$2\pi$, namely%
\begin{equation}%
{\displaystyle\iint\limits_{S_{n}}}
KdS+%
{\displaystyle\int\limits_{\partial S_{n}}}
\kappa_{g}ds+%
{\displaystyle\sum\limits_{\partial S_{n}}}
(\text{jumps)}=2\pi\mathcal{X}(S_{n}). \label{GaussBonnet}%
\end{equation}
This is the Gauss--Bonnet formula, valid for a compact oriented surface $S$
($=S_{n})$ with a piecewise smooth boundary, where the "jumps" disappear if
the boundary is smooth. It is very likely that Gauss had ideas about how to
generalize his formula (\ref{excess}), but he did not publish more on the
topic. Bonnet proved special cases of equation (\ref{GaussBonnet}) in 1848.

\subsection{Surfaces in classical analytic geometry}

\subsubsection{Ruled and developable surfaces}

Let us have a closer look at some types of surfaces frequently encountered in
the classical literature, starting with the ruled surfaces. A surface is
\emph{ruled }if for each point there is a straight line on the surface passing
through the point. Obvious examples are cylinders and cones. It is doubly
ruled if there are two distinct lines passing through each point, and familiar
examples are the hyperbolic paraboloid and the hyperboloid of one sheet:%
\begin{equation}
z=\frac{x^{2}}{a^{2}}-\frac{y^{2}}{b^{2}}\text{, \ \ \ \ \ }\frac{x^{2}}%
{a^{2}}+\frac{y^{2}}{b^{2}}-\frac{z^{2}}{c^{2}}=1. \label{hyperboloid}%
\end{equation}
In fact, these are the only doubly ruled surfaces of degree 2, and the only
triply ruled surface is the plane.

A \emph{ruling }of the surface is a 1-parameter family of straight lines on
the surface, which yields a parametrization of the surface of type
\begin{equation}
(s,t)\rightarrow\gamma(s)+t\delta(s) \label{ruled}%
\end{equation}
where $\gamma(s)\ $and $\delta(s)$ are given space curves. Each of the lines,
parametrized by $t$, is called a \emph{generatix}, whereas the curve
$\gamma(s)$ cutting every line is the \emph{directrix}, \emph{\ }and
$\delta(s)$ is the \emph{director}. Monge and his school gave these surfaces
some thought, but later they played an important role in Pl\"{u}cker's study
of line geometry. Recall that any surface must have curvature $K\leq0$ along a
straight line, consequently a ruled surface satisfies $K\leq0$ everywhere.

Monge's work was sometimes\ preceded by Euler's, but he worked independently
and with his own originality. For example, developable\emph{\ }surfaces were
introduced by Euler and Monge in 1772 and 1785, respectively, using their own
definition. To begin with, let us define a \emph{developable\ }surface in
3-space to be a ruled surface with vanishing Gaussian curvature, $K=0$. Note
that the above hyperboloid (\ref{hyperboloid}) has negative curvature, so it
is an example of a ruled surface which is not developable. For the simpler and
final definition of \textquotedblleft developable\textquotedblright, see
(\ref{equivalences}) below.

Ruled and developable surfaces appear naturally in surface theory, as in the
following construction. For any chosen curve $\gamma$ on the surface $S$, the
family of normal lines to $S$ along $\gamma$ spans a ruled surface $\bar{S}$
with $\gamma$ as its directrix. Then, by a theorem of Bonnet, the curve
$\gamma$ is a line of curvature of $S$ if and only if the associated ruled
surface $\bar{S}$ is developable.

Classically, the developable surfaces were found to belong to four types:%
\begin{equation}
(i)\text{ planes, \ }(ii)\text{ cylinders, \ }(iii)\text{ cones, }(iv)\text{
tangent developables.} \label{list}%
\end{equation}
However, this list is not complete, not even among the analytic developable
ones. As an example, there is an analytic developable, homeomorphic to the
M\"{o}bius strip (see Spivak[1975], vol. 3, p. 355). Ruled surfaces different
from the types (\ref{list}) became also known as \emph{scrolls} (Schraubenfl\"{a}che).

Let us give parametrizations (\ref{ruled}) of the developable surfaces of type
(ii) - (iv) in (\ref{list}). They are actually generated by a single space
curve $\gamma(s)$ with velocity $\gamma^{\prime}(s)\neq0$, on some interval
$s\in(a,b)$, in the following way:

\begin{itemize}
\item Generalized cylinder: Take $\gamma(s)$\ to be a curve in the xy-plane
and let $\delta$ be the constant vector $(0,0,1)$. The surface (\ref{ruled})
has principal curvature $\kappa_{1}=0$ in the vertical direction, and
$\kappa_{2}=\kappa(s)$ is the curvature of $\gamma(s)$, with principal
direction along $\gamma^{\prime}(s)$. Consequently $K=0$ and $H=$ \ $\frac
{1}{2}\kappa(s).$

\item Generalized cone: Choose a vector $v_{0}$ as the vertex of the cone, and
let $\delta(s)=\gamma(s)-v_{0}$. The vertex $v_{0}$, corresponding to $t=-1$
in (\ref{ruled}), is a singular point of the surface. Away from the vertex,
$\gamma(s)-v_{0}$ points in the principal direction \ with $\kappa_{1}=0$, so
$K=\kappa_{1}\kappa_{2}=0.$ \ \ \ 

\item Tangent developable: Assume $\gamma(s)$ is an arc-length
parametrization, $\gamma(s)$ has curvature $\kappa(s)\neq0$, and take
$\delta(s)=\gamma^{\prime}(s)$. The surface (\ref{ruled}) has two sheets which
for $t=0$ meet along the curve $\gamma(s)$ as a cuspical edge. Away from this
curve, the surface has the principal curvature $\kappa_{1}=0$ in the direction
of $\gamma^{\prime}(s)$, so again the total curvature vanishes.
\end{itemize}

Already in 1825 Gauss published a paper on conformal transformations, where he
compared two surfaces whose metric expressions $\mathrm{d}s^{2}$ (see
(\ref{metric1})) differ by a scalar function, say%
\begin{equation}
E\mathrm{d}u^{2}+2F\mathrm{d}u\mathrm{d}v+G\mathrm{d}v^{2}=\mu^{2}(E^{\prime
}\mathrm{d}u^{2}+2F^{\prime}\mathrm{d}u\mathrm{d}v+G^{\prime}\mathrm{d}v^{2}).
\label{metric3}%
\end{equation}
He observed that for $\mu=1$ a complete equality (Gleichheit) holds between
the surfaces, so that one surface can be developed onto the other. Thus he
extended the notion of \textquotedblleft developable\textquotedblright, which
in this context also became known as \emph{applicable}, and thus a major
classical problem has been to determine those surfaces applicable to a given
surface, not just the plane.

Briefly, applicable surfaces have the same metric expression $\mathrm{d}s^{2}
$ (see (\ref{metric1})), and they were regarded as deformations of each other
(locally). In modern terminology, they are (locally) isometric, but not
necessarily \emph{congruent,} that is, transformable to each other by a rigid
motion of the ambient space. For example, the catenoid and the helicoid are
applicable, one can in fact \textquotedblleft stretch and
twist\textquotedblright\ the catenoid continuously to the form of a helicoid
without changing the intrinsic geometry. They are examples of minimal surfaces
discussed below, and certainly they look very different.

Following the footsteps of Gauss, Minding's theorem (1839) says that surfaces
in 3-space with the same constant curvature $K$ are mutually applicable, that
is, locally isometric. In particular, surfaces with vanishing curvature $K=0$
are applicable with the plane and hence they can be constructed by suitably
bending of plane regions. Therefore they are also ruled, namely we have the
following three equivalent conditions for surfaces, briefly said to be
\emph{flat}:%

\begin{equation}
\text{(i) developable, \ (ii) }K=0\text{, \ (iii) locally isometric to a
plane.\ } \label{equivalences}%
\end{equation}

In France, Bonnet and Bour investigated the curvature of special surfaces such
as ruled surfaces or surfaces of rotation, as well as the developability for
ruled surfaces. For example, a theorem of Bour says that a scroll is
applicable to some surface of rotation. But first of all, when it comes to the
classical problem concerning applicable surfaces, Weingarten in Germany was
the first who made a major step forward when he described in 1863 a class of
surfaces applicable to a given surface of rotation.

At this point, let us illustrate with some examples why the study of
asymptotic curves played such an important role in the classical surface
theory. According to a theorem of Bonnet, an isometry between non-ruled
surfaces which maps a family of asymptotic curves into the asymptotic curves
of the other surface must be a rigid motion (cf. Chern [1991]).The problem of
finding another surface applicable to a given surface $S$ amounts analytically
to solving a certain Monge--Amp\`{e}re equation, whose characteristics are the
asymptotic curves of $S$. This system of equations may well become
over-determined if additional constraints are posed on the surface. There is
the following result concerning a given space curve $\gamma$ lying on $S$,
namely if $\gamma$ is asymptotic then one can construct infinitely many
surfaces $S^{\prime}$ through $\gamma$ which are applicable to $S$. On the
other hand, if $\gamma$ is not asymptotic and is kept fixed, then there no
such surface except $S$ itsef.

\subsubsection{Minimal surfaces}

Next, let us turn to minimal surfaces, whose well known physical models are
the soap bubbles suspended between strings and always tend to minimize their
area. Ever since the first discoveries in the 18th century the determination
of minimal surfaces has posed a challenge to geometers up to present time. It
started with Euler's discovery in 1744 that a catenary (\ref{catenary})
rotated around the $x$-axis yields a minimal surface of revolution. This is
the well known \emph{catenoid}, which for any finite segment of the catenary
has the smallest area among all surfaces suspended between the two boundary
circles. By a simple variational analysis on the area of a surface of
rotation, Euler derived a differential equation for the rotated curve having
the catenary as its solution. He may not have proved it rigorously, but with
the original definition the catenoid is, in fact, the only minimal surface of
rotation different from the plane. This was also verified by Bonnet in the 1850's.

Using a more general variational analysis of area, Lagrange (1860) studied
minmal surfaces on the explicit form $z=f(x,y)$. He considered the area
\[
A(t)=%
{\displaystyle\int\limits_{U}}
(1+\tilde{f}_{x}^{2}+\tilde{f}_{y}^{2})^{2}\mathrm{d}x\mathrm{d}y
\]
of a variation $\tilde{f}(x,y;t)\ $of a minimal surface $(t=0)$ over an open
region $U$ in the $xy$-plane, with boundary fixed for all $t$. By demanding
$A^{\prime}(0)=0$ for all variations he derived the associated Euler-Lagrange
equation
\begin{equation}
(1+f_{y}^{2})f_{xx}+(1+f_{x}^{2})f_{yy}-2f_{x}f_{y}f_{xy}=0. \label{minimal}%
\end{equation}
\qquad It is a rewarding exercise at an undergraduate level to deduce Euler's
catenary from this equation by expressing it in polar coordinates $(r,\theta)$
and setting $\theta=0$. The solutions of the resulting reduced differential
equation, $rf^{\prime\prime}+f^{\prime}(1+f^{\prime2})=0$, are those curves in
the $(r,z)$-plane which generate the rotationally symmetric minimal surfaces
$z=f(r)$. The reduced equation is easily solved by quadrature and yields the
catenary family of curves, $a\cosh(\frac{z}{a}+b)-r=0$, as expected.

This simple calculation applied to an example from the early 18th century,
indeed, illustrates the basic idea of \textquotedblleft
equivariant\textquotedblright\ differential geometry, a modern reduction
technique based on the interaction between symmetry and the least action
principle. This lies at the heart of Sophus Lie's symmetry approach to
differential equations, but for many reasons the method was not properly
developed until the late 20th century.

Meusnier discovered in 1776 the\emph{\ helicoid}, which is closely related to
the catenoid. In sylindrical coordinates it has the simple equation
$z=c\theta$. Until about 1830 these two surfaces besides the plane were the
only known minimal surfaces, but at this time H.F. Scherk (1798--1885)
discovered his first surface, implicitly defined by $e^{az}\cos ay-\cos ax=0$.
His discovery was regarded as sensational, and in 1831 he was awarded a prize
at the Jablonowski Society in Leipzig. The reader can find the Scherk surface
by seeking the solutions of equation (\ref{minimal}) of the splitting type
$z=g(x)+h(y)$, by applying the method of separation of variables. However, the
helicoid is the only ruled minimal surface in 3-space other than the plane, as
was shown by E.C. Catalan (1814--1894) in 1842.

On the other hand, Meusnier derived in 1776 also his formula (\ref{HK1}) for
the mean curvature $H$, and by comparison with formula (\ref{minimal}) he thus
observed that $H=0$ is a necessary condition for minimality. In fact, in the
19th century the simple and local condition $H=0$ became the new and modern
definition of a minimal surface. In effect, the original meaning of a minimal
surface, being infinitely extendible and with no boundary curve, was abandoned.

Henceforth minimal surfaces became the natural surface analogue of geodesic
curves, which are locally of shortest length, but perhaps not the shortest
curve between any two of its points. The mathematical problem of existence of
a minimal surface with a given boundary became known as the Plateau problem,
after the Belgian physicist Joseph Plateau (1801--1883), who conducted
extensive studies of soap films. We mention briefly that the existence of a
solution, for a given boundary curve, was not proved until 1931, but little
could be said about the geometric properties of the solution.

During the 19th century complex analysis gradually became an important
direction of mathematics. Then it also turned out that complex analytic
functions have a close connection with minimal surfaces, and in the 1860's
Weierstrass, Riemann and Enneper found representation formulas which
parametrize a minimal surfaces for each pair $(f,g$) of functions. Thus
Enneper and Weierstrass created a whole class of new parametrizations, roughly
by taking two functions on a domain $D$ say, with $fg^{2}$ holomorphic, and
define the surface as the set of points%

\begin{equation}
(x,y,z)=\operatorname{Re}\left(
{\displaystyle\int}
f(1-g^{2})\mathrm{d}\zeta,%
{\displaystyle\int}
if(1+g^{2})\mathrm{d}\zeta,%
{\displaystyle\int}
2fg\mathrm{d}\zeta\right)  . \label{minimal1}%
\end{equation}
The difficulty lies in controlling the global behavior of the surface, which
may have singularities such as self-intersections. We also mention that H.A.
Schwarz and his collaborators solved the Plateau problem in 1865 for special
boundary curves, by finding appropriate functions to be inserted into the
above formulae (\ref{minimal1}).\ As the calculus of variations and
topological methods were developed, the study of minimal surfaces took new
directions in the 20th century, and the methods of Weierstrass and Schwarz
came more in the background.

\subsubsection{Dupin's cyclides and some related topics}

Dupin's cyclides are among those surfaces whose remarkable properties
attracted considerable attention in the early 19th century and, in fact, since
the 1980's they still do. They constitute a 3-parameter family, and originally
the cyclide is defined as the envelope of the 1-parameter family of spheres
tangent to three fixed spheres. In fact, Dupin discovered them still as an
undergraduate student, but they appeared first in his \emph{Developpements }in
1813. Later, in his book "Applications de Geometrie" (1822), he called them
\emph{cyclides.} Algebraically, a cyclide is of order three or four, and
having circular lines of curvature is one of their essential properties.
Moreover, their two focal surfaces $F_{i}$ degenerate to curves of second order.

But there are still many other ways of characterizing them, alternative
definitions are due to Liouville, J.C. Maxwell (1868), Casey (1871) and Cayley
(1873). In recent years there are, in fact, many indepth studies of their
algebraic and geometric properties. For a survey of all this, including modern
applications to geometric design, we refer to Chandru et al. [1989]. On the
other hand, in the mid 1860's Dupin's cyclides became somewhat subordinate to
another family of surfaces also named cyclides, namely the generalized
cyclides\emph{\ }discovered by Darboux and Moutard. As a consequence, the
previous cyclides were rather overshadowed by the new family of surfaces,
having different properties which for the purpose of the Darboux school were
regarded more important. We shall return to them below.

An interesting aspect of Dupin's cyclides is their symmetry properties; they
can, in fact, be generated and described neatly in terms of \emph{conformal}
transformations. To explain this, first recall that these are the angle
preserving transformations, a property which is also evident from their action
on the squared line element (or metric) of the space, namely
\begin{equation}
\mathrm{d}s^{2}\longrightarrow\mu^{2}\mathrm{d}s^{2}\text{ \ \ \ \ (or
\ }(\mathrm{d}s^{2}=0)\text{\ }\longrightarrow(\mathrm{d}s^{2}=0))
\label{conformal}%
\end{equation}
where the function $\mu^{2\text{ }}>0$ may depend on the transformation. As
indicated, in Klein's letters \textquotedblleft conformal\textquotedblright%
\ is also expressed by the invariance of the equation $\mathrm{d}s^{2}=0$.

In our case, the space is the Euclidean space $\mathbb{R}^{3}$ (or a
subregion) with the metric (\ref{Emetric}), and therefore the celebrated
Liouville's theorem (1846) on conformal mappings describes them as the
composition of similarities (that is, Euclidean motions and homotheties) and
inversions. The latter type consists of the classical geometric
transformations called \emph{inversion }with respect to a sphere, and
Liouville called them transformation by \emph{reciprocal radii}. The
transformation interchanges the inside and outside of a fixed sphere and
inverts the radial distance. For a sphere of radius $\rho$ centered at $v_{0}%
$, the vector algebra expression of the inversion is \
\begin{equation}
v\rightarrow v_{0}+\rho^{2}\frac{v-v_{0}}{\left\vert v-v_{0}\right\vert ^{2}}.
\label{inversion}%
\end{equation}
Liouville showed that a cyclide can be obtained from a torus of revolution, or
a circular cylinder or a circular conic, by a suitable choice of inversion.
The cyclide is also \emph{anallagmatic}, in the sense of being its own inverse
with respect to at least one inversion.

The corresponding transformation of the plane is the inversion wih respect to
a circle, defined similarly. If straight lines are regarded as circles of
infinite radius, one can say briefly that the transformation transforms
circles into circles with one circle fixed. The method of inversion seems to
be attributed to Steiner (1824), but it was used in special cases by Poncelet
(1822), Pl\"{u}cker and others, and it was studied most extensively by
M\"{o}bius (1855). It was also discovered through physical considerations; for
example, it appeared as \textquotedblleft the method of
images\textquotedblright\ in the work on electrostatics by W. Thomson (1845).
Inversion was one of the first non-linear transformations to be deeply studied
in geometry. As a natural generalization, Cremona introduced in 1854 the
general birational transformation on the plane, which became known as Cremona
transformations and developed by him in the 1860's. They were found to have
many applications such as the reduction of singularities of curves, and the
study of elliptic integrals and Riemann surfaces.

\ \ \ \ \ 

{\large Remarks on groups}

In retrospect, it is tempting to interpret the above geometric objects in the
light of modern group theory. By adding to $\mathbb{R}^{3}$ a point at
infinity, our space becomes (via stereographic projection) conformally the
same as the 3-sphere $S^{3}$, with its transformation group $CG(3)$ consisting
of the tenfold infinity of conformal transformations. Namely, it is a
10-dimensional Lie group, containing the orthogonal group $SO(4)$ as the group
of isometries of the sphere. Now, the latter clearly contains the
2-dimensional torus group $T$, which together with all its conjugates
$T^{\prime}$ in $CG(3)$ act on the sphere with tori as their orbits. The
images of these tori in $\mathbb{R}^{3}\cup\{\infty\}$ are the Dupin cyclides,
and they can be permuted among themselves by transformations from $CG(3)$.

\ \ \ \ \ \ \ \ \ \ \ \ \ \ \ \ 

The above description of surfaces using groups has a striking similarity with
Klein and Lie's approach to W-surfaces, as part of their study of
W-configurations in 1870--71. We touch this topic only briefly and refer to
Hawkins [2000] \S 1.2, for supporting evidence that it was during the
collaboration on this project, based upon Lie's paper [1870a], that Klein and
Lie developed the basic general principles leading to the idea of continuous
groups of transformations, in analogy with the customary definition of a group
of substitutions (or permutations) in algebra. In this study they were working
with specific projective transformations on complex projective 3-space.
Namely, in the mentioned paper on tetrahedral line complexes, Lie had focused
attention on the totality $G$ of projective transformations fixing the
vertices of a given tetrahedron- in modern terms $G$ is a (complex)
3-dimensional torus. At this time Klein and Lie used terms like
\textquotedblleft cycle\textquotedblright\ or \textquotedblleft closed
system\textquotedblright\ for families of transformations closed under
composition, tacitly assuming the family would also be closed under taking
inverses as in the algebraic case. Since the surfaces they were seeking are
just the orbits of the various 2-dimensional subtori $T$ of $G$, the
classification of these surfaces would necessarily involve a certain group
classification problem. The problem was rather intractable and hence
postponed, but really they did not return to it (cf. note to letter of 12.9.71).

\subsubsection{Pseudospherical surfaces}

Surfaces of constant curvature $K$ in Euclidean 3-space are clearly natural
geometric objects, and for $K=0$, resp. $K>0$, the prototype examples are the
plane and the round sphere respectively. Minding (1839) also posed the
question about the uniqueness of the sphere among closed surfaces with
constant $K>0.$ The affirmative answer was not given until Liebmann's theorem
(1900), valid for non-singular and $C^{2}$-differentiable surfaces. However,
in the case of $K<0$ little was known until Minding made an explicit study of
such surfaces in 1838. He discovered, in fact, three types of surfaces of
rotation, and some non-rotational surfaces as well. The geometrically simplest
one is the \emph{tractoid}, namely the surface obtained by rotating the
tractrix (\ref{tractrix}), and\ it became known as the \emph{pseudosphere. }

By using the expression (\ref{HK1}) for $K$ one is led to the following
differential equation
\begin{equation}
f_{xx}f_{yy}-f_{xy}^{2}+(1+f_{x}^{2}+f_{y}^{2})^{2}=0 \label{tractoid}%
\end{equation}
valid for surfaces on the explicit form $z=f(x,y)$ and with $K=-1$. Then, an
application of the same reduction to (\ref{tractoid}) as was applied to
equation (\ref{minimal}) when we found the catenoid, will also yield a
differential equation for the tractrix which is solvable by quadrature. The
pseudosphere was, in fact, already known to Gauss, who referred to it as the
\textquotedblleft opposite\textquotedblright\ of the sphere in a note written
in the 1820's. For many years, the term pseudosphere was used confusingly in
the literature for any surface of constant negative curvature. After all,
Minding had concluded that all these surfaces are isometric, or more
precisely, applicable to each other. It was Beltrami who finally in 1868
referred to them as \emph{pseudospherical }surfaces, in order to
\textquotedblleft avoid circumlocation\textquotedblright\ as he puts it, 30
years after Minding's results had appeared. We refer to Coddington [1905] for
an interesting account of the historical development of pseudospherical
surfaces during the years 1837--1887.

We also remark that in terms of asymptotic coordinates $(u,v)$, that is, the
coordinate lines are the asymptotic lines, equation (\ref{tractoid}) takes the
original form of the sine-Gordon equation%
\begin{equation}
\psi_{uv}=\sin\psi\label{sine G}%
\end{equation}
with $\psi$ as the angle between the asymptotic curves. As a consequence, the
study of pseudospherical surfaces is euivalent to that of the above equation,
and the equation was much used for that purpose in the 19th century. In modern
mathematical theories the equation is also found to be interesting because it
has soliton solutions.

\subsubsection{Coordinate geometry, triply orthogonal systems, and generalized
cyclides}

The classical study of the differential geometry of 3-space led, in fact, to
several kinds of differential equations with modern applications nowadays. For
us, it is appropriate to recall some of the early developments on triply
orthogonal systems. Although applications of particular examples occurred
already in the works of Leibniz and Euler, most of the early contributions are
due to French geometers such as Lam\'{e}, Dupin, Liouville, Bonnet and Darboux.

A coordinate system in 3-space (or a subregion) amounts to a triple
$S^{(\alpha)}$, $S^{(\beta)}$, $S^{(\gamma)}$ of 1-parameter families of
surfaces, called coordinate surfaces, such that every point $p$ lies on a
unique surface from each family, which yields a bijective correspondence
$p\longleftrightarrow(\alpha,\beta,\gamma)$. The triple is said to be
orthogonal if the coordinate surfaces intersect each other perpendicularly.
Then, according to the celebrated \emph{Dupin's theorem}, two orthogonally
intersecting surfaces must intersect along lines of curvature. Dupin published
his proof in his previously mentioned paper \emph{D\'{e}veloppements }(1813);
for a modern proof, see Spivak[1875], Vol. 3. The orthogonality property of
the coordinates is also reflected by the corresponding expression for the
arc-length element (\ref{Emetric}), namely
\begin{equation}
\mathrm{d}s^{2}=P\mathrm{d}\alpha^{2}+O\mathrm{d}\beta^{2}+R\mathrm{d}%
\gamma^{2} \label{metric2}%
\end{equation}
where the coefficients $P,Q$, and $R$ are functions depending on the geometry
of the surfaces.

During the first half of the 19th century French geometry developed roughly in
two directions, namely with focus on differential geometry or projective
geometry, in the tradition of Monge and Poncelet, respectively. According to
Hawkins[2000: 28], the \textquotedblleft French metrical
geometry\textquotedblright\ --- a term often used by Klein and Lie ---
referred to those geometers who combined concepts from both methodologies. In
fact, their basic approach is within the framework of conformal geometry
rather than projective geometry, partly inspired by Liouville's theorem (1846)
which gives a precise description of the totality of conformal transformations
of 3-dimensional space (or higher). During the 1860's in Paris, a bright
aspiring student arose from the elite schools for mathematical training and
gradually developed his ideas which placed him centrally among the
\textquotedblleft metrical\textquotedblright\ geometers around 1870. His name
was Gaston Darboux, born the same year as Sophus Lie.

Darboux's first two papers, which appeared in the Nouvelles Annales in 1864,
are concerned with his construction of a family of 4th degree curves which he
referred to as \emph{cyclic}. He starts with the planar sections of tori and
the intersection of a sphere with other quadratic surfaces, and he also
includes their inverses with respect to spheres. Among the cyclic curves one
finds many of the classically well known curves, such as the Descartes ovals,
the scissoid, and the lemniscate. But, as a student at \'{E}cole Normale,
Darboux became familiar with the works of Lam\'{e}, Dupin and Bonnet on triply
orthogonal systems of surfaces, and this topic became, in fact, his major
interest for a long time.

At this time the best example of such surfaces was still the confocal surfaces
of degree 2 (see note to letter 13.12.70), and the surfaces were all defined
by a single equation. Years before, Kummer had studied analogous families of
plane curves, $f(x,y,s)=0$, $s$ a parameter, namely for each point there are
two curves passing through it, and they meet orthogonally. He also found that
the curves had to be confocal, that is, they have the same foci. In fact, the
orthogonality and the confocal property amount to the same thing. For example,
the orthogonal family of Descartes ovals have three common foci. Here we shall
simply remark that the definition of foci can also be extended to algebraic
curves of degree $n>2$. On the other hand, the situation is different in
3-space, and Darboux was able to construct a triply orthogonal family which is
not confocal.

On August 1,1864, Darboux presented to the Academy of Science his discovery of
the following triply orthogonal and confocal system of surfaces, expressed by
one algebraic equation of degree $4$%
\begin{equation}
\mu(x^{2}+y^{2}+z^{2})^{2}+\frac{\alpha\lambda-4h}{\alpha-\lambda}x^{2}%
+\frac{\beta\lambda-4h}{\beta-\lambda}y^{2}+\frac{\gamma\lambda-4h}%
{\gamma-\lambda}z^{2}-h=0 \label{cyclides}%
\end{equation}
where $\alpha,\beta,\gamma,h$, and $\mu$ are constants and $\lambda$ is the
parameter. For example, the case $\mu=0$ yields a 1-parameter family of
confocal quadrics. The above surfaces were subsequently called (generalized)
\emph{cyclides}, but their properties are quite different from those of
Dupin's cyclides. Clearly, the intersection of a cyclide with a sphere is a
cyclic curve.

However, on the same day as Darboux, the 15 year older Moutard announced to
the Academy that he had discovered the same system of surfaces. Moutard was an
expert on anallagmatic surfaces, namely surfaces invariant under an inversion,
and for surfaces of degree 4 he found that they were invariant under 5
different inversions if they contain the imaginary circle at infinity as a
double curve. The last statement simply means the 4th order terms have the
same form as in (\ref{cyclides}). But it is a geometric statement in our space
extended to the complex projective 3-space. Here, the portion of our extended
complex surface which lies in the plane at infinity is expressed by the
equation derived from (\ref{cyclides}) by ignoring all terms of degree lower
than 4, namely the equation%
\begin{equation}
x^{2}+y^{2}+z^{2}=0 \label{circle}%
\end{equation}
with muliplicity 2. This equation describes the imaginary circle at infinity,
as it was referred to in classical projective geometry. By seeking the lines
of curvature of his surfaces Moutard found the family of cyclides. We refer to
their papers Darboux [1864], [1865] and Moutard [1864a, b].

From his triply orthogonal system of cyclides Darboux derived the new
coordinate system $(\lambda_{1},\lambda_{2},\lambda_{3})$ of space, a kind of
generalized elliptic coordinates, and expressed the metric $\mathrm{d}s^{2}$
on the form (\ref{metric2}). As in the well known case of the confocal
quadrics, Darboux discovered that the intersection curves on the surfaces,
which by Dupin's famous theorem are lines of curvature, also form an isotermal
system of curves on each surface. But more importantly, he also came across
the following extension of Dupin's theorem, namely when two orthogonal
one-parameter families of surfaces intersect along lines of curvature, there
is always a third family of surfaces intersecting the other two orthogonally.
He used this result to find the condition a given family of surfaces, say
$f(x,y,z)=\lambda$, has to satisfy in order to belong to some triply
orthogonal system. His answer to this was a certain partial differential
equation of order 3 in two independent variables, but it was not calculated
explicitly since it was too complicated. But he did, for example, determine
those orthogonal system for which the lines of curvature are in a plane. All
these results were published in his classic memoir [1866], actually his first
memoir on orthogonal systems, and it was subsequently presented as his
doctoral thesis.

Certainly, orthogonal systems of surfaces is a field with which Darboux's name
will always be associated. But it was rather Cayley (1872) who first succeeded
in finding a tractable differential equation, which in a simple way determines
triply orthogonal systems, maybe of a special kind. However, Darboux quickly
analyzed and realized the essence of Cayley's work, which he extended to the
case of $n$ variables. Now he was able to prove various results on orthogonal
system in higher dimensions, for example, the determination of orthogonal
systems consisting of surfaces of degree two, or orthogonal systems containing
a \ given surface. The topic continued to play a major role in Darboux's
geometric works for many decades, as can be seen from his voluminous treatises
on surfaces (1878 and 1898).

However, as it often happens in mathematics, new techniques and viewpoints may
suddenly lead to substantial simplifications of previous hard work. Thus, it
is appropriate to mention the work of G.M. Green (1891--1919) in USA, who was
partly inspired by Darboux in France and E.J. Wilczynski (1876--1932) in USA.
The young Green wrote a 27 page paper in 1913 which largely overrides
Darboux's many pages devoted to triply orthogonal systems. He used a pair of
simultaneous partial differential equations of order two, but the new idea
comes from the projective geometric setting which he learned from Wilczynski.
In the late 1870's the French G.H. Halphen (1844--89) examined differential
equations invariant under projective transformations, and the topic of his
doctoral dissertation in 1878 was differential invariants. Apart from the
early investigations of Halphen, Wilczynski is largely regarded as the founder
of projective differential geometry, a geometric setting where he was the
first to demonstrate the utility of completely integrable systems of
homogeneous linear differential equations.

\subsection{Riemann and the birth of modern differential geometry}

Riemann's Habilitation lecture at G\"{o}ttingen in 1854 is generally regarded
as the birth of modern differential geometry, but it became generally known
only after its first publication in 1867. It is regarded as a classic of
mathematics, which is even more remarkable, since the audience of Riemann's
lecture was the Philosophical Faculty of G\"{o}ttingen, and its purpose was to
demonstrate lecturing ability. Certainly, the essay [1867] is almost devoid of
explicit mathematical content, but Gregory[1989], for example, argues
convincingly that the paper is better understood if we see Riemann speaking
primarily as a philosopher, but with a rather powerful mathematical methodology.

In fact, as a student Riemann had taken philosophy classes, and he was well
acquainted with, but also disagreed with, the Kantian view of space. As his
influences Riemann names only two persons, namely his supervisor Gauss, who
was in the audience in 1854, and the philosopher Johann F. Herbart
(1776-1841), who is also known as the founder of pedagogy as an academic
discipline. Herbart held the chair after Kant in K\"{o}nigsberg, and became a
professor of philosophy at G\"{o}ttingen University in 1833.

Gauss was probably the only one who realized the depth of Riemann's ideas.
According to Freudenthal[1970--1990], the lecture was too far ahead of its
time to be appreciated in those days, and it was not fully understood until 60
years later, when the mathematical apparatus developed from Riemann's lecture
provided the frame for the physical ideas in Einstein's general theory of
relativity. In the meantime, however, Riemann's ideas were a major source of
inspiration for many upcoming geometers, such as Beltrami, Helmholtz,
Clifford, Cristoffel, Lipschitz, Lie, Klein, Killing (1847--1923),
Poincar\'{e}, Ricci-Curbastro, and Levi-Civita. But Riemann himself did, in
fact, not publish any work on differential geometry in his lifetime. An
English translation by Clifford of Riemann's lecture appeared first in the
journal Nature (1873), but Spivak[1979], Vol.\thinspace2, also has an English
translation, supplemented with comments.

\subsubsection{Riemann's lecture and his metric approach to geometry}

The mathematical contents of Riemann's essay is often summarized roughly as a
generalization to higher dimensions of Gauss's results on the intrinsic
geometry of surfaces. But we shall have a closer look at the essay, which is
divided into three parts. First of all, he states that he is aiming at a
continuous space model, whereas the measurable properties of a discrete space
are simply determined by counting. Next, he observes the confusing status of
non-Euclidean geometry, which was not generally accepted at that time. He
attributes this problem to the fact that geometers do not distinguish clearly
between the topological and metric properties of space, which he discusses in
Part 1 and Part 2 of his address, respectively.

Here we use the modern term \textquotedblleft topology\textquotedblright\ for
the 20th century discipline which grew out from geometry and \textquotedblleft
analysis situs\textquotedblright\ during the decades after Riemann. Many
fundamental ideas in topology date, in fact, back to the works of Riemann, his
Italian friend Enrico Betti, and later also Poincar\'{e}, see for example
Betti\thinspace\lbrack1871], Poincare\thinspace\lbrack1895]. Riemann refers to
the underlying space as a multiply or $n$-fold \emph{extended quantity},
anticipating the modern concept of a (differentiable) $n$-dimensional
\emph{manifold,} and for simplicity we shall also use this term in the sequel.
In particular, the surfaces studied by Gauss are 2-dimensional manifolds. In
contrast to this, it should be noted that the notion of space is undefined in
the axiomatic development of geometry, although its properties are implicitly
described by the axioms.

Now, Riemann argues that topological considerations alone would not be
sufficient, say, to deduce Euclid's parallel postulate, and he points out that
experimental data are needed to determine the actual metric properties of
Physical Space, which he takes up in the final Part 3 of his address. Riemann
had a strong background in theoretical physics, influenced by the physics
professors W.E.\thinspace Weber (1804--91) and J.B.\thinspace Listing
(1808--82). In fact, he was Weber's assistant for 18 months, and besides,
Listing must also be counted among the early pioneers of topology.

According to Riemann, the measurable properties of Space are, after all, the
subject of geometry. The distance between two points is measured by physical
intruments, for example, one uses a rod or some optical instrument, and in
those days the measurable properties of space were found to agree completely
with Euclidean geometry. But the instruments and the notion of a rigid body
lose their validity when it comes to infinitely small distances where, in
fact, the metric may even disagree with the ordinary assumptions of geometry.
This is a possible scenario the physicists should be prepared to meet, and as
a preparation he presents a vastly general vision of geometry in Part 2.

In Part 2 Riemann displays a hypothesis on the metric structure of the space
which is as general as could be imagined at that time. Almost all mathematical
results appear here, but at the same time he is deliberately vague on many
points and avoids technical definitions and calculations, since many in the
audience had little knowledge of mathematics. According to Riemann, a given
manifold can be endowed with many different \emph{metric relations}, each
characterized by the form of the infinitesimal distance $\mathrm{d}s$, which
Riemann postulates to be expressible as a positive quadratic form%
\begin{equation}
\mathrm{d}s^{2}=\sum g_{ij}\mathrm{d}x_{i}\mathrm{d}x_{j} \label{metric}%
\end{equation}
in the differentials $\mathrm{d}x_{i}$ of the local coordinates, where the
coefficient functions $g_{ij}(x)$ may vary from point to point. If they are
constant, then the metric is Euclidean, and a linear change of coordinates
will transform the metric to the standard form%
\begin{equation}
\mathrm{d}s^{2}=\mathrm{d}y_{1}^{2}+\mathrm{d}y_{2}^{2}+...+\mathrm{d}%
y_{n}^{2} \label{diagonal}%
\end{equation}
which may also be regarded as the infinitesimal Pythagorean law. Therefore, by
continuity the geometry of (\ref{metric}) agrees optimally with the Euclidean
geometry in the vicinity of each point. Expressions of the type (\ref{metric})
have become known as a \emph{Riemannian metric}, and by the term
\textquotedblleft space\textquotedblright\ we shall mean a manifold with a
given metric of this kind. The global geometric properties, such as a distance
function measuring the distance between two points, will follow from the
metric by integrating $\mathrm{d}s$ along curves. For Riemann, the physical
space is merely an example of a 3-dimensional space and, contrary to the
Kantian viewpoint, he argues that the actual determination of its metric
(\ref{metric}) is a matter of physical measurements.

Next, Riemann turns to the metric properties behind the basic relation
(\ref{metric}) involving the $n(n+1)/2$ functions $g_{ij}=g_{ji}$. He argues
that the degrees of freedom of the functions are found by subtracting $n$ of
them, due to the freedom in the arbitrary choice of $n$ coordinate functions.
Thus Riemann concludes there is some set of $n(n-1)/2$ functions which
determine the metric completely, and he proposes to choose the functions
\begin{equation}
K_{1,2},K_{1,3},.K_{1,n},K_{2,3},..,K_{n-1,n} \label{curvature}%
\end{equation}
which at each point $p$ give the Gaussian curvature$\ $in $n(n-1)/2$
independent surface directions.

The above \textquotedblleft sectional\textquotedblright\ curvature $K_{ij}$ is
calculated as the Gaussian curvature of a surface $S_{ij}$ with the specified
surface direction at $p$, which the modern reader will interpret as the
2-dimension tangent plane at $p$. Implicit in Riemann's argument is the fact
that the number $K_{ij}$ only depends on the surface direction and not on the
actual choice of $S_{ij}.$ On the other hand, it is a major point of Riemann
that the \textquotedblleft sectional\textquotedblright\ curvatures
(\ref{curvature}), in turn, determine the metric relations (\ref{metric}).

Riemann pays most attention to spaces of \emph{constant curvature}, meaning
that all the numbers $K_{ij}$ are equal to the same constant $K$, valid for
any point. Here he gives the following \textquotedblleft
standard\textquotedblright\ expression for the metric in appropriate
homogeneous coordinates, namely
\begin{equation}
\mathrm{d}s^{2}=\frac{1}{(1+\frac{K}{4}%
{\displaystyle\sum}
x_{i}^{2})^{2}}%
{\displaystyle\sum}
\mathrm{d}x_{i}^{2} \label{K constant}%
\end{equation}
In particular, when $K=0$, as in the Euclidean plane or space, $\mathrm{d}%
s^{2}$ is the sum of squares of complete differentials, and Riemann proposes
to refer to these spaces as \emph{flat}. The round sphere is the familiar
example of a space with constant $K>0$, whereas Minding's pseudosphere has
constant $K<0$, but its connection with hyperbolic geometry was not known to Riemann.

\ In Part 3 Riemann focuses on the possible models for the physical space. He
points out that every determination from experience remains inexact, and this
circumstance becomes important when the empirical determinations are extended
beyond the observational limits into the immeasurably large or the
immeasurably small. Moreover, in the first case Riemann also distinguishes
between unboundedness and infinitude, noting that an unbounded universe may
possibily have finite size (volume).

Riemann argues that Space must be some unbounded 3-dimensional manifold, and
he does not exclude the possibility of having a variable curvature. This was a
radical idea which even many decades later was deemed too speculative by most
\textquotedblleft experts\textquotedblright. But first of all, he points to
further empirical evidence which put the spaces of constant curvature into the
forefront. Loosely speaking, these are the spaces which look the same at each
point and in every direction, and in terms of well known physical terms
Riemann describes this spatial property as \textquotedblleft free movability
of rigid bodies\textquotedblright, in the sense that they can be
\textquotedblleft freely shifted and rotated\textquotedblright. Thus, he
asserts the following two properties of a space are equivalent:%
\begin{align}
&  (i)\text{ \ \ the space has constant curvature;}\label{freemobile}\\
&  (ii)\text{ \ rigid bodies can be freely moved.}\nonumber
\end{align}
But notions such as \textquotedblleft rigid body\textquotedblright\ and
\textquotedblleft free movability\textquotedblright\ are physical rather than
geometrical, and Riemann makes no effort to explain statement (ii) in
geometric terms, as Helmholtz tried with his axiom $H3$, see Section 5.4.

\subsubsection{The beginning of modern differential geometry}

The publication of Riemann's lecture in 1867 was met with widespread acclaim,
and many were influenced by his ideas. Shortly afterwards, Helmholz in Germany
and Clifford in England published their own interpretations and extensions,
which also helped bring the attention to a wider community. As a physicist,
Helmholtz preferred to link the foundations of geometry with the physical
statement (ii) in (\ref{freemobile}), whereas others were challenged by the
necessity of having statement (ii) reexpressed purely in terms of derived
geometric concepts. This is closely related to the Riemann--Helmholtz space
problem, which we shall return to in Section 5.4.

Christoffel, Lipschitz and Schering were the first who started to elaborate
the Riemannian approach to geometry, based upon the postulated infinitesimal
structure of $\mathrm{d}s^{2}$. As a consequence, the transformation theory of
quadratic differential forms like (\ref{metric} ) became a central topic.
Whereas Riemann described in his lecture ---but with no calculations ---the
conditions for $\mathrm{d}s^{2}$ to be transformable to the flat metric,
Christoffel took the next step in 1869 by determining the necessary and
sufficient conditions for a quadratic differential form to be transformable
into another one, by a suitable change of coordinates. Lipschitz treated the
same problem in 1870, but the solution of Christoffel turned out to be more
useful. His analysis led him to the invention of a process which enabled him
to derive a sequence of tensors from a given one. This process was named
\emph{covariant differentiation }by Ricci in 1887. Riemannian geometry was
finally supplemented in the early 20th century by the important notion of
\emph{parallel transport, }which brings tangent vectors along closed curves
and expresses their total change of direction using the curvature tensor. This
notion was developed by Levi-Civita in 1917 and independently by J.A. Schouten
(1883--1971) in 1918. (See also note to letter 22.1.73).\ 

The modern geometer uses the Riemann curvature tensor $\{R_{ijkl}\}$ expressed
in the language of tensor calculus to investigate the curvature properties of
a space. But it is natural to inquire whether Riemann himself was ever in
possession of this tensor. In fact, some of the mathematical analysis
underlying Riemann's address in 1854 can be found in the second part of his
Pariser Arbeit, which is an essay he submitted in 1861 to the Academy of
Science in Paris in competition for a prize, announced in 1858 and relating to
a question on heat conduction. Again there is a quadratic differential
expression like (\ref{metric}), but not interpreted as a metric this time, and
Riemann is seeking the integrability conditions under which it can be
transformed to the simple form (\ref{diagonal}). He introduces expressions
which are essentially the components of the curvature tensor $\{R_{ijkl}\}$,
and he shows the integrability conditions are the vanishing of the components.
In particular, in the two variable case the condition is the vanishing of the
Gauss curvature, $K=0$.

But Riemann was not awarded the prize, perhaps because the way he obtained his
results was not satisfactorily explained. Neither was the prize awarded to
anyone else, and it was finally withdrawn in 1868. An English translation of
an extract from the prize essay's second part can be found in Spivak\thinspace
\lbrack1979], Vol.\thinspace2, pp.179--182. In 1872 Riemann's mathematical
papers came into the hands of Clebsch, who was Riemann's successor in
G\"{o}ttingen. But he died the same year, and they were temporarily passed
over to Dedekind \ and Weber. His Collected Mathematical Works were edited and
published by H.\thinspace Weber in 1876, but the prize essay appeared first in
the 2nd edition in 1892.

\section{Projective geometry}

\subsection{The origins of projective geometry\textbf{\ }}

Projective geometry is the first of the \textquotedblleft
new\textquotedblright\ geometries of the 19th century, naturally arising from
the classical Euclidean geometry. It is also the simplest and most fundamental
of these geometries. In a way, it is concerned with the aspect of figures that
remain unaltered when the observer changes his position, and even goes to
infinity. The geometry originates from the idea of \emph{perspective}, namely
the study of the geometric rules of perspective drawing, by which spatial
objects and relations in 3-space are projected onto a 2-dimensional plane. The
principles behind all this have been explored since antiquity.

The earliest and most basic projective invariant is the \emph{cross-ratio }of
four collinear points, see (\ref{ratios}); here we are using the modern term
introduced by Clifford (1878). This ratio has been shown, more recently, to be
the unique key invariant of projective geometry. In the surviving book after
Menelaus, called \emph{Sphaerica}, there is, in fact, a theorem in spherical
geometry which corresponds to the invariance of the cross-ratio. Appolonius
and Pappus were cognisant of the simpler theorem valid in the plane, but the
origins of that discovery is an open question. The invariance of the
cross-ratio under perspective projection can, in fact, be deduced from the
ancient theorem, already known to Thales (600 BC), saying that a line drawn
parallel with one side of a triangle cuts the other two sides proportionally.
Pappus wrote about Euclid's lost books, the Porisms, he gave 38 different
porisms and also suggested that Euclid knew about the invariance of the cross-ratio.

It has been quite a favorite sport among geometers to reconstruct, with
varying success, the lost ancient works. The work of Chasles, \emph{Les trois
livres de porismes d'Euclid} (1860) is at least recognized as an elaborate and
ingenious "restoration" of the porisms. H. Zeuthen, who studied with Chasles
in the 1860's, shared his interest in ancient Greek mathematics and he
investigated in detail the work of Apollonius on conic sections, such as his
method of determining the foci of a central conic. Moreover, in his treatise
\emph{Die Lehre von den Kegelschnitten im Altertum }(1886), Zeuthen suggests
as a possibility that the Porisms were just a by-product of a fully developed
projective geometry on conics.

\subsubsection{Developments in the 16th and 17th century}

Let us mention two important medieval scholars, both from Nuremberg in
Bavaria, namely Johann Werner (1468-1522) and Albrecht D\"{u}rer (1471-1528),
the latter best known as an artist painter. Werner worked on spherical
trigonometry and was maybe the last writer in the medieval tradition of conic
sections with some original contribution. But it does not seem that he knew
about the cross-ratio. D\"{u}rer was also one of the most important
Renaissance mathematicians; his remarkable achievements were through the
applications of mathematics to art, but he also developed new and important
ideas within mathematics itself. His masterpiece was the descriptive geometric
treatise on the human proportions, finished in 1523 but published
posthumously. The reason for his delay is largely because he felt it was
necessary first to write an educational elementary mathematical treatise,
which he published in 1525 as four books through his own publishing company.
These were, in fact, the first mathematical books published in German. In the
last book, for example, he wrote on regular and semi-regular solids, his own
theory of shadows, and an introduction to the theory of perspective. In 1527
he also published a work on (military) fortification, maybe as a response to
the threat of an invasion by the Turks felt by the people of Germany at that time.

It is remarkable that the celebrated Pappus's theorem (cf. Section 1.1.2)
naturally belongs to projective geometry, although the subject was not
developed until 1500 years later. But during this long time span there were
also a few other momentous discoveries of the same kind, notably the theorems
of Desargues (1639) and Pascal (1640). Inspired by the works of Appolonius and
Pappus, Desargues studied geometric objects such as conics from a new
viewpoint by focusing on\emph{\ }perspective or \emph{central}
\emph{projections}, and properties invariant under these.

To describe the basic theorem of Desargues, take two triangles $ABC$ and
$A^{\prime}B^{\prime}C^{\prime}$ in perspective from a point $O$, that is, the
three points $O,A,A^{\prime}$, resp. $O,B,B^{\prime},\ $resp. $O,C,C^{\prime}$
are collinear. Then the three intersection points of the corresponding sides
of the triangles, $AB$ and $A^{\prime}B^{\prime}$, $BC$ and $B^{\prime
}C^{\prime}$, $AC$ and $A^{\prime}C^{\prime}$, when suitably extended, are
lying on the same line.

Although the work of Desargues did not receive much acclaim by his
contemporaries, there were important exceptions and his lectures in Paris
influenced French geometers such as Descartes and Pascal. The latter published
an essay, at the age of 17, with several projective geometric theorems, such
as the famous Pascal's\emph{\ hexagon theorem} saying that for any hexagon
inscribed in a nondegenerate conic, the three points of intersection of the
opposite sides are collinear. Pappus's theorem is actually the other case
where the conic degenerates into two lines.

A central projection may well map parallel lines to intersecting lines, and
therefore, in the new geometry parallel lines are no longer special. This
suggests that parallell lines behave as if they intersect at some common ideal
point "at infinity", an idea already suggested by work of Kepler (1571--1630)
and Desargues. The introduction of these ideal points, one for each
"direction" of lines, yields the projective plane as an extension of the
Euclidean plane. The basic geometric objects are still the points and lines,
but there is one new line, called the line at infinity, consisting of all
those new points. The theorems of Pappus, Menelaus, Desargues and Pascal may
as well be regarded as statements valid in the extended plane. In fact, here
they become simpler than in the Euclidean plane, with no exceptional cases
since two lines always intersect.

Desargues and Pascal also knew about the invariance of the cross-ratio.
Desargues introduced notions such as harmonic sets and involution, and he
carried the ancient \emph{polar} theory about poles and polars, much further
than Appolonius. In plane geometry, this is a construction which uses a given
conic (ellipse, parabola, or hyperbola) to associate a line (polar) to each
point (pole) and vice versa. Namely, the conic has two tangents, say at
$p_{1}$ and $p_{2}$, passing through a given point $p_{0}$, and the line
through $p_{1}$ and $p_{2}$ is the polar of $p_{0}$. This correspondence is
involutive, and moreover, it has the property that the polars of the points on
a line $\lambda$ constitute the lines passing through the pole of $\lambda$.
Desargues started with a circle as the given conic, and he treated a diameter
as the polar of a point at infinity. Here we mention that if the circle is
mapped onto a conic by a central projection, the image of a pair of mutually
perpendicular diameters will be a pair of \emph{conjugate }diameters in the
sense of Appolonius. Desargues also introduced self-polar or
\emph{self-conjugate }triangles, namely each side of the triangle lies on the
polar of the opposite vertex.

However, the rather weird terminology used by Desargues has not survived,
whereas the terms \emph{pole }and \emph{polar}, still used today, were first
introduced by his later compatriots F. J. Servois (1768--1847) and Gergonne,
in 1811 and 1813 respectively. Before that, however, Euler, Legendre, Monge,
and Brianchon had also used the pole-polar construction. Monge and his
students Servois and Brianchon were in the forefront of projective geometry at
the beginning of the 19th century, but Servois is perhaps better known for
initiating the algebraic theory of operators, and he came close to discovering
the quaternions before Hamilton.

The late 17th century scholar in Paris, Philippe de La Hire (1640--1718), was
originally an artist (painter), and his interests in geometry arose from his
study of perspective in art. Perhaps he deserves to be considered, after
Pascal, a direct disciple of Desargues in projective geometry. His famous
treatise \emph{Nouvelle m\'{e}thode} (1673), which clearly displays the
influence of Desargues, is a broad projective approach to the study of conic
sections. Utilizing his own method of projection, he also reproved all 364
theorems of Appolonius. Somewhat strangely, however, he did not mention
Desargues, claiming he was not aware of the latter's work until after the
publication of his own.

In Italy, the geometer Giovanni Ceva (1647--1734) was interested in the
synthetic geometry of triangles, and he is known for the rediscovery of the
ancient Menelaus's theorem (\ref{Menelaus}) as well as his own celebrated
\emph{Ceva's theorem} from 1678. The latter theorem states that three lines
from the vertices $A,B,C$ of a triangle to points $P,Q,R$ on the opposite
sides, respectively, are concurrent precisely when the product of the ratios
in which the sides are divided is equal to 1:%
\begin{equation}
\frac{AR}{RB}\frac{BP}{PC}\frac{CQ}{QA}=1. \label{Ceva}%
\end{equation}
In fact, Ceva's theorem was known to early Arab mathematicians, and it dates
back to an 11th century king of Saragossa.

\subsubsection{Euler's affine geometry and Monge's descriptive geometry}

The subject of projective geometry dragged along and fell into oblivion, say
for the next 100--150 years after La Hire and Ceva. However, as a prelude to
the resurrection of projective geometry in the early 19th century, it is
timely to recall first the important step made by Euler in the 18th century,
when he investigated those properties of Euclidean plane figures which remain
invariant under \emph{parallel} projection from one plane to another. He
coined the term \emph{affine }for this purpose. Thus he initiated \emph{affine
geometry }as a kind of geometry featuring the parallelism in Euclidean
geometry.The affine and the Euclidean plane are identical as sets, having the
same points and lines. But the notion of angle is undefined, and comparison of
lengths in different directions is also meaningless. The generality of affine
motions, generated by all possible parallel projections, is simply illustrated
in a fixed Cartesian plane model, where the affine motions are coordinate
transformations of type
\begin{equation}
(x,y)\longrightarrow(ax+by+e,cx+dy+f),\text{ \ \ }ad-bc\neq0, \label{affine}%
\end{equation}
whereas suitable (orthogonality) conditions on $a,b,c$, and $d$ are needed for
a Euclidean motion. So, affine geometry is to be regarded as a "relaxation" of
Euclidean geometry.

Next, let us return to descriptive geometry, which originated with the
medieval works of D\"{u}rer, to overcome the problems of projection and to
describe the movement of bodies in space. But his ideas were not put on a
sound mathematical basis until the work of Monge. As a young man in the early
1760's, Monge studied orthogonal projections in 3-space and represented a
figure by its \textquotedblright shadows\textquotedblright\ in mutually
perpendicular planes. Then he devised a method to reconstruct the original
figure from the "shadows". Thus orthographic projection, the graphic method
used in modern mechanical drawing, evolved from Monge's simple scheme, and the
discipline became known as \emph{descriptive geometry.}

The story goes that Monge's ideas originated from a problem on fortifications,
which he solved for the French military. This reminds us of D\"{u}rer's
treatise in 1527, which also dealt with fortifications but probably outdated
by Monge 250 years later? Monge's solutions were so successful that the
military kept the method secret for 30 years and forbade Monge to publish
them. His own account on the subject appeared around 1795, when \'{E}cole
Polytechnique was established. In France, several outstanding geometers were,
in fact, educated at military (artillery) schools and were involved in the
Revolution or the ensuing military activities. In particular, quite many
students at \'{E}cole Polytechnique became soldiers or officers in Napoleon's
army and later played important roles in the political and academic life.

Monge gave lectures on descriptive geometry at \'{E}cole Polytechnique, where
this topic became a permanent part of the curriculum, and he emphasized
geometric visualization of mathematical and physical problems. The French
style soon became a model for other schools and military academies in other
countries, including the United States, and descriptive geometry became a
major topic. However, with regard to visualization of the geometry, T. Olivier
(1793-1853) went even beyond Monge by building a geometric model collection
for pedagogical purposes. Some of the models were ruled surfaces, even with
moving parts to illustrate how the surfaces are generated, others were
designed to illustrate how curves arise as the intersection of certain
surfaces. Selling models, particularly in the United States, became an
"industry" which gave Olivier quite a good income. In Germany this kind of
enterprise was promoted, for example, by Pl\"{u}cker, who had studied a year
in Paris around1823 and was Klein's influential teacher in Bonn in 1867--68.
At Clebsch's suggestion, in 1869 L.C. Wiener (1826--96) constructed plaster
models of cubics surfaces and others, which were later exhibited in London and
Chicago. Klein also became a proponent for plaster and string models in
Germany, and in the early 1870's his interest in models was expressed in
letters to Lie.

Carnot was Monge's early student, who also joined him in establishing the
\'{E}cole Polytechnique. He was teaching there, with a strong engineering
background, but is still best known as a geometer. His desire was to overcome
the increase of generality due to the algebraic methods of Descartes, so he
tried to simplify pure geometry and give it a universal setting. Of particular
interest are his ideas about \emph{correlative figures}, obtained by
continuous deformation of a given figure. By first establishing geometric
relations in the simplest case, where the involved quantities are positive
numbers, with no further restrictions he assumed the relations are identities
which still hold when the figure is replaced by a correlative figure. The
convenience of this principle or rather \textquotedblleft
axiom\textquotedblright\ was demonstrated by deducing generalizations of the
theorems of Menelaus and Ceva, and by establishing several theorems of
Euclid's \emph{Elements} from one single theorem. He published these results
in 1801--1803, and the principle is, in fact, a forerunner of Poncelet's
continuity principle, see below. Carnot's military masterpiece, like that of
D\"{u}rer and Monge, was also on fortification, published in 1809.

\subsection{The rise of projective geometry in the 19th century}

At the end of the 18th century Euclidean geometry was still the basic frame
for geometric thought, but with the turn of the century the situation changed
dramatically. After all, our visual world has the geometry of a projective
rather than a Euclidean space, so many geometers believed in the points at
infinity and perhaps also regarded the fundamental concepts of geometry to be
projective. Gradually, the conspicuous beauty and elegance of projective
geometry made it a favorite study among many geometers, who swarmed into the
new \textquotedblleft gold field\textquotedblright\ and quickly uncovered the
most accessible treasures. Thus, the rise of projective geometry made it
synonymous with modern geometry of the 19th century, and during the first
decades geometers in France and Germany played a leading role. Paris continued
to be at the center of the scene, and the priority of French mathematicians in
the creation of projective geometry cannot be denied. Monge's role was very
influential, both as a director and teacher at \'{E}cole Polytechnique. He can
be said to be the first modern specialist in geometry as a whole. Other
prominent geometers were Carnot, Poncelet, Gergonne and Chasles in France,
whereas M\"{o}bius, Steiner, von Staudt, and Pl\"{u}cker were at the forefront
in Germany.

\subsubsection{Poncelet and the creation of projective geometry}

In Paris around 1800, unexpected geometric ideas dating back to the 17th
century and even back to ancient times, were rediscovered and further
investigated. C.J. Brianchon (1783-1864) discovered the long forgotten
Pascal's hexagon theorem shortly after he entered \'{E}cole Polytechnique.
This led him to another hexagon theorem as well, by a skillful application of
the ancient polarity associated with a conic (cf. Lie-Scheffers [1896], Kap.1,
\S 3). \emph{Brianchon's theorem} (1806) says that for any hexagon
circumscribed about a conic, the three diagonals meet at a common point. In
fact, the theorems of Pascal and Brianchon are the first clear-cut significant
example of a pair of dual theorems. Brianchon's theorem was later generalized
to the case of a $(4n+2)$-gon, by M\"{o}bius (1847). But from the long list of
geometers who originated from the school of Monge, Poncelet ranks first and
foremost (cf. e.g. Darboux[1904: 101]).

Poncelet's early mathematical career is particularly interesting. He had
learned about the works of Monge, Carnot and Brianchon at \'{E}cole
Polytechnique, from which he graduated in 1810, at the age of 22. But he was
older than usual, due to health problems, and now he chose a military career.
Being trained as an engineer he took part in Napoleon's ill-fated invasion of
Russia in 1812, where he participated in the terrible battle at Krasnoy and
barely survived. As a prisoner, kept at Saratov until the defeat of Napoleon
in 1814, Poncelet discovered and wrote down the basic principles of projective
geometry, which were kept in his notebook when he finally returned to France
in the fall of 1814. During the following years he developed his new ideas in
a systematic way, being employed as a Captain of Engineers and a teacher of
mechanics at Metz, until 1825 when he accepted the position as professor of
mechanics. Poncelet published many articles on geometry and mechanics, but he
also continued his miltary career, and he was highly regarded for his
mechanical inventions.

During his early studies as well as in the Saratov notes, Poncelet applied
analytic geometry. However, for some reason, after his return to Paris he
changed his taste and became a staunch advocate of synthetic geometry. His
famous \emph{Trait\'{e} } (1822) ignited a tremendous surge forward in the
geometrical developments, which also took place at major universities abroad.
With his studies of the relationship between a figure and its image by central
projections, Poncelet took the final step towards a precise mathematical
description of the ancient geometric conception of perspectivity. In fact,
Desargues had initiated this subject in 1639, but his forgotten work did not
come to light until 1845.

For example, when a plane figure is illuminated\ by light\ rays emanating from
an outside point, its shadow in any other plane is a projective image of the
figure. As we have already observed, the parallel\emph{\ }projections studied
by Euler are the special case of central projections from ideal points
infinitely far away. By collectively referring to properties of figures
preserved by all these induced transformations as \emph{projective}, Poncelet
actually introduced a completely new discipline called projective geometry. We
mention that M\"{o}bius and Chasles, who developed the theory in a different
way, used the alternative terms \emph{collineation }and \emph{homography} for
a projective transformation, respectively. Other, such as von Staudt, also
used the term collineation.

Whereas M\"{o}bius and Chasles used analytic methods without hesitation,
Poncelet was strongly against the use of coordinates. But to achieve the
generality of analysis he found it necessary to introduce into synthetic
geometry imaginary points as well as ideal points. Thus he made a bold attack
on imaginary points, with a courage and thoroughness far ahead of his
predecessors. For this purpose he built upon Carnot's idea about correlative
figures, and in his \emph{Trait\'{e} }he introduced the \emph{Principle of}
\emph{Continuity}, a term coined by himself. Another important principle,
which also came in the limelight in the 1820's, is the Principle of duality,
and we shall return to both of them below.

Poncelet wrote about imaginary points and lines without having a general
definition, although occasionally he gave a rather complicated geometric
definition. Anyhow, the principle of continuity paved the way for the
introduction of imaginaries\ in geometry, whose geometrical interpretation
sometimes had a quasi-mystical status, exemplified by Steiner's reference to
the \textquotedblright ghosts in the shadowy kingdom of
geometry\textquotedblright\ (cf. Rowe [1989: 212]).

Poncelet announced for the first time one of the basic principles of modern
geometry, namely that every circle in the plane passes through two immovable
imaginary points, known as the \emph{absolute} points or \emph{circular
}points at infinity. They are common to all circles in the Euclidean plane. In
the same vein he also introduced the \emph{spherical }circle at infinity,
which all spheres in Euclidean 3-space have in common. But it is, indeed, more
of a triumph for the analyst, such as Pl\"{u}cker, that he can easily
"calculate" these points using complex numbers. Namely, in terms of
homogeneous coordinates $x_{i}$ of the projectively extended plane or space,
the two loci of points are typically given by (see (\ref{A2}))
\begin{equation}
\text{(i) \ }x_{3}=0,\text{ }x_{1}^{2}+x_{2}^{2}=0,\text{ \ \ (ii) \ }%
x_{4}=0,\text{ }x_{1}^{2}+x_{2}^{2}+x_{3}^{2}=0 \label{circular}%
\end{equation}
Here and in the sequel we shall adhere to the usual meaning of homogeneous
coordinates, namely they are only determined up to a common non-zero factor.

Moreover, according to Poncelet, two conics intersect in four points, real or
complex, and two real conics with no real common point have two common
imaginary chords. As a curiosity, let us also mention the Poncelet--Steiner
theorem, postulated by Poncelet and verified by Steiner, that all Euclidean
constructions (with ruler and compass) can be carried out with ruler alone
plus a single circle and its center.

\subsubsection{The Principle of Continuity and the Principle of Duality}

The Principle of Continuity has, indeed, a long history, and it was observed
or enunciated by various scholars before Poncelet, in one form or another. It
dates at least back to Kepler in the early 16th century, and Leibniz stated it
as a general law applicable in a broad philosphical sense, see Kleiner [2006].
For example, Boscovich (1711-1787) enunciated and used the principle, and it
was used by Monge and Carnot before Poncelet applied it in 1813. Among
geometers, however, it was not generally accepted until 1822, when it was
formulated by Poncelet in his \emph{Trait\'{e}}.

Kline\thinspace\lbrack1972: 844]) contends that the principle of continuity
was, in fact, accepted during the 19th century as intuitively clear and
therefore had the status of an axiom. It was freely used by geometers and they
never deemed that it required a proof. A modern formulation amounts to the
statement that if an analytic identity involving a finite number of variables
holds for all real values, then it also holds by analytic continuation for all
complex values (cf. Bell [1945: 340]). However, although Poncelet could not
justify his use of the principle, he refused to present it as a simple
consequence of analysis. As a result, he became involved in lengthy
controversies with other mathematicians such a A.L. Cauchy (1789-1867), the
pioneer of real and complex analysis.

Let us return to the beginning of the 19th century, when Brianchon discovered
two closely related theorems. At that time neither Brianchon nor his
contemporaries realized the general principle underlying his discovery, namely
the \emph{Principle of} \emph{duality}. This is one of the corner stones in
projective geometry, and the use of it virtually doubles the geometric
\textquotedblleft harvest\textquotedblright,\ at one stroke and without extra
labour. In the projective plane the principle relies on the fact that two
lines have a unique intersection point, and conversely, two points determine a
unique line. In effect, by formally replacing \textquotedblright
point\textquotedblright\ by \textquotedblright line\textquotedblright\ and
vice versa, in any theorem involving only points, lines and the incidence
relation, one obtains an equally valid statement called the dual theorem, as
exemplified by the theorems of Pascal and Brianchon.

The duality principle was first questioned by Brianchon, but as a new
discovery the principle was claimed by both Poncelet and Gergonne. The latter
had discovered the duality principle by observing the symmetry of the
incidence relations between points and lines, thus anticipating the
self-duality of projective geometry which is so evident from the modern
axiomatic viewpoint. In the plane and the 3-space, the principle will apply to
all statements which do not involve metric properties, and the term
\emph{duality} was introduced by Gergonne to denote the relationship between
the original and the dual theorem. He was so obsessed by duality that he
modified submitted papers to his Annales, worked out the dual versions of
theorems and presented the mutually dual theorems side by side in two columns.
It is possible that Gergonne, after all, saw a deep independent geometric
principle, but he did nothing to establish a logical basis for it.

In Monge's \emph{G\'{e}om\'{e}trie descriptive} (1799), one of the topics is
the pole and polar construction in plane geometry, and Monge also gives
elegant proofs bringing in the third dimension. On the other hand, Poncelet
(1826) was the first to establish the duality principle using this theory, and
he exploited the principle to its limit. Depending on a given conic, the
\emph{polarity }construction amounts to a geometric transformation, namely the
\emph{polar} transformation which in a certain way maps points to lines and
lines to points and thus realizes the duality principle. This is, indeed, the
first example of a geometric transformation taking geometric objects of one
type to objects of another type. The correspondence between points and lines
was called a \emph{correlation} by Chasles, a term still used today, but more
generally, this is a projective mapping of the projective n-space onto its
"dual" projective n-space.

However, Gergonne had discovered the same principle, but in a different way,
and for him the conics behind Poncelet's geometric resiprocation were of
subordinate importance. Now both claimed priority for the discovery of the
principle. But the underlying reason for it became a controversial topic and
there was a debate going on for many years, involving Poncelet, Gergonne,
Pl\"{u}cker, M\"{o}bius, Chasles and others, before the principle was finally clarified.

In fact, duality in the projective plane also extends to a duality between
plane curves, because the points along a given curve yield by polarity a
family of lines which envelop another curve, called the \emph{dual curve}. It
should be noted that the dual of a conic is also a conic, because the conic
has two tangents passing through any point off the curve, which by polarity
are mapped to two collinear points on the dual curve.

We mention that Pl\"{u}cker used the term \emph{reciprocity} when he referred
to duality. In the second volume of [1828-31] he discussed this topic, based
on his new idea of taking lines as the fundamental geometric \emph{Elements}.
He displayed the reciprocity at work in the geometry of conics, treated as
envelopes of lines and expressed in terms of homogeneous line coordinates. It
is also noteworthy that decades later, Pl\"{u}cker and Lie constructed many
types of \emph{reciprocities} or geometric transformations which generalize
the above duality principle. The simple idea, related to the same principle,
is that families of points and lines may be intimately related and
reciprocally associated with respect to a curve. The latter may be swept out
in two ways, (i) either as generated by the motion of a point, or (ii)
enveloped by the turning motion of a straight line (the tangent). In
particular, a curve will be dual to another curve, in a way preserving
tangency, and this leads us to \emph{contact transformations}, contact
geometry, and the geometrical works of Lie in the early 1870's. In his more
general setting the old principle appears as a simple special case, namely a
\textquotedblright linear\textquotedblright\ version which Lie always referred
to as the Poncelet-Gergonne's reciprocity principle.\emph{\ }

\subsubsection{The cross-ratio and harmonic sets}

As Appolonius must have observed, in Poncelet's geometry ellipses, parabolas
and hyperbolas are congruent\ since they are obtained by central projection
from a circle. In such a generality, one may wonder what properties of figures
are invariant under all projective mappings, and hence are independent of
magnitude, distance and angle. To this end, consider the following two types
of ratios, involving three (resp. four) different collinear points: \emph{\ }
\begin{equation}
(P,Q;R)=\frac{PR}{QR},\text{ \ \ \ }(P,Q;R,S)=\frac{(P,Q;R)}{(P,Q;S)}%
=\frac{PR}{QR}\frac{QS}{PS}\text{\ \ } \label{ratios}%
\end{equation}
First of all, the 3-point ratio $(P,Q;R)$ is easily seen to be invariant under
affine, but not projective transformations. On the other hand, the ancients
(e.g. Appolonius, Pappus) knew about the double ratio $(P,Q;R,S)$, namely the
\emph{cross-ratio} , and its invariance under projections. In France, this
ratio was also used by Desargues and Pascal in the 17th century, but with the
exception of Brianchon, leading French pioneers of projective geometry such as
Poncelet and Chasles were probably unaware of this invariant at the time when
it appeared in the works of their German counterparts Steiner and M\"{o}bius.

Chasles got the idea of cross-ratio through his attempts to understand
Euclid's lost works, and he knew that Pappus also had the idea. Concerning the
17th century scholars, he had seen what La Hire had written, but the works of
Desargues had escaped him first. Independent proofs of the invariance of the
cross-ratio were given by M\"{o}bius (1827), Steiner (1832), and Chasles
(1837). M\"{o}bius and Chasles referred to the ratio as the
\emph{Doppelverh\"{a}ltniss }or \emph{anharmonic ratio} (rapport
anharmonique), respectively.

M\"{o}bius, Steiner and Chasles also introduced the cross-ratio of four lines
of a pencil, that is, lines in a plane meeting at a point $O$. In fact, this
approach is seen to be more basic and it also yields a simple proof of the
projective invariance of the cross-ratio in both cases. Suppose a line $l$
intersects the lines $p,q,r,s$ from the pencil at the points $P,Q,R,S$, and
let $pq$ denote the angle between $p$ and $q$ at $O$. Then by simple
trigonometry
\begin{equation}
(P,Q;R,S)=\frac{PR}{QR}\frac{QS}{PS}=\frac{\sin pq}{\sin pr}\frac{\sin
rs}{\sin qs} \label{cross}%
\end{equation}
and the right hand side is defined to be the cross-ratio of the lines. On the
other hand, the invariance of the cross-ratio under projection also follows
from the expression (\ref{cross}). Steiner did not consider negative
quantities in his geometry, so for him the angles in (\ref{cross}) are
positive, whereas M\"{o}bius and Chasles considered segments to be oriented
and the cross-ratio to be a signed quantity.

Both Steiner and Chasles used the cross-ratio as a basic tool to characterize
the distinguished family of curves of order two, namely conics, and their
results were essentially the same. In his \emph{Trait\'{e} de sections
coniques} (1865), Chasles considers the cross-ratio of four lines passing
through four given points on a conic and a fifth point. He finds that the
cross-ratio has a fixed value as long as the fifth point also lies on the
conic, which enableds him to give a projective definition of conics in terms
of the cross-ratio. As a consequence, a homography must map a conic to another conic.

The calculation of the (signed) cross-ratio (\ref{ratios}) involves
measurement of (directed) segments along a fixed line. The crucial property is
that a projective transfomation of the line does not affect the cross-ratio.
However, its value $\lambda$ obviously depends on the order of the four given
points; in fact, there are up to six different values and they are $\neq0,1$,
namely
\begin{equation}
\lambda,1-\lambda,\frac{1}{\lambda},\frac{1}{1-\lambda},\frac{\lambda
-1}{\lambda},\frac{\lambda}{\lambda-1}, \label{values}%
\end{equation}
The set $\left\{  P,Q,R,S\right\}  $ is called a \emph{harmonic} quadruple, or
the points are said to be in harmonic position, if the cross-ratio is $-1,2$
or $1/2$. More specifically and fixing the ratio to be $-1$, the pairs $(P,Q)
$ and $(R,S)$ are said be to harmonically related to each other, or $S$ is
said to be the \emph{harmonic conjugate} of $R$ with respect to $(P,Q).$

Poncelet developed his theory of harmonic separation at great length, although
the invariance of the cross-ratio itself had somehow escaped him, at least he
did not use it. For given points $P,Q,R$, there are, indeed, simple coordinate
free constructions of the 4th harmonic point $S$, making the cross-ratio equal
$-1$, see for example Gray [2007: 26]. Using analytic geometry, on the other
hand, with $0,q,r,s$ as the coordinates of the points $P,Q,R,S$ on the real
line, a simple calculation yields
\begin{equation}
(P,Q;R,S)=-1\text{ \ }\Longrightarrow s=\frac{qr}{2r-q}. \label{harmonic}%
\end{equation}
Von Staudt also made extensive use of harmonic sets; for example, he simply
defined projectivity between two lines as a map preserving harmonicity.
However, this was a daring step, and critical questions about his approach
came up in the 1870's, involving Zeuthen, Klein, Darboux and others.

The invariance of the cross-ratio under projective transformations was a
useful fact for geometers who developed projective geometry. Chasles
introduced the term homography to describe a transformation between planes
which carries points into points and lines into lines. To ensure that such a
transformation is also projective, in the sense of Poncelet, he added the
extra condition that the cross-ratio must be preserved. Rather surprisingly,
however, this was later shown to be a superfluous condition.

By definition (\ref{ratios}), the cross-ratio involves the concept of length
of segments of which the ratio is compounded, whereas projective geometry was
supposed to be more fundamental than Euclidean geometry and hence should not
involve the concept of distance at all. The fundamental criticism of the work
of M\"{o}bius and Chasles also referred to this fact. However, some geometers
realized that direct usage of the general cross-ratio (\ref{ratios}), which
involves measurements, could be circumvented by successive construction of 4th
harmonic points. Still, geometric purists such as von Staudt found the
dependence of the cross-ratio on metric concepts intolerable, and he was the
first to advance the study of projective geometry in a way independent of
metric considerations. After his death in 1867, it turned out that his
pioneering work had taken him close to a final solution.

\subsection{ The analyst and the synthesist}

Since the beginning of analytic geometry in the 17th century, geometric
objects could be studied using tools from algebra and analysis, as an
alternative to the synthetic approach in the tradition of Euclid's
\emph{Elements}.\emph{\ }The two approaches --- analytic contra synthetic---
persisted side by side in the development of projective geometry in the 19th
century, and thus a distinction was made between analytic or algebraic
geometry on the one side, and synthetic geometry on the other side. The
\emph{analysts }would gladly use analytic or algebraic techniques from other
areas of mathematics, and they would express geometric relations in terms of
coordinates and equations. The \emph{synthesists}, on the other hand, inspired
by the ancient Appolonius, were the advocates of the purely geometric methods,
with intuition as a guide and logic as the instrument for a strict formal
reasoning, avoiding measurements and algebra.

Some of the leading geometers were, in fact, clearly favoring one of the two
major trends. Prominent synthesists were Carnot, Poncelet, Steiner, von Staudt
and Cremona, whereas M\"{o}bius, Chasles, Pl\"{u}cker, Cayley and Salmon were
heading the analytic trend. To the latter we may also include Grassmann, Hesse
and Clebsch, who pioneered the use of algebra and analysis in a modern sense.
Gergonne attributed the modern rise of coordinate geometry to Monge, but the
latter was also well acquanted with pure geometry, for example through his
treatise on descriptive geometry (1799), and wisely he chose to remain neutral
with regard to the ensuing controversy between the purists and the analysts.

However, the controversy also caused some bitter rivalry and priority
quarrels, and there were upcoming rumours, even with undertones flavored by
nationalism. But apparently, M\"{o}bius and Chasles were among the more
generous and diplomatic ones, and M\"{o}bius and von Staudt would stay aloof
from discussions not of a purely scientific nature. We should add that
Chasles, allegedly an analyst, also defended pure geometry. According to
Kline\thinspace\lbrack1972: 850], Chasles thought analytically but presented
his results geometrically. He calls this approach the \textquotedblleft mixed
method\textquotedblright, and says it was used later by others. Presumably,
this also includes eclecticists such as Sophus Lie, who thought synthetically
but wanted to present his results analytically.

\subsubsection{Poncelet, Gergonne, and Chasles}

As we have seen, Poncelet himself initiated the purely geometric approach,
having returned from the Russian battlefields and soon became convinced that
analytical methods are inferior to the synthetic ones. He set himself to undo
everything the successors of Descartes had done, and then he would reprove or
improve everything. But his dubious continuity principle actually promoted the
analytic trend as well, and in the later years of his career he seemed to be
more interested in mechanics and became more of an analyst.

There was also the saying that Poncelet and Steiner occasionally concealed in
synthesis what they had discovered by analysis. In fact, for some reason or
other, 50 years after they were composed, Poncelet decided to publish his
Saratov notes, and they appeared in two volumes in 1862--64, in their original
analytic form. In his address given at the Congress of Science and Arts in St.
Louis (1904), Darboux "let the analytic cat out of the synthetic bag" by
refering to the "unfortunate publication of the Saratov manuscripts", which
showed that the principles which served as a foundation for the
\emph{Trait\'{e} }(1822) were established by the aid of the Cartesian Analysis
(cf. Darboux [1904], Bell [1945: 339]).

Gergonne's mathematical career had been delayed due to his previous military
life, but after 1810 his journal and editorial position made him very
influential. He was advocating the use of coordinates, and to demonstrate the
power of the analytic approach he asked for proofs of classical problems of
synthetic geometry, such as the famous Appolonius problem: the construction of
a circle tangent to three given ones. His own solution in 1816 became known as
Gergonne's construction, and later many new elegant solutions were also also
found, with or without analysis, by Poisson, Pl\"{u}cker, Chasles, Poncelet,
Steiner, and others. Gergonne himself wrote about 200 articles, not
necessarily on geometry.

Chasles was the only follower of Poncelet of major importance in France, and
he devoted his entire life to projective geometry, following up the works of
Poncelet and Steiner. He was only four years younger than Poncelet, but from
his analytic approach to geometry he seemed to belong to a different
scientific epoch. Chasles was also a judicious historian of geometry, but in
his first major work \emph{Aper\c{c}ue historique} (1837), which is a classic
mathematical historiography, he admits that he had neglected the German
writers since he did not know the language. His interest in the past also made
him the first who fully appreciated the forgotten works of Desargues and
Pascal, but surprisingly, Poncelet disliked his appraisal of those 17th
century scholars. For some reason or other, the relationship between Poncelet
and Chasles was rather hostile.

Like many others from Monge's school, Chasles was also a professor at
\'{E}cole Polytechnique, for about ten years, and in 1846 a chair of higher
geometry at Sorbonne was specially created for him. He wrote a very important
text showing the power of synthetic geomertry, while his general study of
geometry, with all the concepts he had intoduced, such as cross-ratio,
pencils, and involution, appeared in his second major work, \emph{Trait\'{e}
de g\'{e}om\'{e}tri\'{e} superieure} (1852). It seems that he rediscovered or
superseded many of Steiner's results, with his own analytic approach, but
unintentionally since allegedly he did not know Steiner's papers.

Many classical counting problems are concerned with the enumeration of conics,
and some even date back to Appolonius. In 1848 Steiner posed, but wrongly
solved the problem of enumerating all conics tangent to five given conics. His
answer was 7776, but Chasles developed his theory of \emph{characteristics} to
solve the problem, and he gave the correct answer 3264 in 1864. At this time,
a Danish student named H.G. Zeuthen was in Paris for two years and studied
geometry with Chasles, having received a scholarship in 1863 after graduation
from the University of Copenhagen. His doctoral dissertation (in Copenhagen)
in 1865 was, in fact, on a new method to determine characteristics, and for
the next ten years he worked mostly on enumerative geometry. In 1868 Sophus
Lie met Zeuthen for the first time, a meeting which was very decisive for his
choice of future career as a mathematician. Both Lie and Klein developed a
lasting friendship with Zeuthen, who became Denmark's first internationally
acknowledged mathematician.

\subsubsection{M\"{o}bius, Pl\"{u}cker, and Steiner}

In Germany, M\"{o}bius was the first pioneer of projective geometry. After
finishing his studies in Leipzig he moved to G\"{o}ttingen in 1813 to study
astronomy under Gauss. In fact, throughout his career he taught mechanics and
had the title of astronomer. In the 1840's M\"{o}bius became a professor in
astronomy as well as the director of the Observatory in Leipzig. But he is
best known for his results in pure mathematics. M\"{o}bius had heard about the
geometric works of the French pioneers Poncelet, Gergonne and Chasles and he
also acknowledged his indebtedness to them. But M\"{o}bius rather followed his
own approach, which was also largely adopted by Chasles. These two
mathematicians, of very different personality, were regarded as scientific
equals and they headed the analytic trend in projective geometry for many years.

M\"{o}bius's basic work \emph{Der barycentrische Calcul}\textit{\ }(1827) on
analytic geometry became a classic. Here he presents many of his results on
affine and projective geometry, and he discusses projective mappings and other
geometric transformations. But most importantly, he introduces his
\emph{barycentric} coordinates, which are clearly inspired by the notion of
center of mass in classical mechanics. The key idea is to assign to each point
$p$ in the plane a triple of numbers $(m_{0},m_{1},m_{2})$, depending on a
fixed triangle with vertices $p_{0},p_{1},p_{2}$, such that $p$ will be the
center of mass of the triangle when the mass $m_{i}$ is assigned to the vertex
$p_{i}$. Clearly, the numbers are unique up to scaling and, moreover, some of
them must be negative if $p$ lies outside the triangle. In fact, the $m_{i}$'s
are \emph{homogeneous} coordinates applied to projective geometry for the
first time. For analysts such as Pl\"{u}cker, homogeneous coordinates became a
flexible tool in many types of coordinate descriptions, and with Pl\"{u}cker
the leadership of analytic geometry was definitely in Germany.

Like Felix Klein, Pl\"{u}cker also had his basic education in D\"{u}sseldorff,
but contrary to Klein's experience he was inspired by his gymnasium teacher to
study mathematics. The young student moved around to various major
universities such as in Bonn, Heidelberg and Berlin, and after completion of
his doctoral thesis in Marburg in 1823, at the age of 21, he went to Paris.
Here he attended courses in geometry and came under the influence of the great
school of French geometers, in the spirit of its founder Monge, with Poncelet,
Gergonne and Chasles as the leading figures in the development of projective
geometry. Returning to Germany, Pl\"{u}cker submitted in 1824 his habilitation
thesis to the University of Bonn, where he continued his research and further
qualified to become a university teacher (extraordinary professor) in 1829.
One of his great achievements dates back to this time, namely he proposed the
revolutionary idea that the straight line rather than the point may be used as
the fundamental geometric element. Thereby he initiated the new discipline
called \emph{line geometry}, a 4-dimensional geometry which represents a new
way to study the various geometric configurations in 3-space. His first memoir
on the subject was published in the Philosophical Transactions of the Royal
Society in London, where many of his later publications also appeared.

Pl\"{u}cker's first paper, in Gergonne's Annales (1826), was in fact a
synthetic approach to the tangents of conics, a favorite topic at the time.
But this also drew him into the crossfire between Gergonne and Poncelet, in
particular their dispute over the discovery of duality, a topic which also
appealed to Pl\"{u}cker. However, at this time he switched completely to the
analytic approach and soon he became the leading expert on both analytic and
algebraic geometry. He did not aim at collecting existing results, exploiting
existing principles, but modestly building analytic geometry anew and along
the lines suggested by Monge. In addition to several published books,
Pl\"{u}cker contributed many important papers to various periodicals in
Germany, France, England, and Italy, and terms such as "new method" or "new
geometry" typically appear in the titles, subtitles, or prefaces. Undoubtedly,
no single person contributed more to analytic geometry than Pl\"{u}cker, with
regard to both volume and power.

In Germany, however, the relations between Pl\"{u}cker and Steiner were far
from being friendly, and they became engaged in an endless feud. Having such
an extraordinary geometric intuition, Steiner was said to be the greatest
geometer since the legendary Appolonius, and he became the obvious leader of
the German school of synthetic geometry. In 1834 a chair of geometry was
established for him in at the University of Berlin, which he kept until his
death in 1863. His teaching was so influential that courses in projective
geometry at many universities have up to present time been based on his
outlines, and even some of his \textquotedblright
old-fashioned\textquotedblright\ terminology has survived. On the other hand,
in 1834 Pl\"{u}cker had just spent one year as an extraordinary professor in
Berlin, and it was perhaps expected that he too would try to make his career
here. However, now he quickly decided to leave Berlin and was offered the
position as an ordinary professor at Halle, where he stayed for two years
before he finally returned to Bonn in 1836.

The story about Jakob Steiner is rather peculiar, indeed. It starts with a
Swiss shepherd and farmer\ son who first went to school at the age of 18,
where his great talent for geometry was discovered, and later as a student at
German universities he managed to support himself precariously as a tutor. In
Berlin he became acquainted with the prosperous engineer and largely
self-taught mathematician Leopold Crelle, who had a special ability to spot
exceptionally talented young mathematicians and generously offering them his
friendship and support. Another of these young men was the Norwegian Niels H.
Abel, and together with Steiner they strongly encouraged Crelle in the
founding of his mathematical journal in 1826, the first journal in Germany
devoted exclusively to mathematics. The first volume came out in 1827, filled
with many original works of Abel and Steiner. Here we find, for example,
Steiner's proof of the formula $(n^{3}+5n)/6+1$ for the number of pieces space
can be divided into by n planes. Steiner published altogether 62 papers in the journal.

Crelle's Journal became an important publication channel for mathematics in
general, and also the analyst geometer Pl\"{u}cker submitted papers to this
Berliner journal. In Bonn, Pl\"{u}cker was applying his analytic methods and
thus promoted an independent development of modern geometry. His success was
in a way comparable with that of his great contemporaries Poncelet and
Steiner, who continued to cultivate geometry in its purely synthetic form.
However, at his chair of geometry in Berlin, Steiner became increasingly more
influential and authoritative, and Crelle seemed to favour Steiner over
Pl\"{u}cker with regard to their personal conflict. Perhaps somewhat
dubitable, but the story goes that Steiner would no longer write in the
journal if Pl\"{u}cker did so, and for many years Crelle was in reality forced
to deny Pl\"{u}cker access to the journal.

Pl\"{u}cker felt disappointed with the receptions his geometric works were
judged in Germany, where the impact of the synthetic approach of Poncelet and
Steiner was regarded to be more useful. Perhaps this was the reason why
Pl\"{u}cker switched over to natural science when he accepted the chair of
physics at Bonn in 1847. For another explanation of the switch, perhaps
filling the chair of physics with a mathematician would have been untenable.

It is likely that Pl\"{u}cker's accomplishments both as a mathematician and
physicist were rather unacknowledged in Germany during his lifetime, and
certainly English scientists appreciated his work more than his compatriots
did. In England, his reputation as a profound geometer flourished, and he was
encouraged by Cayley and Sylvester who dominated British pure mathematics in
the second half of the 19th century. For example, at a meeting of the British
Association for the Advancement of Science in 1848, Sylvester hailed
Pl\"{u}cker as "the master" of English mathematicians and, moreover, there was
none between Pl\"{u}cker and Descartes when it came to the relation of
geometry to analysis (cf. Gray [2007]). In 1863 Steiner died, and a couple of
years later Pl\"{u}cker again turned his attention to mathematics. We shall
return to some of his mathematical accomplishments later.

\subsubsection{Steiner, von Staudt, and Cremona}

Steiner was the first of the German school of geometry who took over the
French ideas, following the strict synthetic approach of Poncelet. He was
undoubtedly the most extreme synthesist at all, he was said to be hating
analysis and he would teach geometry without using figures. But he came to
impressive geometric results, especially in his younger days, and with his
first book \emph{Systematische Entwickelungen }(1832) he made a grand effort
to unify classical geometry, based on a new conception of projective geometry
and a new approach to conic sections. Thus he would seek the common roots and
uncover the fundamental properties of the classical geometry. Much of his
later works aimed at encompassing more recent results as well, due to himself
or others, into his synthetic geometric framework. Quite often, therefore,
Steiner was able to construct purely geometric proofs of results first
discovered in analytic or algebraic geometry.

A principal idea of Steiner is to use the simple projective concepts such as
points, lines, planes, pencils of lines or planes etc. to build up more
complicated structures, using the principle of duality and the cross-ratio
(\ref{cross}) as basic tools. Steiner's approach to duality is in the
tradition of Monge, Brianchon and Poncelet. First of all, he establishes a
powerful theorem on conic sections, which in fact yields a new defnition of
these curves. Namely, starting from two projectively related lines, the lines
between corresponding points on the two lines envelop a "conic" which is
tangent to the two lines. Moreover, from the proof of the theorem it is clear
that the "conics" are projections of a circle and hence are actually conics.
Now, from the pole and polar theory he can construct the polar map and realize
the duality principle. In particular, starting with two projectively related
pencils of lines, the locus of all intersections of corresponding lines will
be a conic passing through the centers of the pencils. Steiner's synthetic
theory of conic sections was one of his chief accomplishment in projective geometry.

Steiner investigated also algebraic curves and surfaces with his synthetic
approach. Let us first recall some basic properties of quadrics. A quadric has
the property that all its plane sections are conics, and moreover, if it
contains a line, then it is a ruled surface and hence contains infinitely many
lines. In fact, we know which quadrics are ruled (see Section 2.3.1). Now,
Steiner was looking for a surface of degree
$>$%
2 with the property that its plane sections are still conics, for example two
conics. On a trip to Rome in 1844 Steiner discovered his Roman surface, which
is such a surface and it has degree 4. In fact, it has the peculiar property
that each tangent plane cuts the surface in two conics and therefore a double
infinity of conics is lying on it (see Gray [2007: 270]).

Turning to (real) cubic surfaces, Cayley first showed that the number of
straight lines on the surface must be finite, which enabled Salmon to prove
there are exactly 27 lines in general, some of which may be complex. This is
the content of the Cayley-Salmon theorem, published in 1849. Even 20 years
later Cayley wrote on the topic of cubic surfaces, see his memoir [1869b], and
this may also have inspired Klein to his own studies in 1872-73. Both
Sylvester and Salmon had given a pure geometric construction of a nonsingular
cubic surface and its 27 lines, but without any indication of how to make a
model (thread or plaster), which was, in fact, constructed for the first time
by L.C. Wiener in 1869.

In 1856 Steiner wrote an important paper which was the basis for his purely
geometric approach to cubic surfaces and, for example, now he gave a synthetic
proof of the Cayley-Salmon theorem. However, he also stated many other results
without proofs, which were in fact supplied 7-10 years later by Cremona and
Sturm. They were highly regarded geometers in the generation after Steiner and
Staudt, working in the same spirit of pure geometry. Sturm studied the cubic
surface problems in his dissertation (Breslau, 1863) and continued with this
work, whereas Cremona during his years 1860-67 at Bologna did important work
on transformations of plane curves, birational transformations, and wrote a
long memoir on cubic surfaces which appeared in Crelle's Journal in 1868.\ As
was perhaps expected, in 1866 Cremona and Sturm were jointly awarded the
Steiner-Preis of the Berlin Academy. Cremona's transformation theory had many
links to Lie's work in the 1870's.

Von Staudt was the second major synthesist in Germany, a contemporary of
Steiner but his opposite in many respects. He studied under Gauss i
G\"{o}ttingen during 1818-22, and based on his work on the determination of
the orbit of a comet he received his doctorate from the University of Erlangen
in 1822. In 1835 he was appointed to this university and stayed there to the
end of his life in 1867. His lifestyle was modest, he communicated hardly with
people but continued peacefully with his own rigorous solitary research.

In his 1822 treatise, Poncelet had pointed out the distinction between
projective and metric properties of figures, namely the projective properties
are logically more fundamental, but it was von Staudt, rather than Poncelet
himself or Steiner, who began to build up projective geometry as a subject
independent of distance and hence without reference to length or angle size.
With his book \emph{Geometrie der Lage} (1847), supplemented by three booklets
(1856-60), von Staudt presented the essense of his project, aimed at
introducing an analogue of length on a projective basis. For information on
the publication history of this peculiar book, see Hartshorne [2008].

A major problem von Staudt was facing was how to define the cross-ratio, the
fundamental invariant of projective geometry, in an intrinsic way. Its
importance was also made clear in the works of M\"{o}bius, Steiner and
Chasles, but the cross-ratio involves four lengths and length is a metric
concept. Consequently, numbers must be associated to the points of a line, but
the numbers and their algebraic operations must be defined purely
geometrically. In the opposite direction, this would also pave the way for the
idea of having non-metric geometry on which a notion of distance can be
defined. In fact, the idea of a projective theory of metric geometry was first
elaborated by Cayley in his\emph{\ Sixth Memoir }[1858], see below.

From a modern viewpoint, it is a general coordinatization problem to represent
a geometric lattice by "closed subsets" of some algebraic structure, and von
Staudt was the first to realize the possibility of such a coordinatization in
projective geometry (cf. Lashkhi [1995]). It was for this purpose he created
his "Wurf algebra" or the "algebra of throws", where a "throw" is the
geometric construction which attaches a symbol (number) to any point. Starting
from three points on a line $l$, von Staudt constructs a \emph{harmonic chain}
of points on $l$, by first constructing the 4th harmonic point and then
iterating the construction indefinitely by choosing at each step three points
from the previously constructed points. Similarly, starting from four coplanar
points, with no three on the same line, he constructs a \emph{harmonic net}
(or M\"{o}bius net) of points in the plane. It is clear, however, that the
successive 4th harmonic construction was essentially the same as the ancient
construction of commensurable sets, so that the transition from the harmonic
chain (or net) to the whole line (or plane) would be similar to proving that
things valid for commensurable ratios hold as well for all ratios.

The coordinates (symbols) of the harmonic chain on $l$ can be taken to be
rational numbers,\ with their algebraic operations given by the appropriate
geometric constructions. Thus by identifying the points $P,Q,R,S$ in
(\ref{cross}) with their coordinates $\overline{p},\overline{q},\overline
{r},\overline{s}$, their cross-ratio can be calculated in the usual way as%
\[
(\overline{p},\overline{q};\overline{r},\overline{s})=\frac{(\overline
{r}-\overline{p})(\overline{s}-\overline{q})}{(\overline{r}-\overline
{q})(\overline{s}-\overline{p})}%
\]
We refer to Kline [1972: 850], which also illustrates the construction of the
point labelled $2$ on the line $l$, starting from chosen points labelled
$0,1,\infty$, where $\infty$ is the point at infinity. Since $(0,1;\infty
,2)=1/2$, this amounts, indeed, to a construction of a 4th harmonic point, see
(\ref{values}).

M\"{o}bius and von Staudt defined a general collineation to be a one-to-one
transformation which takes points to points, lines to lines (and planes to
planes, in dimension 3). Then they asked the question whether these are
necessarily "projective" transformations in the appropriate sense. In the
plane, M\"{o}bius assumed continuity, which enabled him to show the
transformation is a composition of perspectivities and hence is projective in
the sense of Poncelet; in particular it preserves the cross-ratio. On the
other hand, von Staudt had no continuity assumption, but he showed harmonicity
is preserved, which would enable him to conclude the cross-ratio is preserved
as well. But how could this be true ? What about the lowest dimensional case,
namely transformations on a single line ?

Most likely, the work of von Staudt was poorly understood and little
appreciated during his lifetime, although Reye lectured on his approach and
even published the lectures in 1866. It seems to be Klein, with his interest
in the foundations of geometry, who first focused attention on von Staudt and
saw his work in a new light in the early 1870's. Klein learned from his friend
Otto Stolz, who was close to von Staudt's spirit and was well acquainted with
the new geometric ideas of both von Staudt and Lobachevsky. Soon it became
clear that von Staudt had, in fact, made implicit use of the Euclidean
parallel postulate, which is a blemish since parallelism is not a projective
invariant. But the problem was not deeply rooted and Klein was able to remove
it, as explained in Klein's second paper [1872c] on non-Euclidean geometry.

However, in [1872c] Klein was also alerted by a more serious gap in von
Staudt's work, namely he did not really show that a harmonic chain fills the
line "densely", penetrating any small interval. This was needed to ensure that
a projective mapping of a line was uniquely determined by its restriction to
the harmonic chain. More precisely, the key result asserted by von Staudt
amounts to a statement nowadays referred to as the Fundamental Theorem of
projective geometry:

\emph{A projective transformation of a line to itself is uniquely determined
by its values at three different points. Alternatively, if the transformation
fixes three points, then it is the identity transformation.}

The transformations of the line considered by von Staudt were those preserving
harmonicity. Therefore, starting from three given points, the image of the 4th
harmonic point is already determined, and by repeated usage of the 4th
harmonic point construction the transformation is determined on the whole
harmonic chain. From this he simply concluded the transformation was
determined on the whole line.

We recall that Eudoxus had tacitly used a basic property of line segments,
referred to as the Archimedian axiom since the 19th century, which granted the
possibility of fencing an incommensurable ratio by commensurable ones, and the
validity of this step was in fact clarified in 1872 by the Cantor--Dedekind
axiom. Now Klein insisted in his article that the "density" property of the
harmonic chains might as well be postulated as a similar continuity axiom.
Klein's article drew immediate responses from Cantor, L\"{u}roth, and Zeuthen,
which Klein described and commented in his subsequent paper [1874a].

Of particular interest was the answer from Zeuthen, which included a proof of
the fact that a harmonic chain is, indeed, a dense set of points on the
projective line. Still, however, Klein insisted that a continuity assumption
was needed in von Staudt's definition of a collineation, for the same reason
that a function defined on the rational numbers may not be extendible to a
continuous function on the whole real axis. But the topological concept of
continuity was not so well understood in the 1870's, and earlier it had only
been handled in a confusing manner.

Finally, in 1880 Klein received a letter from Darboux with his proof that no
extra condition of continuity was, after all, needed in von Staudt's
definition of collineations. The situation was similar to other theorems which
were known to be true without an explicit continuity condition, such as for
example transformations of the plane which send all circles and lines to
themselves. M\"{o}bius had assumed a continuity condition to conclude that
these transformations are the same as inversions (\ref{inversion}) and their
compositions, but the extra continuity assumption is in fact superfluous. We
refer to Klein [1880] for Darboux's continuity argument, and to Gray [2007:
341--42]) for a readable account of the above events, including the arguments
of Zeuthen and Darboux.

Briefly, it may be said that the success of Steiner and von Staudt was due to
their unification and purification of geometry, respectively. In the 19th
century Steiner was usually ranked higher than his successor von Staudt. But
since the early 20th century Steiner seems to have fallen below his successor,
whose originality and depth have been hailed by many writers starting with
Klein [1928] and Coolidge [1934]. However, some new viewpoints are presented
in the recent paper Bl\aa sj\"{o} [ 2009].

Cremona studied mathematics under Brioschi and others in his hometown Pavia,
where he took a doctorate in civil engineering in 1853. But due to his
previous military activity against the ruling Austrian government he was
prevented from obtaining a position, so for the following three years he made
his living as a private tutor of mathematics. The situation improved in 1856,
with his second paper published, and in early 1857 he was secured a full
teaching position at the scientific high school in the capital Cremona of the
Lombardy region. He wrote a number of papers in the following years, examining
curves with projective methods which later became characteristic for his more
important mathematical works. In 1860 he was appointed professor at Bologna,
in fact, by Victor Emmanuel II who was proclaimed king of the recently united
Italy in 1861.

Cremona left Bologna in 1867, when he was appointed to the Polytechnic
Institute of Milan, on Brioschi's recommendation, and received the title of
Professor in 1872. But the next year he moved to Rome, as director as well as
professor of graphic statics at the newly established Polytechnic School of
Engineering. The years he spent in Milan were, however, the most creative time
of his mathematical career. Projective geometry was a theme which was present
in almost all his works, and it all culminated with his monumental book [1873]
on the \emph{Elements} of projective geometry. From now on his administrative
and teaching duties began to put an effective end to his research career.
Although he was appointed to the chair of higher mathematics in Rome in 1877,
the political pressure on him finally persuaded him to serve the new Italian
State henceforth at the political level.

The method of graphic statics, founded by C. Cullman at Z\"{u}rich around
1860, applies projective geometry, rather than analytic methods, to the design
and analysis of stationary mechanical frames or equilibrium of systems of
forces. Cremona gave important contributions to this topic, and improved
previous work of Maxwell. For example, Maxwell's notion of reciprocal figures,
which appeared in an engineering journal (1867), was interpreted by Cremona
(1872) as duality in projective 3-space. To begin with, Cremona's geometric
views were largely influenced by Chasles, but later he became more attracted
to the synthetic approach of Steiner and von Staudt.

Cremona [1873] first reviews the developments of projective geometry and the
contributions of the various pioneers. Then, based on his undestanding of its
present state and the previous efforts to establish projective geometry as the
most fundamental geometry, Cremona aims at pulling it all together anew, in
purely projective terms without any resort to Euclidean geometry. Much of the
basic work was, after all, already provided by Chasles, Steiner, and von
Staudt, but the great merit of Cremona is that he finally succeeds in
presenting projective geometry as an independent geometric system, released
from the embarassment of its Euclidean origins (cf. Gray [2007: 245]).

\subsection{Some basic analytic developments}

In the classical description of projective spaces $P^{n}$, say in dimensions
$n=1,2,3$, there is always the distinction between ordinary and ideal points,
sometimes called finite points and points at infinity, respectively. This is
because $P^{n}$ was regarded as an extension of the Euclidean n-space $E^{n}$,
or rather the affine n-space $A^{n}$ (since metric properties are
non-projective). The set $A^{n}$ consists of the ordinary points, whereas the
points at infinity, as it turned out, actually constitute a projective space
of one dimension lower. For example, the ordinary line $A^{1}$ extended by an
ideal point $\infty$ becomes the projective line $P^{1}$, so by writing
$P^{0}=\left\{  \infty\right\}  ,$
\begin{equation}
P^{1}=A^{1}\cup P^{0}\ \text{, \ \ }P^{2}=A^{2}\cup P^{1},\text{\ \ }%
P^{3}=A^{3}\cup P^{2} \label{extensions}%
\end{equation}

While the synthesists were developing projective geometry, the analysts
treated the same subject as coordinate geometry and pursued their own methods.
To be more specific, let us consider the projective plane $P^{2}$, where
xy-coordinates are introduced in the ordinary plane $A^{2}$ and thus we
identify it with the Cartesian plane $\mathbb{R}^{2}$. In this xy-plane model,
the affine transformations have the simple algebraic expressions
(\ref{affine}), whereas a general projective transformation, the joint effect
of both parallel and central projections, expresses as%

\begin{equation}
(x,y)\longmapsto\left(  \frac{a_{21}+a_{22}x+a_{23}y}{a_{11}+a_{12}x+a_{13}%
y},\frac{a_{31}+a_{32}x+a_{33}y}{a_{11}+a_{12}x+a_{13}y}\right)
,\ \det(a_{ij})\neq0 \label{projective}%
\end{equation}
Therefore, unless the transformation is affine there are two lines in the $xy
$-plane
\begin{equation}
l_{A}:a_{11}+a_{12}x+a_{13}y=0\text{, \ \ \ \ }l_{B}:b_{11}+b_{12}x+b_{13}y=0
\label{lines}%
\end{equation}
such that the mapping (\ref{projective}) is undefined along $l_{A}$ and its
image is the complement of $l_{B}$, which remains to be determined. This also
illustrates the role of the extended plane $P^{2}$ and its ideal line
$l_{\infty}=P^{1}$in (\ref{extensions}). Namely, the mapping (\ref{projective}%
) extends to a transformation which maps the line$\ l_{A}$ to$\ l_{\infty}$
and $l_{\infty}$ to $l_{B}$. In effect, we shall obtain an invertible mapping
of $P^{2}$, whose expression in terms of homogeneous coordinates becomes,
indeed, very simple, see (\ref{trans}) below.

This classical picture is somewhat misleading, however, because the
homogeneous structure of $P^{n}$ gives no preference to any specific
hyperplane $P^{n-1}$; in fact, any hyperplane can be mapped to any other one
by a suitable projective transformation. Therefore, one may choose a
hyperplane $\simeq P^{n-1}$ in $P^{n}$ and refer to it as the points at
infinity. Then the complementary set identifies with an affine subspace
$A^{n-1}$, and we recover the classical picture (\ref{extensions}). This
aspect of projective spaces became better understood with the introduction of
homogeneous coordinates. Moreover, parallel with the developments of linear
algebra and group theory came the modern vector space model of a projective
space, which has largely reduced projective geometry to the setting of
algebraic geometry. Below we shall \ explain the vector space model, see (78)- (81).

Projective geometry in higher dimensions $n>3$ was first introduced around
1850, in particular in its analytic form as created by M\"{o}bius\textbf{\ }%
and\textbf{\ }Pl\"{u}cker. They also introduced complex coordinates, extending
$P^{n}$ to the complex projective n-space $\mathbb{C}P^{n}$, an extension of
the Cartesian space $\mathbb{C}^{n}$ which is completely similar to the real
case extension of $\mathbb{R}^{n}=A^{n}$ to $P^{n}$.

\subsubsection{Algebraization and homogeneous coordinates}

The modern rise of coordinate geometry was due to Monge, according to
Gergonne, and a major difficulty was the search for the "best" coordinate
system, with regard to a specific type of problems. To demonstrate the
analytic approach Gergonne asked for proofs of classical problems in synthetic
geometry, such as the famous Appolonius problem: the construction of a circle
tangent to three given ones. His own solution became known as Gergonne's
construction, and later many new elegant solutions, analytic or purely
synthetic, were also found by Poisson, Pl\"{u}cker, Chasles, Poncelet,
Steiner, and others.

The analytic approach to projective geometry became, in fact, a major
geometric discipline during the first decades of the 19th century. Up to 1830
or so, the coordinates being used had typical geometric interpretations such
as length, angle, area, or volume, but this pattern was broken and a new era
began with the introduction of homogeneous coordinates. This also opened the
gate to algebraic geometry and was a major step towards the complete
arithmetization of the geometry, as seen from a modern viewpoint. It is
remarkable that the idea of homogeneous coordinates was discovered almost
simultaneously around 1827 and independently by four geometers, namely
Bobillier, M\"{o}bius, Feuerbach, and Pl\"{u}cker. However, Bobillier and
Feuerbach really did not have the coordinate concept in mind. Let us have a
closer look at some of these early developments.

The early stage of the algebraization was closely related to the use of
\emph{abridged notation,} a method\emph{\ }which, in fact, seeks to avoid the
use of coordinates, elimination of variables, and related messy calculations.
The idea is to represent algebraic equations or their loci by single symbols
$C,$ $C^{\prime}$ etc. or $C=0,C^{\prime}=0$ etc., and apply algebraic
operations to the symbols themselves rather than the coordinates. For example,
it follows from the identity
\begin{equation}
(C_{1}-C_{2})+(C_{2}-C_{3})+(C_{3}-C_{1})=0\text{\ } \label{concur}%
\end{equation}
that the three angular bisectors of a triangle are concurrent, and similarly
the three chords each of which is common to two of three given circles, are
concurrent. One verifies this by writing the lines (resp. circles) on the
normal form $C=px+qy+r$, with $p^{2}+q^{2}=1$ (resp. $C=x^{2}+y^{2}%
+ax+by+c=0)$. This simple proof is also found in many textbooks nowadays.

Applications of this kind, pioneered by Gergonne and Lam\'{e} in 1816, began
to appear in Gergonne's Annales and other journals, as systems of curves or
surfaces became a subject of study. Consider, for example, the widely used
Gergonne's "lambdalization", which for given curves (resp. surfaces) $E=0$,
$E^{\prime}=0$ of degree $n$ constructs the one-parameter family (called
pencil)
\begin{equation}
E+\lambda E^{\prime}=0 \label{lambda}%
\end{equation}
which represents curves (resp. surfaces) of degree $n$ passing through the
intersection of the two given loci.

The foremost early user of abridged notation was E. Bobillier (1798-1840), who
explained the method and published extensive applications in 1827-28. As an
illustrating example, given the three edges $C_{i}=0$ of a triangle, he
considered the two families of equations%

\[
\text{(i) \ }aC_{1}+bC_{2}+cC_{3}=0,\ \ \text{(ii) \ }aC_{1}C_{2}+bC_{2}%
C_{3}+cC_{3}C_{1}=0
\]
with variable parameters $a,b,c$, which in case (i) yield all lines in the
plane, and in case (ii) yield all conics circumscribed about the triangle. In
particular, in case (i) the parameters $a,b,c$ serve as homogeneous line
coordinates in the plane since their mutual ratios determine a unique line.
Recall that M\"{o}bius also studied figures by relating them to a given
triangle, namely he introduced barycentric coordinates relative to the
triangle, but these are homogeneous point coordinates rather than line coordinates.

Next, let us see how K.W. Feuerbach (1800-1834) came up with a system of
homogeneous coordinates, not by using the abridged notation but by the elegant
analytic methods of Lagrange in solid geometry. Feuerbach's investigations of
tetrahedrons (triangular pyramids) in 1827 are similar to what M\"{o}bius did
in the plane the same year, but his approach was purely geometrical rather
than mechanical. He considered a fixed plane $\Pi$ and five generic points
$p_{1},.,p_{5}$ in 3-space, the five tetrahedra $T_{i}$ with the vertices
$p_{j}$ for $j\neq i$, and the five distances $d_{i}$ from $p_{i}$ to $\Pi$.
Then he observed the relation
\[%
{\displaystyle\sum\limits_{i=1}^{5}}
d_{i}Vol(T_{i})=0
\]
where the distances and volumes are signed quantities. Next, by letting the
point $p_{5}$ be a variable point $p$ he introduced the quadruple $(\tau
_{1},\tau_{2},\tau_{3},\tau_{4})$ depending on $p$, where $\tau_{i}%
=Vol(T_{i})$, which yield a kind of homogeneous coordinates for $p$, closely
related to the barycentric coordinates of M\"{o}bius. However, Feuerbach was
not concerned with new coordinates but rather with new theorems in
"tetranometry"; he aimed at expressing all geometric invariants (around 45,
say) of a tetrahedron in terms of six basic ones, such as its edges.

Finally, we turn to Pl\"{u}cker, whose first major work was the two volumes
[1828-1831] on the development of analytic geometry, based on his lectures at
Bonn. Apparently unaware of the works of Lam\'{e}, Gergonne and Bobillier,
Pl\"{u}cker had himself discovered several important aspects of the abridged
notation, which he presented in the first volume and elevated the method to
the status of a principle. Now he became the real expert and made the widest
and most effective use of these ideas. One also speaks of Pl\"{u}cker's
abridged notation, namely when Gergonne's letter $\lambda$ is replaced by
Pl\"{u}cker's $\mu$ in the pencil (\ref{lambda}). For example, an elegant
proof of Pascal's hexagon theorem follows by analysis of the pencil
(\ref{lambda}) with cubic polynomials $E=pqr$, $E^{\prime}=p^{\prime}%
q^{\prime}r^{\prime}$, where $p,q,..,r^{\prime}$ are the six lines of the
hexagon suitably partitioned into two triples as indicated.

The first homogenenous point coordinates proposed by Pl\"{u}cker are the
socalled \emph{trilinear }coordinates. Namely, for a fixed reference triangle,
a point $p$ is assigned the ordered triple $(d_{1},d_{2},d_{3})$ of signed
distances from $p$ to the side lines. Fixed multiples $k_{i}d_{i}$ of the
functions $d_{i}$ are still referred to as trilinear coordinates, and the
barycentric coordinates of M\"{o}bius are recovered as a special case.
Henceforth, we shall use the notation $(d_{1}:d_{2}:d_{3})$ to stress the
homogeneity of the coordinates, that is, they are determined modulo a common
multiple. But Pl\"{u}cker also took the crucial step towards the complete
algebraization of projective geometry by introducing homogeneous coordinates
$(x_{1}:x_{2}:x_{3})$ without any geometric interpretation, namely as the
result of applying any invertible linear substitution to the above coordinates
$d_{i}$.

The projective lines $\simeq P^{1}$ in $P^{2}$ are given by linear equations
$ax_{1}+bx_{2}+cx_{3}=0$, and any of them may serve as the line at infinity
$l_{\infty}$. Henceforth, we shall assume homogeneous coordinates in the plane
chosen so that $l_{\infty}$ is the line $x_{3}=0$:%

\begin{equation}
P^{1}=l_{\infty}:x_{3}=0\text{, \ }(x_{1}:x_{2}:0)\longleftrightarrow
(x_{1}:x_{2}) \label{P1}%
\end{equation}
\
\begin{equation}
A^{2}:x_{3}\neq0\text{, \ }(x_{1}:x_{2}:x_{3})=(\frac{x_{1}}{x_{3}}%
:\frac{x_{2}}{x_{3}}:1\ )\longleftrightarrow(\frac{x_{1}}{x_{3}},\frac{x_{2}%
}{x_{3}})=(x,y) \label{A2}%
\end{equation}
and thus the (ordinary) plane $A^{2}$ is identified with the $xy$-plane. For
example, in M\"{o}bius's barycentric coordinates $(m_{1}:m_{2}:m_{3})$, the
equation of $l_{\infty}$ is $m_{1}+m_{2}+m_{3}=0$, so by setting $x_{1}=m_{1}%
$, $x_{2}=m_{2}$, $x_{3}=m_{1}+m_{2}+m_{3}$, we obtain new coordinates in
terms of which $P^{1}$ and $A^{2}$ are characterized as in (\ref{P1}) and
(\ref{A2}).

Projective coordinates are just homogeneous coordinates naturally rising from
the vector space model of the projective space, as follows. To a given
(n+1)-dimensional vector space $V\simeq\mathbb{R}^{n+1}$ we can associate a
projective space $P^{n}=P(V)$ whose points (resp. lines) are defined to be the
1-dimensional (resp. 2-dimensional) subspaces of $V$, and more generally a
k-plane $\simeq P^{k}$ in $P^{n}$ is a $(k+1)$-dimensional subspace of $V$.
For example, the statement that two points $p_{1},p_{2}$ in $P^{n}$ span a
unique line $l\simeq P^{1}$ is the statement that two distinct central lines
in $V$ span a unique 2-dimensional subspace $\Lambda$. In particular, a vector
$\mathbf{v}\neq0$ and any of its nonzero multiples represent a single point
$p=[\mathbf{v}]$ in $P^{n}$. By choosing a basis $\mathbf{v}_{1}%
\mathbf{,v}_{2}\mathbf{,..,v}_{n+1}$ for $V$, every nonzero vector
$\mathbf{v}=%
{\displaystyle\sum}
x_{i}\mathbf{v}_{i}$ is assigned a nonzero coordinate (row) vector
\begin{equation}
\mathbf{x}=(x_{1},x_{2},..,x_{n+1})\text{,} \label{coord vector}%
\end{equation}
This procedure provides $P^{n}$ with projective coordinate $x_{i}$, and we
shall write
\begin{equation}
p=(x_{1}:x_{2}:x_{3}:..:x_{n+1})\in P^{n} \label{homo}%
\end{equation}
Conversely, a coordinate polyhedron in $P^{n}$ consists of n+1 vertices
$p_{i}=[\mathbf{v}_{i}],i=1,..n+1$, which are the image of some basis of $V$,
and thus the polyhedron gives rise to a projective coordinate system as above.

For example, in projective 3-space, four points $p_{1},..,p_{4}$ in general
position span a coordinate tetrahedron. With respect to the associated
coordinate system the vertices $p_{i}$ of the tetrahedron play the role of the
"standard basis"
\begin{equation}
p_{1}=(1:0:0:0),\text{ \ }p_{2}=(0:1:0:0)\text{, \ }p_{3}=(0:0:1:0)\text{,
\ }p_{4}=(0:0:0:1) \label{vertice}%
\end{equation}
and the opposite face of $p_{i}$ is the plane $(x_{i}=0)$. Then, if we return
to the classical picture (\ref{extensions}) and identifies $A^{3}%
\simeq\mathbb{R}^{3}$ using Cartesian coordinates $(x,y,z)$ as indicated in
(\ref{A2}), the four points in (\ref{vertice}) have the following
interpretation: $\ p_{4}$ is the origin of $\mathbb{R}^{3}$, and $p_{1}%
,p_{2},p_{3}$ are the points at infinity representing the directions of the
three coordinate axes of $\mathbb{R}^{3}.$

In homogeneous coordinates, projective transformations of $P^{n}$ are simply
linear substitutions, namely induced from linear transformations $V\rightarrow
V$ in the above vector space model of $P^{n}$. Therefore, each transformation
$\varphi$ has an invertible matrix $A=(a_{ij})$ of dimension $n+1$, unique up
to a non-zero multiple, so that%
\begin{equation}
\varphi=[A]:(x_{1}:x_{2}..:x_{n+1})\longmapsto(x_{1}^{\prime}:x_{2}^{\prime
}..:x_{n+1}^{\prime})\text{, \ }x_{i}^{\prime}=%
{\displaystyle\sum_{i,j=1}^{n+1}}
a_{ij}x_{j} \label{trans}%
\end{equation}
In the language of groups, the set of these matrices $A$ is the matrix group
usually denoted by $GL(n+1)$, and the non-zero multiples $kI$ of the identity
$I$ is the subgroup $Z$ of all matrices $A$ which yield the identity
transformation $[A]=Id$ of $P^{n}$. Therefore, the full group of projective
transformations can be expressed as the quotient group
\begin{equation}
G(P^{n})=PGL(n+1)=GL(n+1)/Z \label{proj}%
\end{equation}

In the case $n=2$, by restricting the transformation (\ref{trans}) to the
affine subspace $A^{2}\simeq\mathbb{R}^{2}$ with coordinates as in (\ref{A2}),
we calculate
\begin{equation}
(x,y)\longrightarrow(x:y:1)\longrightarrow(x_{1}^{\prime}:x_{2}^{\prime}%
:x_{3}^{\prime})=(\frac{x_{1}^{\prime}}{x_{3}^{\prime}}:\frac{x_{2}^{\prime}%
}{x_{3}^{\prime}}:1)\longrightarrow(x^{\prime},y^{\prime}) \label{trans2}%
\end{equation}
where $(x^{\prime},y^{\prime})$ are exactly the rational expressions in
(\ref{projective}). Moreover, the line $l_{A}$ in (\ref{lines}) expresses as
$x_{3}^{\prime}=0$ and is therefore mapped to the line $l_{\infty}$ at
infinity, whereas the line $l_{B}$ in (\ref{lines}) is mapped to $l_{\infty} $
by the inverse transformation, whose matrix is $B=A^{-1}$.

\ In the projective plane $P^{2}$ with homogeneous coordinates $x_{i}$, a
second order curve (or conic) is the zero set of a quadratic form
\begin{equation}
F(\mathbf{x,x})=\mathbf{x}A\mathbf{x}^{T}=%
{\displaystyle\sum\limits_{i,j=1}^{3}}
a_{ij}x_{i}x_{j}=0 \label{conic}%
\end{equation}
where $A=(a_{ij})$ is a symmetric matrix of dimension 3. Since all matrices of
this kind constitute a 6-dimensional vector space $W$, and matrices differing
by a nonzero scalar represent the same conic, the set of conics naturally
identifies with the projective 5-space $P^{5}=P(W).$

\subsubsection{M\"{o}bius and his approach to duality}

Duality in the projective plane say, is realized by constructing a polar
transformation $\pi$ which is a one-to-one correspondence between points
(pole) and lines (polar) respecting the incidence relation. For this purpose,
pioneers such as Poncelet made use of a chosen conic (\ref{conic}), and their
construction is purely geometric.\ For example, let us choose a circle and
determine the polar of a point $p$ outside the circle. There are exactly two
tangents of the circle passing through $p$, say tangents at the points $p_{1}$
and $p_{2}$, and the polar $\pi(p)$ is the line joining $p_{1}$ and $p_{2}$.

The approach discovered independently by M\"{o}bius (1827) and Pl\"{u}cker
(1829), establishes the analytic counterpart of the geometric principle of
duality, in a new and simple algebraic way, as follows. When a line $l$ \ in
$P^{2}$ is given by a homogeneous linear equation
\begin{equation}
\mathbf{a\cdot x}=a_{1}x_{1}+a_{2}x_{2}+a_{3}x_{3}=0, \label{linear}%
\end{equation}
the three coefficients $a_{i}$ may be regarded as homogeneous \emph{line
coordinates} for $l$, and we shall write $[a_{1}:a_{2}:a_{3}]$ to distinguish
them from point coordinates. Then the equation (\ref{linear}) also provides a
one-to-one correspondence, namely a duality between lines and points%
\begin{equation}
\lbrack a_{1}:a_{2}:a_{3}]\longleftrightarrow(a_{1}:a_{2}:a_{3}). \label{dual}%
\end{equation}
On the one hand, the equation (\ref{linear}) determines all points
$p=(x_{1}:x_{2}:x_{3})$ on the line $l$, but on the other hand, it also
determines the pencil of lines $[a_{1}:a_{2}:a_{3}]$ passing through a fixed
point $p$. We mention that modern mathematicians rather regard $a_{1}%
,a_{2},a_{3}$ as homogeneous coordinates of another projective plane, namely
the dual $P^{2\ast}$ of our plane $P^{2}$. Anyhow, we shall regard
$\mathbf{a}=(a_{1},a_{2},a_{3})$ as the associated coordinate vector
(\ref{coord vector}) of a point in a projective plane.

More generally, let us choose a nondegenerate conic (\ref{conic}), namely the
matrix $A$ is invertible. If points and lines are represented by column
vectors $\mathbf{x},\mathbf{a}$, respectively, the corresponding polar
transformation can be expressed neatly by matrix multiplication as
\begin{equation}
\pi:\mathbf{x\longrightarrow}A\mathbf{x=a,}\text{ \ or }%
\mathbf{a\longrightarrow A}^{-1}\mathbf{a}=\mathbf{x}\text{\ \ }
\label{duality}%
\end{equation}
Two points $p,q$ are said to be \emph{conjugate} if one point lies on the
polar of the other, and this symmetric relation expresses as $A\mathbf{x\cdot
y}=0$ in terms of the coordinate vectors $\mathbf{x}$, $\mathbf{y}$ of the
points. In particular, the self-conjugate points $p$ are those lying on the
conic (\ref{conic}). In view of this, however, the polar map (\ref{dual}) is
rather remarkable, since by (\ref{linear}) and (\ref{dual}) the conjugacy
between $p$ and $q$ simply expresses as $\mathbf{x\cdot y}=0$. This amounts to
using the matrix $A=I$ associated with the conic
\begin{equation}
x_{1}^{2}+x_{2}^{2}+x_{3}^{2}=0 \label{circle0}%
\end{equation}
whose real locus is the empty set. Thus there are no self-conjugate points at
all, but years later the analyst M\"{o}bius identified the locus with an
imaginary circle in the complex extension $\mathbb{C}P^{2}$ of the projective plane.

We have seen how the duality principle can be realized by constructing polar
maps using symmetric invertible matrices $A$, as in (\ref{duality}), and this
works also in higher dimensions. The geometric construction makes use of the
associated quadric hypersurface $A\mathbf{x}\cdot\mathbf{x}=0$ in $P^{n}$. In
his study of geometric mechanics M\"{o}bius made another remarkable discovery,
namely that a skew-symmetric matrix $A$ can also be used to construct a
duality, via the equation $A\mathbf{x\cdot y}=0$ as before. However, in this
case there is no underlying quadric, since $A\mathbf{x}\cdot\mathbf{x}=0$
holds for all $\mathbf{x}$ and so all points are self-conjugate. On the other
hand, invertible skew-symmetric matrices exist only in even dimensions,
$n+1=2k$, so there will be a duality of this kind in $P^{3}$ but not in
$P^{2}$. Thus, M\"{o}bius also settled a dispute between Gergonne and Poncelet
on the nature of dualities in the projective plane; now it turned out that all
of them are actually associated with conics and symmetric matrices as in
(\ref{conic}) (cf. also Gray [2007: 151]).

Pl\"{u}cker presented his work on reciprocity and homogeneous (line)
coordinates in the second volume of [1828-31]. Around 1830 he was maybe not
aware of the work of M\"{o}bius, but proper credit must be given to
M\"{o}bius, whose geometric results were largely overtaken and extended by
Pl\"{u}cker in the 1830's. This is also the decade that Pl\"{u}cker devoted
primarily to an indepth study of algebraic curves in the plane, the major
topic of his books [1835], [1839].

\subsubsection{Classical algebraic geometry and Pl\"{u}cker's formulae for
plane curves}

Modern algebraic geometry arose from the classical studies of curves and
surfaces in the Cartesian plane or space, defined by algebraic (polynomial)
equations, and besides classification of these objects a central topic has
been the behavior of their intersections. For this purpose many concepts and
numerical invariants have been gradually introduced to describe and
distinguish the various types and their possible singularities. The following
review of this topic is largely centered around the achievements of
Pl\"{u}cker during the 1830's, with Boyer\thinspace\lbrack1956] and
Gray\thinspace\lbrack2007] as our major references.

Important properties common to curves in the plane were discovered by Decartes
and Newton, but the scientific foundations of the theory of plane curves are
due to Euler and G.\thinspace Cramer (1704--1752) around 1750. A
classification of the algebraic curves was attempted by Euler, who
distinguished them from the transcendental ones, and Cramer's initial study of
their singularities\ was continued in a modern geometric sense by Poncelet. A
very basic result from the 18th century is the celebrated B\'{e}zout's
theorem, claiming that the number of common points to two plane algebraic
curves with no common component is equal to the product of their degrees. In
fact, special cases such as the intersection of lines, conics, and cubics,
were known already in the 17th century.

E. B\'{e}zout (1730--1783) published his theorem in 1776, based upon heuristic
reasoning and cumbersome calculations, but precise conditions for the theorem
to hold were not formulated. In general, for the validity of B\'{e}zout's
theorem and its generalization to higher dimensions, imaginary points and
points at infinity must also be considered. But the most delicate part is the
assigning of proper intersection multiplicities. We shall not describe this
procedure, but we mention that a common tangent point of two plane curves has
intersection multiplicity at least two.

In the sequel we shall recall some of the classical theory of plane curves,
due to Pl\"{u}cker and others. Assuming projective coordinates $(x:y:z)$, say
$z=0$ is the line at infinity, an algebraic curve $C_{n}$ of degree (or order)
$n$ is the locus of a polynomial equation
\begin{equation}
C_{n}:f(x,y,z)=%
{\displaystyle\sum\limits_{i+j+k=n}}
q_{i,j,k}x^{i}y^{j}z^{k}=0 \label{curve}%
\end{equation}
and so its equation in the Cartesian $xy$-plane is found by setting $z=1$.
Clearly, the latter equation still has degree $n$, unless $f(x,y,z)$ is
reducible with $z$ as a factor. The abridged notation $C_{n}$ for nth order
curves is rather typical in the classical literature, and it is also used in
Klein's letters to Lie. Similarly, $F_{n}$ is used to denote a surface
(Fl\"{a}che) of degree $n$.

The number of coefficients in (\ref{curve}) is $1+n(n+3)/2$, consequently a
given curve $\bar{C}_{n}$ is uniquely determined by%
\begin{equation}
\mu_{n}=\frac{n(n+3)}{2} \label{points}%
\end{equation}
suitably chosen points. Let $C_{n}$ and $C_{n}^{\prime}$ be distinct curves
having $\mu_{n}-1$ of these points, but missing the last point $\bar{p}$, say.
Then a principle due to Lam\'{e} states that the equations
\begin{equation}
C_{n}+\mu C_{n}^{\prime}=0 \label{Lame}%
\end{equation}
represent all curves of degree $n$ passing through the $\mu_{n}-1$ given
points. Clearly, for a specific value of $\mu$ we recover the given curve
$\bar{C}_{n}$.

On the other hand, Cramer and Euler had observed the paradoxical fact, that
although a cubic curve $C_{3}$ in general is uniquely determined by $\mu
_{3}=9$ points, two cubic curves still have $n^{3}=9$ points in common, by
B\'{e}zout's theorem. They both realized that, somehow, interdependence of
points was involved, and this puzzle became known as Cramer's paradox. It was
first resolved by Pl\"{u}cker, who gave a clearer answer which, in fact,
follows from the simple observation that the curves $C_{n}$, $C_{n}^{\prime} $
in (\ref{Lame}) with $\mu_{n}-1$ given points in common, actually have $\ $%
\begin{equation}
\delta_{n}=n^{2}-(\mu_{n}-1)=\frac{1}{2}(n-1)(n-2) \label{delta}%
\end{equation}
additional points in common, by B\'{e}zout's theorem. Therefore, the whole
curve family (\ref{Lame}) have the $n^{2}$ points in common. For example, in
the case of quartic curves, they will have 3 extra common points which are
determined by the 13 common given points. Similar explanations of Cramer's
paradox were also given by Gergonne, Jacobi, and Lam\'{e}.

In Pl\"{u}cker's analysis of plane curves the principle of duality plays a
crucial role, since any curve $C=C_{n}$ has a dual curve $C^{\ast}%
=C_{m}^{\prime}$ of some degree $m$, called the \emph{class} of $C$. An
alternative definition of the class $m$ was given by Gergonne in 1826, namely
it is the number of tangents to $C$ passing through a fixed (generic) point.
The two definitions would seem to be equivalent, since lines and points are
dual to each other and the fixed point becomes a line cutting $C^{\ast}$ in
$m$ points. So, by B\'{e}zout's theorem $m$ is also the degree of $C^{\ast}$.

Now, Gergonne mistakenly assumed $n=m$, despite the fact that Monge had
earlier estimated the number of tangents to be $n(n-1)$, in a theorem which
had been overlooked. Thus, Monge's result seemed to imply $m=n(n-1)$ for all
$n$; in fact, Poncelet confirmed that a curve of degree $n$ is generally of
class $m=n(n-1).$ But on the other hand, the dual of $C^{\ast}$ is
$C^{\ast\ast}=C$, which clearly leads to a contradiction when $n\neq2$. This
explains the so-called \emph{Duality paradox}, which was still unsettled in
the early 1830's. But Pl\"{u}cker and others must have realized that the
correct value of $m$ is also sensitive to the behavior of $C$ and $C^{\ast}$
at their \emph{singular} points. Since he was using line coordinates as well
as point coordinates, Pl\"{u}cker was in a better position to resolve the
paradox by simultaneously analyzing the curve and its dual curve.

The advantage of developing point and line conceptions simultaneously had also
been noted by Brianchon and Poncelet. In fact, the notion of a curve as the
envelope of its tangent lines had been proposed already in the 17th century,
and Leibniz (1692) gave rules for calculating envelopes. However, it was
M\"{o}bius (1827) who determined the condition%
\begin{equation}
f^{\ast}(u,v,w)=0 \label{dual2}%
\end{equation}
for a line $ux+vy+wz=0$ to be tangent to the curve $C_{n}$ (\ref{curve}). But
he did not express the idea that $[u:v:w]$ are coordinates of the line, nor
did he associate (\ref{dual2}) with the dual curve and hence overlooked the
possibility of finding the class $m$ of $C_{n}$ by calculation of the degree
of the curve (\ref{dual2}). This discovery is essentially due to Pl\"{u}cker,
who clearly understood that all curves, except points and lines, have both
point equations and line equations. In a paper published in Crelle's Journal
(1830) he made the prescient remark that the general theory of curves should
be developed together with the idea of singular points and singular tangents.

With his new insight Pl\"{u}cker clarified the correspondence between tangent
singularities of a curve $C$ and point singularities of its dual curve
$C^{\ast}$, and vice versa. This enabled him to set forth the celebrated
formulae (\ref{formulae}), from which a resolution of the duality paradox
follows as a simple consequence. He communicated the formulae in the first
place to Crelle's Journal (1834); a sketchy version of his theory appeared in
[1835] and a complete account and further extension in [1839]. Pl\"{u}cker
considered two types of point singularities and two types of tangent
singularities, namely
\begin{align}
\text{double point }(\delta)\text{ }  &  \longleftrightarrow\text{\ bitangent
}(\tau)\label{sing}\\
\text{cusp }(\kappa)  &  \longleftrightarrow\text{stationary tangent }%
(\iota)\nonumber
\end{align}
where the letters $\delta$, $\kappa$, $\tau$, $\iota$ count the number of each
kind. Alternative terms used in the literature are node, double tangent,
triple point, and inflectional tangent, respectively. For example, a double
point (i) and a cusp (ii) at the origin is illustrated by the cubic curves:%
\[
(i)\text{ }x^{3}+y^{3}-5xy=0\text{, \ \ }(ii)\text{ }x^{2}-y^{3}=0.
\]
\ 

As indiated in (\ref{sing}), double points and bitangents are dual to each
other, and similarly cusps and stationary tangents are dual. That is, via the
duality map $C\rightarrow C^{\ast}$ a double point of $C$ becomes a bitangent
of the dual curve, and conversely. Hence, if the symbols $\delta^{\ast}%
,\tau^{\ast}$ etc. count the singularities of $C^{\ast}$, then $\delta^{\ast
}=\tau$, $\tau^{\ast}=\delta$, $\kappa^{\ast}=\iota$, $\kappa=\iota^{\ast}$.
The six numerical invariants associated with the plane curve $C_{n}$ are
constrained by Pl\"{u}cker 's formulae, which can be stated as follows:%

\begin{align}
\text{\ }m  &  =n(n-1)-2\delta-3\kappa\text{, \ \ \ \ }n=m(m-1)-2\tau
-3\iota\text{\ }\label{formulae}\\
i  &  =3n(n-2)-6\delta-8\kappa\text{, \ \ }\kappa=3m(m-2)-6\tau-8\iota
\nonumber
\end{align}

As expected, a conic ($n=2$) has no singular points or tangents, and its dual
is still a conic. A curve is said to be \emph{non-singular} if $\delta
=\kappa=0$, but we note that such a curve has necessarily singular tangents
(if $n>2)$, and consequently its dual curve is singular. Point singularities
had been studied at least in the previous 200 years, and it was known that the
number of such points is limited by the degree $n$ of the curve. For example,
C. MacLaurin (1698-1746) showed the number of double points is limited by the
constant (\ref{delta}). The finiteness of the number of bitangents was first
suggested by Poncelet (1832), but shortly afterwards Pl\"{u}cker established
an upper bound $\tau_{n}$, namely we have
\begin{equation}
\delta\leq\delta_{n}=\frac{1}{2}(n-1)(n-2)\text{, \ \ \ }\tau\leq\tau
_{n}=\frac{1}{2}n(n-2)(n^{2}-9) \label{limits}%
\end{equation}

B\'{e}zout's theorem is the basic tool for the estimation of the above
numerical invariants, namely the singular points in question should be
recovered by intersecting the given curve $C_{n}$ with some suitably related
curve $C^{\prime}$, for example the Hessian of $C_{n}$. Its equation is
$H(f)=0$, where $H(f)$ is the Hessian determinant of $f$ (\ref{curve}), which
has degree $3(n-2).$ Hence, the curve and its Hessian have $3n(n-2)$ common
points, and they include the inflection points of $C_{n}$, that is, the points
where the curvature vanishes. Consequently,%

\[
\iota\leq\iota_{n}=3n(n-2)
\]
and moreover, equality holds if $C_{n}$ is non-singular. In general, the
Hessian also passes through any double point or cusp, in fact, it has 6-fold
and 8-fold contact with such points, respectively. So, this reduces the
maximal value $\iota_{n}$ of $\iota$ by $6\delta+8\kappa$ when the curve is
singular, in agreement with (\ref{formulae}).

Pl\"{u}cker estimated the number of bitangents of a non-singular curve $C_{n}
$ to be the upper bound $\tau_{n}$ in (\ref{limits}), as follows. Since the
dual curve $C^{\ast}$ has degree $m=n(n-1)$, and the dual of $C^{\ast}$ has
degree $n$, one of Pl\"{u}cker's formulae amounts to the equation
\[
n=n(n-1)(n(n-1)-1)-2\tau-3\iota_{n}%
\]
and as an equation for $\tau$ this\ has the unique solution $\tau=\tau_{n} $.
However, he did not give any proof of this independent of his formulae
(\ref{formulae}). Such a proof was given by Jacobi (1850), which in addition
to several other results confirmed the validity of Pl\"{u}cker's approach.

Pl\"{u}cker's formulae provided an effective tool for the determination of the
possible values of the above numerical invariant, say for $n,m\leq10$. But
first of all, they enabled him to progress more deeply into the study of cubic
and quartic curves. However, Pl\"{u}cker did not proceed to the more general
case of higher orders, where more complicated "non-Pl\"{u}ckerian"
singularities must also be considered, but we shall leave this topic here and
refer to Gray [2007: 169]. In the following decades, however, quite new
techniques of algebraic geometry were gradually developed, such as resolution
of singularities of plane curves as well as space curves. Cremona
transformations were found to be useful for this purpose, in particular, they
are effective in the reduction of singularities of curves to double points
with distinct tangents.

Around 1860 or so, Clebsch came across the following numerical invariant for
plane curves%
\begin{equation}
g=\delta_{n}-\delta-\kappa=\frac{1}{2}(n-1)(n-2)-(\delta+\kappa) \label{genus}%
\end{equation}
which takes the same value for both $C_{n}$ and its dual curve. The term
\emph{deficiency} was originally used, since it counts the maximal number of
double points reduced by the actual number of double points and cusps. The
same invariant was in fact studied by Riemann, who also realized its
importance as a topological invariant. It was renamed the \emph{genus }of the
curve, which is also the modern term.

\subsubsection{\textquotedblleft Projective geometry is all
geometry\textquotedblright\ (Cayley 1858)}

The work of von Staudt (with some later corrections) made it clear that
projective geometry can be built up without dependence on Euclidean metric
concepts such as angle and distance. But conversely, it is also rather
surprising that it is possible to express these metric quantities on the basis
of purely projective concepts. For example, the angle in radians between two
intersecting lines $l_{1}$ and $l_{2}$ in the Euclidean plane can be expressed
in terms of the cross-ratio of four lines%
\begin{equation}
\psi=\measuredangle(l_{1},l_{2})=\frac{i}{2}\log(l_{1},l_{2};\omega_{1}%
,\omega_{2}) \label{Laguerre}%
\end{equation}
where $\omega_{1},\omega_{2}$ are the two (imaginary) lines passing through
the vertex of the angle and the two circular points at infinity, $(1:i:0)$ and
$(1:-i:0)$ respectively, see equation (i) in (\ref{circular}). Here we regard
(as usual) the Euclidean plane as part of the projective plane, which in turn
extends to the complex projective plane, namely $E^{2}\subset P^{2}%
\subset\mathbb{C}P^{2}$. The coefficient $i/2$ is needed to make the angle
real and to ensure a right angle has value $\pi/2$. This formula was first
discovered by Laguerre in 1853. It does not seem that he was looking for a
similar formula for the distance between two points.

Independent of Laguerre, Cayley also wanted to show that angle and distance,
namely metric notions in Euclidean geometr, can be formulated in projective
terms. For this purpose he introduced a conic (resp. a quadric) in the case of
plane (resp. space) geometry, which he referred to as the \emph{absolute }
figure. As pointed out in his Sixth Memoir upon Quantics [1858], the two cases
are similar, so let us follow Cayley's approach in the projective plane, with
homogeneous coordinates $(x_{1}:x_{2}:x_{3})$ and an absolute conic
(\ref{conic}) with associated real bilinear form
\begin{equation}
F(\mathbf{x},\mathbf{y)=}\sum_{1}^{3}a_{ij}x_{i}y_{j} \label{bilinear}%
\end{equation}
In terms of line coordinates $[u_{1}:u_{2}:u_{3}]$ (see (\ref{linear})), the
equation $F^{\ast}(\mathbf{u,u})=0$ also describes the above conic, where
$F^{\ast}(\mathbf{u,v})$ is the bilinear form whose matrix $A^{\ast}$ is the
cofactor matrix (or adjoint) of the real symmetric matrix $A=(a_{ij}).$ \ 

Now, Cayley defined the distance between $\mathbf{x}$ and $\mathbf{y}$ by%
\begin{equation}
\delta(\mathbf{x,y)=\arccos}\frac{F(\mathbf{x,y})}{\sqrt{F(\mathbf{x,x})}%
\sqrt{F(\mathbf{y,y})}} \label{distance}%
\end{equation}
and for the angle $\psi=\measuredangle(\mathbf{u},\mathbf{v})$ between lines
$\mathbf{u}$ and $\mathbf{v}$ he defined $\cos\psi$ in the same way using
$F^{\ast}$ instead of $F$.

Depending on the absolute, $\delta$ is a generalization of the modern notion
of a distance function. The projective space has lines, and for three points
$\mathbf{x,y,z}$ on a line the function is additive%
\begin{equation}
\delta(\mathbf{x,y)+}\delta(\mathbf{y,z)}=\delta(\mathbf{x,z)} \label{add}%
\end{equation}
Notice that the value of the expression (\ref{distance}) lies in the interval
$[0,\pi]$, and the formula for $\delta$ only gives the distance modulo $\pi$.
Hence, starting from $\mathbf{x}$ and moving along a line towards $\mathbf{y}%
$, we see from (\ref{add}) that multiples of $\pi$ may accumulate before
$\mathbf{y}$ is reached. In fact, in 1870 Klein modified Cayley's definition
(\ref{distance}) of distance and removed the use of the $\arccos$ function, by
expressing the distance between two points as a line linegral.

Familiar expressions for $\delta$ are obtained, for example when $A=I$
(identity). But similar to Laguerre, Cayley's absolute could also be the two
circular points at infinity, viewed as a degenerate conic, and then he came
out with the expression
\begin{equation}
\delta(\mathbf{x,y)}^{2}=\frac{(x_{1}y_{3}-y_{1}x_{3})^{2}+(x_{2}y_{3}%
-y_{2}x_{3})^{2}}{x_{3}^{2}y_{3}^{2}} \label{distance1}%
\end{equation}
It follows that $\delta$ restricted to the affine plane $A^{2}=(x_{3}%
=y_{3}=1)$is, in fact, the usual Euclidean distance function, so the induced
geometry on $A^{2}$ is that of the Euclidean plane $E^{2}.$

With his paper [1858] Cayley had reduced metric geometry to projective
geometry, but Cayley himself only showed how Euclidean geometry can be
reinterpreted in terms of projective properties. In 1870 Klein suggested the
non-Euclidean geometries are related to the projective metric as well, by an
appropriate choice of the absolute figure as a standard of reference, and with
the paper [1871g] and its second part the following year he set forth his new
ideas. Instead of his previous integral expressions Klein defined the distance
between $\mathbf{x}$ and $\mathbf{y}$, and the angle between $\mathbf{u}$ and
$\mathbf{v}$, by taking the logarithm of a cross-ratio, namely
\begin{equation}
\delta(\mathbf{x,y)}=c\log(\mathbf{x,y;a,b})\text{, \ \ }\psi=c^{\prime}%
\log(\mathbf{u,v;w}_{1}\mathbf{,w}_{2}) \label{distance2}%
\end{equation}
where $\mathbf{a}$, $\mathbf{b}$ denote the two intersection points of the
absolute conic (or quadric) and the line through $\mathbf{x}$ and $\mathbf{y}
$, $\mathbf{w}_{1}$ and $\mathbf{w}_{2}$ are the tangent lines of the conic
through the intersection point, and $c$ and $c^{\prime}$ are appropriate constants.

In plane geometry, Klein showed the geometry is hyperbolic, spherical, or
Euclidean, according to whether the conic is real, imaginary, or degenerate.
For example, the imaginary circle (\ref{circle0}) yields elliptic geometry.
Beltrami's disk model for hyperbolic geometry (see Section 4.4.3) can be
derived from Klein's procedure by taking the following circle as the absolute
conic%
\begin{equation}
x_{1}^{2}+x_{2}^{2}-x_{3}^{2}=0\text{ \ \ or \ }x^{2}+y^{2}=1\text{ \ \ cf.
}(\ref{A2}) \label{circle2}%
\end{equation}
and the interior of the circle as the hyperbolic plane $H^{2}$. Beltrami did
not derive a distance function or formula for angles in this disk. However, in
general one can calculate the cross-ratio expressions in (\ref{distance2}) as
\begin{equation}
(\mathbf{x,y;a,b})=\frac{F(\mathbf{x,y)+}\sqrt{F(\mathbf{x,y)}^{2}%
-F(\mathbf{x,x)}F(\mathbf{y,y)}}}{F(\mathbf{x,y)-}\sqrt{F(\mathbf{x,y)}%
^{2}-F(\mathbf{x,x)}F(\mathbf{y,y)}}} \label{cross1}%
\end{equation}
and similarly for the cross-ratio of lines using $F^{\ast}$, cf. also Kline
[1972: 911]. In the case (\ref{circle2}) we have, of course, $F=F^{\ast}$ and
$\ $%
\begin{equation}
F(\mathbf{x,x)}=x_{1}^{2}+x_{2}^{2}-x_{3}^{2} \label{form}%
\end{equation}

Some writers on the topic (wrongfully) attribute Klein's formulae
(\ref{distance2}) to Cayley, and in fact they turn out to be just a
reformulation of Cayley's formulae. In his Collected Math. Papers, Cayley
added important comments to his Sixth Memoir article [1858] and its connection
with Klein's work, stating that he regarded Klein's approach as an improvement
of his. However, he expressed disagreement with Klein's interpretation of
non-Euclidean geometry and its relation to the projective metric.

But Cayley neglected exploring his metric to the fullest; in particular, he
did not relate it to non-Euclidean geometry, perhaps because of his rather
ambivalent attitude toward this geometry. Instead, Klein discovered many
different metric subgeometries of projective geometry, depending on the choice
of the absolute, a valuable experience which contributed to the formulation of
his Erlanger Programm in 1872.

\ \ \ \ \ \ \ \ \ 

{\large Remarks on groups}{\Large \ }The Cayley-Klein approach also describes
the relationship between motions (isometries) of the metric subgeometry and
projective transformations of the ambient space. The general principle, which
also appeared in Klein's Erlanger Programm, is that the motions are the
restriction of those projective transformations leaving the absolute figure
invariant, not necessarily pointwise. Assuming this, let us determine the
isometry group of the hyperbolic plane, by calculating the group of projective
transformations of $P^{2}$ leaving the circle (\ref{circle2}) invariant. In
the notation of (\ref{proj}), the matrix subgroup of $GL(3)$ leaving the
bilinear form $F$ invariant is usually denoted $O(2,1) $, and extension of
this by the scaling group $Z$ yields the matrix group leaving the form
invariant modulo scaling, consequently
\begin{equation}
Iso(H^{2})=\frac{O(2,1)\cdot Z}{Z}=\frac{O(2,1)}{\left\{  \pm I\right\}
}\simeq SO(2,1) \label{Isohyp}%
\end{equation}

\subsection{Line geometry}

Line geometry is an approach to geometry where geometric objects in projective
3-space $P^{3}$ (or n-space in general) are studied by considering the
straight lines, rather than the points, as the basic geometric \emph{Elements}%
. Thus the geometric objects\ are represented by appropriate configurations of
lines, and so line geometry becomes a branch of projective geometry. It was
initiated in the early 19th century, by mathematicians such as Monge, Malus,
M\"{o}bius, and Hamilton, who studied families of lines in 3-space and their
geometric properties, often motivated by studies and experiments in optics. We
shall use the terms "ray" and "line" synonymously.

Although line geometry went out of fashion in the early 20th century, the
study was revived with modern techniques at the end of the century; for
example, we refer to the survey on low order congruences in Arrondo [2002],
and to Pottman-Wallner [2001] on computational line geometry and its modern
applications. In the sequel we shall assume the underlying 3-space is the
complex projective space $\mathbb{C}P^{3}$, in the tradition of Pl\"{u}cker,
Kummer, Klein and Lie. For surveys of the classical line geometry, see
Lie-Scheffers [1896], Jessop [1903], Rowe [1989].

\subsubsection{Ray systems and focal surfaces}

Recall that rays parallel to the axis of a parabolic surface are reflected
into rays passing through the focal point. In general, however, light rays
will not meet at a single focus after reflection or refraction, so there will
be overlapping and the rays envelop and thus create an interesting geometric
pattern, such as a\emph{\ caustic curve} or \emph{caustic surface. }An example
is provided by a \emph{normal congruence}, namely\ the lines perpendicular to
a given surface. In these examples the rays constitute a 2-parameter family of
rays called a \emph{ray system}. Caustics were maybe first introduced by E.L.
Malus (1775-1812), a student of Monge and Fourier, who began publishing papers
on optics at \'{E}cole Polythechnique in 1808. It was Pl\"{u}cker who
introduced the modern term \emph{line congruence} for a ray system. The
caustic surface became known as the \emph{Brennfl\"{a}che} or \emph{focal
surface} in the German or English literature, respectively. Generally the
surface has two components and the line congruence consists of their common
tangent lines.

In the 1840's Pl\"{u}cker was a leading figure in line geometry, and it was
after his renewal of the theory in the 1860's that the discipline became, in
fact, a major topic in algebraic geometry. Pl\"{u}cker's approach provided, in
fact, much of the geometric framework for the studies of Klein and Lie during
the first years of their career. But their interests were also greatly
stimulated by recent results of Kummer on specific surfaces arising as the
focal surface of certain line congruences. These surfaces have been known as
\emph{Kummer surfaces} ever since. Algebraically, they are quartic surfaces in
projective 3-space with 16 double points, and it seems that they arose from
his interests in Dupin's cyclides and the optical properties of biaxial
crystals. Needless to say, the relationship between the Kummer surfaces, the
theta-function, and quotients of abelian surfaces, was only perceived much
later in the development of algebraic geometry.

As a leading algebraic number theorist in 1855, Kummer became a professor at
the University in Berlin, and in 1857 he was awarded the Grand Prix at the
Academy of Science in Paris for his fundamental work relating to Fermat's last
theorem. At the end of the decade, however, his interests drifted towards
geometry and he became interested in the ray systems examined by Hamilton. In
fact, the first fundamental paper on line congruences was written by Kummer in
1859, and published in Crelle's Journal the next year. He introduced, for
example, the notion of a density function, which for a normal congruence (see
above) equals the Gaussian curvature of the underlying surface. Kummer's
geometric period lasted (at least) through the 1860's, and at his seminar in
Berlin during the fall 1869, with Klein and Lie among the participants, the
topic was in fact line geometry.

\subsubsection{Basic ideas and definitions}

It had been known long before Pl\"{u}cker that a line in 3-space has 4
independent degrees of freedom. In his initial paper [1846] on the subject,
Pl\"{u}cker described the family of all lines in xyz-space in terms of four
line coordinates $(r,s,\rho,\sigma)$ by the equations%
\[
x=rz+\ \rho,\text{ \ }y=sz+\sigma
\]
Although this simple and naive 4-parameter representation has singularities
and exceptional lines (which can be avoided), this did not prevent Pl\"{u}cker
from developing many of the basic concepts and properties of line geometry. He
defined a \textit{line complex }to be a 3-parameter family of lines given by
an algebraic equation $F(r,s,\rho,\sigma)=0$. However, this approach
encountered difficulties since the degree of the defining equation is not
invariant under linear transformations of the variables $x,y,z$. He solved the
problem in his "English" paper [1865] by introducing the auxiliary line
coordinate $\eta=r\sigma-s\rho$, which enabled him to define the
\textit{order} (or degree) of a line complex $\mathcal{C}$ to be the degree of
the defining equation
\begin{equation}
\mathcal{C}:F(r,s,\rho,\sigma,\eta)=0 \label{complex}%
\end{equation}

For example, a linear complex is a complex of order one. A\emph{\ special}
line complex consists of the lines intersecting a given line or curve, called
the \emph{directrix}, but the complex of lines tangent to a given
(non-planary) surface is also called special. These complexes will be of the
same order as the curve or surface. For example, there is a special linear
complex associated with each line, and these complexes were, in fact, studied
by M\"{o}bius, who called them \emph{null systems.} Here the lines through a
point $p$ lie in a fixed plane $\mathcal{C}_{p}$ and constitute the full
pencil of lines, and there is a dual relationship $p\longleftrightarrow
\mathcal{C}_{p}$ between points and planes called the null polarity. First of
all, however, it was Pl\"{u}cker who called attention to the applications of
line complexes of degree one and two to mechanics and optics.

A 1-parameter family $\mathcal{R}$ of lines represents a ruled surface in
3-space, and a 2-parameter family $\mathcal{K}$ is called a line congruence.
They can be expessed as an appropriate intersection of two or three line
complexes
\begin{equation}
\mathcal{K}=\mathcal{C}_{1}\cap\mathcal{C}_{2},\text{ \ \ \ }\mathcal{R}%
\text{\ }=\mathcal{C}_{1}\cap\mathcal{C}_{2}\cap\mathcal{C}_{3}\text{\ \ }
\label{line-2}%
\end{equation}
The \textit{order} of a line congruence $\mathcal{K}$ is defined to be the
number of lines passing through a general point. In (\ref{line-2}) this is the
product $n_{1}n_{2}$ of the orders of $\mathcal{C}_{1}$ and $\mathcal{C}_{2}$.
Dually, the \emph{class}\textit{\ }of the congruence is the number of its
lines lying in a general plane. The two numbers may be different, but they are
equal if $\mathcal{K}$ belongs to a linear complex, say $\mathcal{C}_{1}$ in
(\ref{line-2}) is linear.

The degree of the ruled surface formed by $\mathcal{R}$ in (\ref{line-2})
equals $2n_{1}n_{2}n_{3}$. In the lowest degree case it is, in fact, a doubly
ruled quadric since through each of its points there pass two distinct lines
lying on the surface. So the surface is a hyperbolic paraboloid or a
hyperboloid of one sheet, see (\ref{hyperboloid}), unless it degenerates to a
plane or two planes.

Recall from Section 3.5.1, a line congruence gives rise to a specific surface,
namely the focal surface enveloped by its lines. In fact, Pl\"{u}cker
associated with a given line complex $\mathcal{C}$ infinitely many surfaces of
this kind - the \emph{complex surfaces -} namely the focal surfaces of all
congruences $\mathcal{C}\cap\mathcal{C}_{1}$, where $\mathcal{C}_{1}$ is the
special linear complex with a given line $l$ as directrix. In his book [1868],
Part III, he classified these surfaces into seven families, for complexes
$\mathcal{C}$ of order 2 (cf. letter of 29 July, 1871).

Pl\"{u}cker also used the terminology \emph{complex lines} and \emph{complex
curves}, respectively, for the lines belonging to\ $\mathcal{C}$ and the
curves whose tangent lines belong to $\mathcal{C}$. Moreover, at each point
$p$ there is the \textit{complex cone} $\mathcal{C}_{p}$, namely the
1-parameter family of all complex lines passing through $p.$Thus, an
alternative description of a complex surface is that it is enveloped by the
cones $\mathcal{C}_{p}$ as $p$ runs through a fixed line $l$. We also mention
that Lie, in his study of tetrahedral line complexes in 1869-70 (cf. [1869],
[1870a]) used his own definition of complex surfaces, and they were
expressible as solutions of a certain differential equation. According to Lie,
a tetrahedral complex consists of those lines whose four intersection points
with the planes of a given (coordinate) tetrahedron have a fixed cross-ratio.
However, based on another definition these line complexes had, in fact, also
been studied by T. Reye in Z\"{u}rich a few years earlier.

The complex cones $\mathcal{C}_{p}$ describe the local geometry of a line
complex of degree $n$. At a non-singular point $p$ the cone $\mathcal{C}_{p}$
is non-degenerate, and it has degree $n$ in the sense that it cuts any plane
$P^{2}$ along an algebraic curve $\gamma=$ $\mathcal{C}_{p}\cap P^{2}$ of
degree $n$. Namely, $\gamma$ intersects a general line in $P^{2}$ at n points.
Dually, associated with the given line complex are also the non-singular
planes $\simeq P^{2}$ in 3-space, with the property that the complex lines in
the plane envelop a curve $\Gamma$ of \textit{class }n. >From algebraic
geometry, the curve is said to be of class $n$ since it has n tangents passing
through an arbitrary point in the plane.

On the other hand, the singular points (or planes) of a line complex
constitute the singular locus, whose geometry is generally complicated and
hard to visualize. But it has attracted some attention in special cases, such
as complexes of order two. These are the quadratic complexes, whose cone
$\mathcal{C}_{p}$ at a singular point \ $p$ degenerates into two planes
intersecting along a \emph{singular line}. Thus, the locus of singular points
is the \emph{singularity surface}, enveloped by the congruence of singular lines.

Dually, one also arrives at the singularity surface by considering the
\emph{singular planes} of the complex. In such a plane the class curve
$\Gamma$ degenerates into two points $p_{1}$ and $p_{2},$ and the complex
lines in the plane are those passing through either $p_{1}$ or $p_{2.}$ The
line joining $p_{1}$ and $p_{2}$ is the singular line, and again the singular
lines form a congruence whose focal surface is the singularity\emph{\ }surface
of the quadratic line complex.\emph{\ }

With the paper [1866] Kummer gave a classification of linear line congruences
(missing one case, see ref. in Arrondo [2002]). But more importantly, he also
made a deeper study of the focal surface associated with a line congruence of
order 2 and class 2, which generally is a surface of order 4 and class 4. In
1864 he had shown that this is a surface in $\mathbb{C}P^{3}$ with16 double
points and 16 double tangent planes, the maximum number possible for quartic
surfaces. These are the \emph{Kummer surfaces}, which due to their nice
properties and dominant role among algebraic surfaces of degree 4, became an
interesting study object for geometers and algebraists in the ensuing years.
For a detailed exposition of the topic,\ see Hudson [1990].

The Kummer surface (and its generalization to higher dimensions) developed
into an intricate speciality among French, German, British, and Italian
geometers well into the 20th century, but the interest declined rapidly in the
1920's. Today these surfaces are regarded as a 3-parameter family which play
an important role in the modern theory of so-called K3 surfaces. For example,
the following algebraic equation
\begin{equation}
x^{4}+y^{4}+z^{4}-(x^{2}y^{2}+x^{2}z^{2}+y^{2}z^{2})-(x^{2}+y^{2}+z^{2})+1=0
\label{Kummer}%
\end{equation}
represents a typical Kummer surface, symmetric with respect to permutations
and change of sign of the coordinates,consequently the surface has the
symmetries of an embedded regular octahedron.

\subsubsection{Pl\"{u}cker's new appoach}

Now, let us return to the late 1860's when Pl\"{u}cker was preparing his
2-volume study \emph{Neue Geometrie des Raumes...} with focus on linear and
quadratic line complexes, which appeared in 1868-69. Grassman and Cayley were
in fact forerunners of the new approach, based upon the use of homogeneous
coordinates and \emph{Elements} from exterior linear algebra.

Consider projective 3-space $P^{3}$ with homogeneous coordinates $x_{i}$
relative to a given coordinate tetradron with vertices $p_{1},..,p_{4}$, as in
(\ref{homo}). This tetrahedron was also referred to as the \emph{fundamental
tetrahedron.} A basic observation is that the line $\overline{\mathbf{xy}}$
between two points $\mathbf{x}$ and $\mathbf{y}$ is determined by the
$2$-minors $p_{ij}$ of the following matrix
\begin{equation}
\ [\mathbf{x},\mathbf{y}]=\left[
\begin{array}
[c]{cccc}%
x_{1} & x_{2} & x_{3} & x_{4}\\
y_{1} & y_{2} & y_{3} & y_{4}%
\end{array}
\right]  ,\text{ \ }p_{ij}=\left\vert
\begin{array}
[c]{cc}%
x_{i} & x_{j}\\
y_{i} & y_{j}%
\end{array}
\right\vert =x_{i}y_{j}-x_{j}y_{i} \label{line-3}%
\end{equation}
Therefore, one can take the 6-tuple of Pl\"{u}cker coordinates
\begin{equation}
\ (p)=(p_{12}:p_{13}:p_{14}:p_{23}:p_{42}:p_{34}) \label{(p)}%
\end{equation}
as homogeneous coordinates for the line $\overline{\mathbf{xy}}$. Note,
however, the numbers $p_{ij}$ are constrained by the \emph{Pl\"{u}cker
relation }
\begin{equation}
\mathcal{P}=p_{12}p_{34}+p_{13}p_{42}+p_{14}p_{23}=0 \label{line-4}%
\end{equation}
and conversely, it is not difficult to see that each 6-tuple $(p)$ satisfying
the condition $\mathcal{P}=0$ does, indeed, represent a line $l$ in $P^{3}$.
Let us also define the \emph{dual }line of $l$ to be the line $l^{\ast}$ with
Pl\"{u}cker coordinates
\[
(p^{\ast})=(p_{12}^{\ast}...:p_{34}^{\ast})=(p_{34}:p_{42}:p_{23}%
:p_{14}:p_{13}:p_{12})
\]

On the other hand, following M\"{o}bius and Pl\"{u}cker (see (\ref{linear})),
in the equation $\mathbf{u\cdot x}=0$ for a plane $U$ in $P^{3}$, the plane
and point coordinates $u_{i}$ and $x_{i}$ appear symmetrically. This suggests
a simple realization of the Poncelet-Gergonne duality principle as the
following correspondence between planes and points%

\begin{equation}
U=[u_{1}:u_{2}:u_{3}:u_{4}]\longleftrightarrow(u_{1}:u_{2}:u_{3}%
:u_{4})=\mathbf{u}\text{, \ cf. }(\ref{dual}) \label{dual1}%
\end{equation}
Let $l=U\cap V$ be the line of intersection of two planes, say $l=$
$\overline{\mathbf{xy}}$ has Pl\"{u}cker coordinates $(p)$ as above. However,
the two points $\mathbf{u}$ and $\mathbf{v}$ dual to $U$ and $V$ also define a
line $l^{\prime}=$ $\overline{\mathbf{uv}}$ where planes $\mathbf{u}$ and
$\mathbf{v}$, say with Pl\"{u}cker coordinates
\[
(q)=(q_{12}:q_{13}...:q_{34})\text{, \ }q_{ij}=u_{i}v_{j}-u_{j}v_{i}\text{,}%
\]
and the important fact is that $l^{\prime}$ is the dual line $l^{\ast}$.
Namely, we have $(q)=(p^{\ast})$, so the dual Pl\"{u}cker coordinates are the
Pl\"{u}cker coordinates of the dual line, and thus the duality map
(\ref{dual1}) extended to lines is the involutive map $l\rightarrow l^{\ast}$.

Pl\"{u}cker[1869] also introduced another type of mappings of lines to lines,
namely a general polar relationship depending on a given quadratic line
complex $\mathcal{C}.$ Recall from Section 3.5.2 that the complex lines in a
non-singular plane $\Pi\simeq P^{2}$ envelop a curve of degree 2, which yields
a polarity in the plane in the usual sense. Thus, for a given line $l$, if we
consider the pencil of planes $\pi$ containing $l$, then $l$ has a pole in
each of the (non-singular) planes $\pi$. Pl\"{u}cker showed these poles
actually lie on a common line $l^{\prime}$. However, the correspondence
$l\rightarrow l^{\prime}$ is not involutive (or reciprocal ) since
$l^{\prime\prime}$ is generally not equal to $l$.

The celebrated Pl\"{u}cker imbedding is the construction which regards $(p)$
in (\ref{(p)}) as a point in projective 5-space $P^{5}$, whereby the set of
lines in $P^{3}$ becomes a quadratic hypersurface
\begin{equation}
M_{2}^{4}\subset P^{5}:\mathcal{P}=0 \label{imbedding}%
\end{equation}
One also observes that the coordinates $p_{ij}$ in (\ref{(p)})\ are associated
with the coordinate polyhedron of $P^{5}$ whose six vertices are the edges
$p_{i}p_{j}$ of the fundamental tetrahedron in $P^{3}$, simply because the
only nonzero coordinate of the line $\overline{p_{i}p_{j}}$ is $p_{ij}$.
Moreover, given two lines $\overline{\mathbf{xy}}=(p)$ and $\overline
{\mathbf{zw}}=(q)$, expansion of the determinant of the 4$\times$4-matrix
$[\mathbf{x,y,z,w}]$ into 2-minors yields
\[
\det[\mathbf{x,y,z,w}]=\sum p_{ij}q_{ij}^{\ast}%
\]
so the vanishing of this expression is the condition that the lines intersect.

The 4-dimensional variety $M_{2}^{4}$ is usually referred to as the
Pl\"{u}cker quadric, but sometimes also as the Klein quadric. Indeed, it was
Klein who completed Pl\"{u}cker's volume [1869], and with his continuing works
in 1869 he added an extra geometric touch to line geometry, relating it to
projective geometry in dimension 5. For example, a line complex of order $n$
is defined by a homogeneous polynomial $X=X(p_{ij})$ of degree n, and the line
complex identifies with the 3-dimensional variety \{$\mathcal{P}=0,X=0$\} in
$M_{2}^{4}$, which by B\'{e}zout's theorem is of degree $2n$ in general.
Similarly, the intersection of three line complexes $X_{i}=0$ is generally an
algebraic curve in $P^{5}$
\[
\mathcal{\ P}=0,\ X_{1}=0,X_{2}=0,X_{3}=0
\]
of degree $2n_{1}n_{2}n_{3}$, again by Bezout's theorem, and the curve
represents a ruled surface in 3-space $P^{3}$.

As an example of a 2nd order special line complex, consider the totality of
lines which intersect the spherical circle at infinity (\ref{circular}). It is
also the family of lines satisfying Monge's equation (in Cartesian
coordinates)
\begin{equation}
dx^{2}+dy^{2}+dz^{2}=0 \label{differential}%
\end{equation}
In Pl\"{u}cker coordinates the complex is given by $p_{14}^{2}+p_{24}%
^{2}+p_{34}^{2}=0$. In the beginning of his career (1868-1870), Lie jokingly
referred to the lines as "verr\"{u}ckten Geraden" (crazy lines), as they
seemed to have many paradoxical properties, such as having zero length and
being perpendicular to themselves. Later he called them \emph{minimal lines},
but the term \emph{isotropic lines}, due to the French geometer Ribeaucour,
was commonly used.

\section{Non-Euclidean geometry}

One of the greatest mathematical discoveries in the 19th century is that of
non-Euclidean geometry, which did so profoundly affect our conception of space
and the entire foundation of geometry. The historical developments prior to
this achievement centered around the ancient problem of parallels and the
"truth" of the parallel postulate. From the days of Euclid to the middle of
the 19th century, many prominent scholars have taken risky steps towards a
possible solution of the problem, only to have their name added to the long
list of past failures. So, although there were many actors in this historical
drama, the actual discovery of the new geometry is attributed to Gauss, Bolyai
and Lobachevsky. Their revolutionary findings in the first third of the
century, based on the counter-intuitive assumption that the Euclidean parallel
postulate is actually false, were unique in the history of mathematics. \ 

\subsection{The discoveries of Bolyai, Lobachevsky, and Gauss}

The Hungarian J\'{a}nos Bolyai did most of his creative work in the 1820's, as
a military engineer in the Austrian army. He was the son of a mathematics
professor, Wolfgang (Farkas) Bolyai (1775-1856), who was familiar with the
problem of parallels and strongly warned his son against wasting his life on
this problem "as 100 geometers before had done". But already in 1823 the young
man wrote back to his father that he had "created a new and different world
out of nothing". Unfortunately, he was not able to publish his treatise until
1832, when it appeared as a 28-page Appendix in a two-volume geometry textbook
by his father. Wolfgang had been a friend of Gauss since their student days in
G\"{o}ttingen in the 1790's, and a copy of the book was sent to Gauss.
However, the approval from Gauss was not so well received by the Bolyais, and
it seemed to have a devastating effect on the young Bolyai. His mathematical
career almost ceased, mainly due to\ discouragements and mental depression
and, in fact, during his lifetime he received no public recognition for his work.

Let us also mention that J\'{a}nos Bolyai proved the interesting result that
in the new geometry it is possible to construct by ruler and compass a square
of area equal to that of a circle of radius 1. This is the "quadrature of the
circle", one of the most celebrated ancient problems, which already appeared
in our oldest mathematical document, the Papyrus Rhind (2000 BC). Its
impossibility in Euclidean geometry was finally settled by F. Lindemann in
1882, as a consequence of his proof of the transcendence of $\pi$.

The approach of Nicolai I. Lobachevsky, at the University of Kazan in Russia,
was amazingly similar to that of Bolyai, but probably he never knew about him.
For various reasons, Lobachevsky firmly believed the foundations of the
Euclidean geometry were flawed. The first ideas of his alternative
\emph{imaginary }geometry were apparently set forth in a public lecture in
1826, on the principles of geometry and the theory of parallels, but he failed
to obtain a publication out of it. On the other hand, his long paper \emph{On
the Elements of geometry} which appeared in two parts in Kazan Vestnik in
1829-30, is probably the first publication ever on non-Euclidean geometry. As
we shall see, however, \emph{Elements} from this geometry had appeared already
in printed books by Saccheri and Taurinus.

Lobachevsky wrote several papers on the topic and therefore went further, but
not necessarily deeper, than Bolyai. He made\ two attempts to convey his ideas
outside Russia, namely with two publications in Berlin in 1837 and 1840. The
first was an account in French, \emph{G\'{e}om\'{e}trie imaginaire, }published
in the new Crelle's Journal, heavy with $\ $formulae and dependent on his
papers from 1829-30, so it was virtually impossible to read. The second was
the more readable booklet \emph{Geometrische Untersuchungen}, written in
German, and he sent a copy to Gauss, without knowing about his interest in
this topic. As with Bolyai, however, Lobachevsky's work was not so appreciated
in his lifetime, and in 1846 he was even fired from the university. Indeed,
the only acclaim he was to receive during his lifetime was the appointment in
1842, recommended by Gauss, of his membership at the Academy of Science in G\"{o}ttingen.

Gauss, the foremost mathematician of his time and professor at G\"{o}ttingen
since 1807, had certainly anticipated some of the results of Bolyai and
Lobachevsky on non-Euclidean geometry. In letter correspondences he praised
their talents and the geometric spirit of their work. However, as to the
extent of his own investigations we can only judge from various remarks in
private letters, unpublished notes, and his reviews of books relating to the
theory of parallels. On the basis of this, it is amazing that despite his
great reputation he was afraid of making public his own geometric discoveries
on the subject. But these were ideas which undoubtedly would have refuted
Kant's position on the nature of space and the unique role of the Euclidean space.

\subsection{Absolute geometry and the Euclidean parallel postulate}

For more than 2000 years, Euclidean geometry was the true and real geometry,
namely the Geometry which was regarded as the science of the Space we live in.
Euclid's \emph{Elements} had introduced five basic postulates expressing
self-evident \ properties of points, lines, right angles etc., and proceeded
to deduce altogether 465 propositions by mathematical reasoning. In Section
1.1 we stated these postulates in the spirit of Euclid and denoted them
respectively by $E1,E2,E3,E4,E5$, but with modern critical eyes the first four
postulates are certainly vague statements. Moreover, they are incomplete as
far as rigorous proofs of the propositions are concerned, since many
additional "evident" assumptions must have been tacitly used as well. Keeping
this in mind, we recall that many scholars, even long before Gauss, Bolyai and
Lobachevsky, were led to investigate the restricted geometric content based
upon $E1-E4$, namely without assuming $E5$. In the 1820's Bolyai referred to
this geometry as \emph{absolute geometry} and, for convenience, in the sequel
we shall also do so.

The logical status and geometric implications of the parallel postulate $E5$
had remained a challenge since the days of Euclid. Recall that two lines in
the plane are said to be \emph{parallel} if they do not intersect. It is worth
noticing that Euclid deduced the first 28 propositions without using $E5$, in
other words, they are statements of absolute geometry. In this geometry it
follows, for example, that for every line $l$ and point $p$ outside $l$ there
is at least one parallel line passing through $p$.

Now, the last postulate $E5$, which distinguishes absolute geometry from
Euclidean geometry, is in fact equivalent to the statement that the above
parallel line through $p$ is unique. This alternative version of $E5$ dates
back to Proclus (411-485), but it is known today as Playfair's axiom, after
John Playfair (1748-1819) in Edinburgh, who published in 1795 a new edition of
Book I-VI of the \emph{Elements}.

By assuming the postulate $E5$ one can, for example, prove the following two
basic results about triangles, namely the theorem of Pythagoras and the
statement that the angle sum is $2R$\emph{,} where $R$ denotes a right angle.
Conversely, by accepting these results as obvious or experimental facts, one
is also accepting the "truth" of $E5$. This may explain why the search for a
"proof" of $E5$ became such an important issue, and there seemed to be two
directions to proceed, either

\begin{enumerate}
\item to prove $E5$ as a logical consequence of the postulates $E1-E4$, or

\item to establish the truth of $E5$ from the laws of nature.
\end{enumerate}

Since antiquity many outstanding scholars have followed the first approach,
hoping to show that $E5$ is a superfluous postulate. If the scholar succeeded,
he would have proved that absolute geometry is the same as Euclidean geometry.
But, as we shall see, the scholar would typically invoke another
"self-evident" assumption, and in the misbelief that it was true in absolute
geometry he would use it to deduce $E5$. On the other hand, despite the long
lasting belief that the geometric truths are encoded into the nature, the
second direction was not seriously considered until the 19th century, perhaps
as a last effort.

\subsection{Attempted proofs of the parallel postulate}

The first known attempted proofs date back to the ancient scholars Ptolemy
(ca. 85-165 AD) and Proclus. Of course, their arguments were flawed since the
decisive assumption was merely a disguised version of $E5$ itself. Later
examples of this kind are, besides Playfair's axiom, the axioms named after C.
Clavius (1537-1612), R. Simson (1687-1768), J. Wallis (1616-1703), and
Clairaut (cf. e.g. Greenberg [2001], Gray [2007]). For example, the axioms of
Clavius and Clairaut assert the existence of two equidistant lines or a
rectangle, respectively. Moreover, the existence of a triangle with angle sum
$2R$ also yields $E5$.$\ $Some geometers such as Clairaut or Wallis, perhaps
after realizing how hopeless it was to prove $E5$, proposed to replace $E5$ by
their own axiom in order to improve or simplify Euclid's geometry.

The Persian scholars Omar Khayy\'{a}m (1048-1131) and Nasir Eddin al-Tusi
(1201-1274) are also known for their noteworthy analysis of the parallel
postulate. Their ideas made their way to Europe and may, in fact, have
contributed to the development of non-Euclidean geometry many centuries later.
For example, the quadrilateral named after the Jesuit priest and logician G.
Saccheri (1667-1733), was already introduced by them many centuries earlier.
Saccheri was a student of Giovanni Ceva's brother Tommasco (1648-1737), who
was a professor of mathematics and rhetoric at a Jesuit college in Milan.
Saccheri tried to prove $E5$ by reductio ad absurdum, assuming the negation of
$E5$, and so he attempted to deduce a contradiction from the ensuing bulk of
non-Euclidean results. His little book \emph{Euclid Freed of Every Flaw},
which appeared a few months before he died, created something of a sensation
and was examined by leading mathematicians of the day. But it fell into
oblivion after some years, until it was rediscovered by Beltrami in 1889.

The work of Saccheri was most likely known to the Swiss-German scholar and
leading mathematician of the 18th century, J.H. Lambert (1728-1777), who
proceeded similarly and explored more deeply the consequences of the negation
of the parallel postulate. He was a man of extraordinary insight and published
more than 150 works in various areas. In 1764 he was invited to become a
colleague of Euler and Lagrange at the Preussian Academy of Science\ in
Berlin. Euler had established in 1737 that $e$ and $e^{2}$ are irrational
numbers, and in an outstanding paper of 1768 Lambert showed that $\pi$ is also
irrational. He further conjectured that $e$ and $\pi$ are even transcendental,
but it was in the next century that Hermite and Lindemann, respectively,
verified this conjecture.

Lambert's work \emph{Theory of Parallels} was written in 1766, but perhaps due
to his unsatisfaction with the work it was only published posthumously, by
Johann Bernoulli III 20 years later. Saccheri and Lambert must have regarded
the new geometry as ficticious and without reality. Saccheri even believed,
but wrongly of course, that he had established the truth of $E5$, whereas
Lambert admitted that his attempts had failed. Like ancient geometers and
Clavius as well, Lambert had wrongly assumed that the curve of points
equidistant from a given line and on the same side of the line, is itself a line.

In France, the leading analyst A.M. Legendre , influential in the
restructuring of higher education during and after the revolution, was also
confronted with the parallel postulate and its role in geometry. He became
obsessed with proving its truth, and during the years 1794 to 1833 his 12
different attempts appeared, one after another, in the appendix of the revised
editions of his highly successful \ textbook \emph{\'{E}l\'{e}ments de
G\'{e}om\'{e}trie.}

A famous mistake of Legendre was his assumption that through any interior
point of an angle one can always draw a line which cuts both sides of the
angle, but Legendre never realized that this is just another disguised version
of the postulate $E5$. On the other hand, 100 years after Saccheri, Legendre
rediscovered Saccheri's results in absolute geometry, for example, that the
angle sum is at most $2R$ for any triangle.

Truly, in the 18th century, elementary geometry was rather engulfed in the
problems raised by the parallel postulate. The situation was well illustrated
by the thesis of G.S. Kl\"{u}gel in G\"{o}ttingen in 1763, who described the
flaws of 28 different attempted proofs of $E5$. In 1767 the leading French
scholar J. L.R. d'Alembert (1717-83) referred to the accumulation of false
proofs and lack of progress as the \emph{Scandal of geometry}. Although
Kl\"{u}gel expressed doubt that $E5$ could ever be proved, he did not scare
off but rather inspired scholars like Lambert to try their fortune. Even Gauss
had been working on the parallel postulate since 1792, at the age of 15, but
having made little progress by 1813 he wrote:
\begin{align*}
&  \text{\emph{In the theory of parallels we are even now not further than
Euclid.}}\\
&  \text{\emph{This is a} \emph{shameful part of mathematics}}\emph{..}%
\end{align*}

Gradually but slowly, scholars became convinced that Euclid's parallel
postulate cannot be proved, namely it is independent of the other Euclidean
postulates. For example, Gauss expressed his conviction in 1817 in a letter to
the astronomer H.W.M. Olbers (1758-1840). However, the truth of the parallel
postulate continued to be an unsettled question for quite another reason, due
to the poorly understood relation between abstract "mathematical geometry "
and "physical geometry". So,\ let us turn to the second direction of attempted
proofs of $E5$, as pointed out above.

The first ancient cosmological model of the universe bounded by the celestial
sphere dates back to Eudoxus, who was also the founder of theoretical
astronomy, and this model was refined in Ptolemy's 13 volume treatise
\emph{Almagest.} Even Kepler and Galileo regarded the universe to be limited.
But a physical interpretation of Euclid's postulates would certainly be false
if the geometry was spherical. In fact, the viewpoint that we live in
infinite, unbounded Euclidean space dates back no further than to Descartes in
the 17th century. This model became the geometric frame for Newtonian
mechanics and was later adopted by the Kantian philosophy.

But now, with the emerging non-Kantian ideas about geometry and space, the
question came up whether we rather live in an infinite and unbounded
non-Euclidean space. The real Geometry was supposed to represent physical
space, and therefore, as in the case of other laws of nature, the truth could
be established from physical experiments such as astronomical observations. As
a geometer, Lambert knew that similar triangles would, in fact, be congruent
if the postulate $E5$ was wrong, and moreover, the angle sum of a triangle
would be less than $2R$. In that case, as an astronomer he would worry about
the countless inconveniences, and astronomy would be "an evil task".

Lobachevsky measured in 1829 the parallax of stars, which is almost
negligible, so his observations were inaccurate and hence inconclusive. Gauss
shared Lobachevsky's view in his letter correspondences with his many
astronomer friends, such as Olbers, W. Bessel (1784-1836) and C.L.
Gerling(1788-1864). He discussed with them the possibility that physical space
was not necessarily Euclidean. Indeed, while surveying the estates of Hanover
he set up theodolites on three mountain peaks to test the non-Euclidean
hypothesis experimentally, cf. Coxeter [1998].

Gauss wrote to Olbers that we should not put geometry on a par with arithmetic
that exists purely a priori, but rather with mechanics. The viewpoint of
geometry as an empirical science is also exemplified by the leading French
mathematicians Lagrange and Fourier, who tried to deduce the parallel
postulate from the law of the lever in statics. Lobachevsky was even more
extreme, and being sceptical to the very foundations of Euclidean geometry, he
believed that knowledge about the motion of bodies would help building up the
concepts of geometry, such as ideas about the straight line.

\subsection{The emergence of non-Euclidean geometry}

Let us first have a closer look at Saccheri's approach. On the basis of
absolute geometry, he focused attention on a special quadrilateral $ABCD$,
where the opposite sides $AD$ and $BC$ are congruent and perpendicular to the
base $AB$. Denoting by $\alpha$ the angles at $C$ and $D$ (which are
congruent), Saccheri considered the three possible cases
\begin{equation}
(i)\text{ }\alpha=R\text{, \ \ \ }(ii)\text{ \ }\alpha>R\text{, \ \ }%
(iii)\text{ \ }\alpha<R\text{, } \label{Saccheri}%
\end{equation}
referred to as the hypothesis of the right, obtuse, and acute angle,
respectively. Saccheri rightly argued that case (ii) is impossible, being
incompatible with postulate $E2$ on the indefinite extension of lines. Case
(i) simply means the quadrilateral is a rectangle, namely its angle sum is
$4R$, and consequently the geometry must be Euclidean. Thus the negation of
the parallel postulate $E5$ amounts to the acute angle hypothesis, $\alpha<R
$, from\ which Saccheri deduced many "strange" geometric results. Regrettably,
as Lambert must have observed, Saccheri concluded with an obscure argument
that he had produced a contradiction, so contrary to his own belief he failed
to establish the truth of $E5$. However, with all the "strange" geometric
results Saccheri had actually discovered non-Euclidean geometry, although he
did not recognize it as such. For example, he deduced the existence of lines
approaching each other infinitely close without having intersection -- later
these\ became known as \emph{asymptotic} \emph{parallels}.

\subsubsection{The trigonometry of Lambert, Schweikart, and Taurinus}

Lambert was baffled by the observation that in non-Euclidean geometry there
would be an absolute measure of length, analogous to the measure of angles.
Therefore, similar triangles would also be congruent. For example, all
equilateral triangles with angle $\alpha=50^{\circ}$ say, are congruent, so
their side $s_{0}$ would yield a natural unit of length. Lambert further
noticed that the area of a triangle with angles $\alpha,\beta,\gamma$ (in
radian measure) is proportional to $(\pi-\alpha-\beta-\gamma)$. Knowing that
the area of a triangle on a sphere of radius $r$ is
\begin{equation}
A_{1}=r^{2}(\alpha+\beta+\gamma-\pi), \label{area1}%
\end{equation}
Lambert expressed the area in the former case as
\begin{equation}
A_{2}=k^{2}(\pi-\alpha-\beta-\gamma)=r^{2}(\alpha+\beta+\gamma-\pi)\text{,
\ \ \ where }r=k\sqrt{-1}\ \label{area}%
\end{equation}
and $k$ is a positive constant depending on $s_{0}$. Then he proclaimed that
the geometry behaves like an \emph{imaginary }sphere of radius $r$.

Lambert's ideas were continued by two amateur geometers, the law professor
F.K. Schweikart (1780-1859) and his nephew F.A.Taurinus (1794-1874), who
progressed further using analysis rather than following the classical
approach. From the outset Schweikart accepted with no prejudice the new
geometry where the angle sum of triangles is less than two right angles. He
was not looking for a contradiction, but rather speculated if his geometry,
which he referred to as \emph{Astral Geometry}, would be appropriate for the
study of the physical space.

In a noteworthy memorandum, communicated in 1818 to Gauss via the astronomer
Gerling, Schweikart introduced his own constant $\mathfrak{C}$, in analogy
with Lambert's constant $k$ in (\ref{area}), to be the maximal height of any
right-angled isosceles triangle. He could not determine its value, and perhaps
regarded it as a parameter for many possible astral geometries, but he pointed
out that the geometry would be Euclidean if $\mathfrak{C}$ is infinite. Gauss
complimented him on his results, and remarked that the area of a triangle
would have the upper bound $\pi k^{2}$ when $\mathfrak{C}$ is expressed as
\begin{equation}
\mathfrak{C}=k\log(1+\sqrt{2}). \label{C}%
\end{equation}
In letters to Taurinus and Bessel in the 1820's, Gauss expressed the view that
the non-Euclidean geometry is self-consistent and entirely satisfactory, but
the parameter cannot be determined a priori. As an astronomer, Bessel bravely
suggested to Gauss in 1829 that physical space is maybe slightly non-Euclidean.

On the other hand, Taurinus firmly believed in the truth of Euclidean
geometry, and his work was motivated by his desire to prove the parallel
postulate. Therefore, he continued with his uncle's work to prepare himself
and to better understand geometry in general. Up to 1825, when his first
booklet appeared, he still believed Euclidean geometry was the unique
geometry, but in his second booklet in 1826 he accepted the internal
consistency and lack of contradictions in his uncle's astral geometry, and
even vaguely suggested it might be the geometry of some surface. With due
regard to Saccheri's book, the above booklets, published in Cologne
(C\"{o}ln), seem to be the first printed expositions on the \emph{Elements} of
non-Euclidean geometry.

Most noticeable, Taurinus broke with the traditional synthetic approach and
introduced trigonometry as a method in non-Euclidean geometry, with the usual
trigonometric functions replaced by Lambert's \emph{hyperbolic} functions
\[
\sinh x=\frac{1}{2}(e^{x}-e^{-x})=i\sin\frac{x}{i}\text{, \ }\cosh x=\frac
{1}{2}(e^{x}+e^{-x})=\cos\frac{x}{i}%
\]
Indeed, Lambert himself had missed the connection which the functions provide
between analysis and non-Euclidean geometry. Taurinus started with the
trigonometric relations for a triangle $ABC$ on a sphere of radius $k$, for
example the cosine law
\begin{equation}
\text{ \ }\cos\frac{a}{k}=\cos\frac{b}{k}\cos\frac{c}{k}+\sin\frac{b}{k}%
\sin\frac{c}{k}\cos A \label{cos}%
\end{equation}
where $a$ is the side opposite to the vertex $A$ etc., and by recalling
Lambert's obscure idea of an imaginary sphere, he formally substituted
$k\rightarrow$ $ik$ into the equations. Thus, for example, the law (\ref{cos})
becomes the hyperbolic cosine law
\begin{equation}
\cosh\frac{a}{k}=\cosh\frac{b}{k}\cosh\frac{c}{k}-\sinh\frac{b}{k}\sinh
\frac{c}{k}\cos A\text{\ } \label{cosh}%
\end{equation}
Similarly he obtained the identity%
\begin{equation}
\cos A=\sin B\sin C\cosh\frac{a}{k}-\cos B\cos C \label{cosA}%
\end{equation}

By a simple calculation the limiting case as $k\rightarrow\infty$ is readily
seen to be the Euclidean geometry; for example, formula (\ref{cosh}) yields
the (usual) cosine law
\[
a^{2}=b^{2}+c^{2}-2bc\cos A
\]
Taurinus called the geometry satisfying the resulting hyperbolic laws for
\emph{log-spherical geometry}. He found that it agreed with Schweikart's
Astral Geometry; he also expressed his uncle's constant $\mathfrak{C}$ in
terms of $k$, and it agreed with the above formula (\ref{C}) of Gauss.

Still, Taurinus was not convinced by the nice analytic results, which appeared
in the second booklet. After all he had only worked out the geometry of an
imaginary sphere, perhaps with no real content. He had corresponded with Gauss
about the geometric ideas, but now he also insisted that Gauss should publish
his own ideas on the subject. As a result, however, the angered Gauss ended
his correspondence with Taurinus, and after this, the story goes, Taurinus
ceased his geometric studies, bought up his own published booklets and burned them.

\subsubsection{On the work of Bolyai and Lobachevsky}

In the 1820's there were also mathematicians with the conviction that there is
a plane geometry in which Euclid's parallel postulate is wrong, and Bolyai and
Lobachevsky were the first ones with enough confidence in their geometric
ideas to publish them. Probably they never knew about Taurinus, but his work
would have provided the analytical basis of their geometric study. Instead,
they derived in their own way trigonometric $\ $formulae including those found
by Taurinus. Bolyai also became quite an expert on\emph{\ absolute} geometry,
that is, the aggregate of those propositions of Euclidean geometry which are
independent of the parallel postulate and hence common to both Euclidean and
non-Euclidean geometry. His unifying sine law for triangles $ABC$%
\begin{equation}
\frac{\sin A}{\bigcirc_{a}}=\frac{\sin B}{\bigcirc_{b}}=\frac{\sin C}%
{\bigcirc_{c}}, \label{sine law}%
\end{equation}
where $\bigcirc_{a}$ denotes the circumference of a circle of radius $a$, is
valid for absolute geometry and spherical geometry as well, see also
(\ref{circle1}).

Bolyai and Lobachevsky were also the first who worked in more generality on
the non-Euclidean geometry in three dimensions, perhaps having in mind that it
could be the geometry of physical space. Therefore, they also initiated the
task of reformulating the basic laws of classical mechanics in the new
geometric setting. For example, in 1835 Lobachevsky considered the
gravitational law and defined the Kepler problem with an attractive force
inversely proportional to the area of the 2-dimensional sphere of radius equal
to the distance between the bodies. Bolyai came up with similar ideas about
the same time (see Diacu et al [2008]). But now, let us briefly recall their
basic geometric approach.

A \emph{horocycle} is the plane curve perpendicular to a family of asymptotic
parallel lines; it is also the limiting curve of an expanding circle whose
radius tends to infinity. In the 3-space there are surfaces defined similarly,
called \emph{horospheres} by Lobachevsky. Bolyai and Lobachevsky made use of
this type of curves and surfaces\emph{\ }in their analysis; in fact, they both
attacked plane geometry via the horosphere. They made the remarkable
observation that the induced geometry of a horosphere is that of a Euclidean
plane, whose straight lines correspond to the horocycles lying on the surface.
In fact, this result dates back to F.L. Wachter (1792-1817), a pupil of Gauss
in 1816, who inquired about the limiting form of a sphere in non-Euclidean
3-space when its radius becomes infinite.

Let us consider a triangle $PQR$ with vertices $P,Q,R$ and opposite sides of
length\ $p,q,r$, respectively. We also assume $QR$ lies on the line $l$ and
$PQ$ is perpendicular to $l$, so $r$ is the distance from $P$ to $l$. Now,
keep the vertices $P$ and $Q$ fixed and let $R$ move far away along $l$. Then
the angle at $P$ increases towards a limiting value $\Pi(r)$, depending only
on $r$, and the angle at $R$ tends to zero. Moreover, the side $PR$ approaches
the asymptotic parallel ray of $l$ through $P$, so $\Pi(r)$ is also referred
to as the angle of \emph{parallelism} at $P$. Bolyai and Lobachevsky
discovered the following nice formula for $\Pi(r)$,\ measured in radians,
\begin{equation}
\Pi(r)=2\arctan e^{-r/k} \label{parallelism}%
\end{equation}
where $k$ is the constant appearing in the area formula (\ref{area}), namely
$\pi k^{2}$ is the maximal area of any triangle. In particular, we recall that
Schweikart's constant $\mathfrak{C}$ is the value of $r$ when $\Pi(r)$ is half
of a right angle, so the identity (\ref{C}) found by Gauss and Taurinus is
just a special case of the identity (\ref{parallelism}). In fact, Taurinus
could also have derived the formula (\ref{parallelism}) from the identity
(\ref{cosA}) applied to the above triangle $PQR$.

Next, let us recall Lobachevsky's crucial experiment in 1856 to test whether
his geometry would fit better than the Euclidean one. This test differs from
his first effort in 1829 and is described in \emph{Pang\'{e}om\'{e}trie},
which he dictated as a blind man the year he died. With reference to the above
triangle $PQR$, imagine a distant star at $R$ which is observed at the two
diametrically opposite positions $P$ and $Q$ of the Earth's orbit. Then the
observed parallax should exceed the angle $\pi/2-\Pi(r)$. However, due to the
experimental errors the finiteness of the unknown absolute unit $k$ could not
be established, only that $k$ would be several million times larger than the
diameter$\ r$ of the orbit.

Lobachevsky strongly believed that his geometry was consistent, that is, it
would not lead to any contradictions. For him, geometry was concerned with
measurements and numbers, related by $\ $formulae whose validity is rather an
algebraic problem. He argued that the geometry is based on $\ $formulae for a
triangle, and they would yield the familiar formulae for a spherical triangle
when the sides $a,b,c$ are replaced by $ia,ib,ic$, or when the parameter $k$
is replaced by $ik$. Moreover, the formulae would still make sense and
describe the Euclidean geometry in the limit when $k$ tends to infinity.

The above arguments of Lobachevsky were, in fact, insufficient (see e.g.
Rosenfeld [1988:227]). On the other hand, Gauss and Bolyai were also familiar
with the above trigonometric relations, but they did not regard them as
evidence for the logical consistency of the non-Euclidean geometry. In his
letter to Taurinus in 1824, Gauss had expressed his belief that the new
geometry is self-consistent, but on the other hand, Bolyai was quite disturbed
by his inability to settle this matter. In fact, it remained an open question
until 1868, when Beltrami came across the crucial idea of constructing a
concrete \emph{model} for the geometry, which simply reduced the whole
question to the consistency of Euclidean geometry itself.

Certainly, the recognition of the work of Bolyai and Lobachevsky was hampered
by the fact that none of them were successful in establishing the truth, or
the physical existence, of their geometry. After all, almost everyone believed
that the physical space is Euclidean, so how could an alternative geometry
also be true? However, after Gauss's death in 1855, his unpublished notebooks
and letter correspondences became known among\ mathematicians. In particular,
the publication of his letters to H.C. Schumacher (1780-1850) in 1864, in
which he praised Lobachevsky's work and expressed his own belief that Space
might be non-Euclidean, made a strong impression on European mathematicians
(Rosenfeld [1988]). As a consequence, during the next few years the original
papers of Lobachevsky and Bolyai's Appendix from 1832 attracted considerable
attention and were translated to German, French, Italian, and Russian. Let us
briefly recall what actually happened during these years.

In 1866, when R. Baltzer in Germany was preparing the second edition (1867) of
his textbook \emph{Die Elemente der Mathematik}, he also included a favorable
mention of the discoveries of Bolyai and Lobachevsky. He informed Ho\"{u}el in
France, who issued the same year a French translation of Lobachevsky's German
memoir from 1840, together with excerpts from the Gauss-Schumacher
correspondence. Moreover, in 1867 Ho\"{u}el also wrote a book in order to
explain Lobachevsky's geometry. But on the other hand, Bolyai's Appendix was
not so easily available, so the translation of that paper was delayed until
1868. An historical note on the life and work of Lobachevsky, by E.P.
Yanisevskii in Kazan, appeared also translated to French and Italian in 1868.

\subsubsection{Beltrami and his model of non-Euclidean geometry}

In Italy, Battaglini, Beltrami\ and Cremona were active proponents for the new
geometric ideas. The journal Giornale di matematiche, which Battaglini founded
in 1863, became the major Italian publication channel for papers on
non-Euclidean geometry. Battaglini published in 1867 a paper on Lobachevsky's
imaginary geometry and a translation of Lobachevsky's last paper
\emph{Pang\'{e}om\'{e}trie}, as well as a translation of Bolyai's Appendix
in1868. For geometry in Italy, Cremona was certainly very influential, but his
own works were mostly within projective geometry.

Beltrami had learned about non-Euclidean geometry by reading the French
translations of Ho\"{u}el in 1866, and in 1868, shortly after the publication
of Riemann's famous habilitation lecture in 1854, Beltrami quickly developed a
deeper understanding of the topic by following Riemann's approach to geometry.
Thus appeared his two seminal papers [1868a,b], which settled for the first
time the ancient question of a proof of the parallel postulate, namely he
demonstrated that no proof is possible. He achieved this by exhibiting a
Euclidean model of the geometry of Bolya and Lobachevsky.

By taking a closer look at Beltrami's paper [1868a], one will find that he
gives (at least) two different Euclidean models for the hyperbolic plane, with
explicit formulae for the metric in the sense of Riemann. For simplicity,
consider the upper hemisphere $z>0$ of radius 1 in Euclidean 3-space with
rectangular coordinates $(x,y,z)$. Beltrami modifies its spherical metric by
multiplying it with $z^{-2}$, and therefore, in terms of spherical polar
coordinates $(\phi,\theta$), centered at the north pole so that $z=\cos\phi$,
the metric becomes
\begin{equation}
ds^{2}=\frac{d\phi^{2}+\sin^{2}\phi\text{ }d\theta^{2}}{\cos^{2}\phi}
\label{beltrami1}%
\end{equation}

This is, in fact, a model of the hyperbolic plane with curvature $K=-1$, whose
geodesics (lines) are the semicircles on the sphere lying in planes
perpendicular to the $xy$-plane. Moreover, by vertical projection of the
hemisphere onto the unit disk $x^{2}+y^{2}<1$, Beltrami obtains the so-called
disk model of hyperbolic geometry, where the hyperbolic lines are, indeed, the
Euclidean line segments. Finally, Beltrami further transforms the metric to
hyperbolic polar coordinates $(\rho,\theta)$ centeret at the midpoint of the
disk, namely $\rho$ measures the radial distance from the center of the disk.
The connection between $\phi$ and $\rho$ is found to be $\sinh\rho=\tan\phi$,
which yields the very simple expression
\begin{subequations}
\begin{equation}
ds^{2}=d\rho^{2}+\sinh^{2}\rho\text{ }d\theta^{2} \label{beltrami2}%
\end{equation}

In England, the algebraist and geometer Arthur Cayley was a leading figure,
and with his numerous papers he treated nearly every subject of pure
mathematics. In the paper entitled \emph{Notes on Lobachevsky's Imaginary
geometry }(1865) he made a comparison of the spherical trigonometric formulae
with those of the Lobachevskian geometry, but apparently he overlooked the
essense of Lobachevsky's discovery, although his writings helped promoting the
new geometric ideas (cf. Rosenfeld [1988: 220]). On the other hand, Cayley's
important theory of projective metrics in his \emph{Sixth Memoir} [1859],
where he proposed a generalized definition of distance, was grasped by the
young Felix Klein in 1871, who caught an important idea which was, in fact,
overlooked by both Cayley and Beltrami, see Klein [1871g], [1872c], Stillwell [1996].

In summary, due to the above mentioned publicity, around 1870 the ideas of
Bolyai and Lobachevsky were known to geometers at the major universities in
Europe. Beltrami (1868) and Klein (1871), the latter from a projective
geometric viewpoint, had finally completed the last step needed for
non-Euclidean geometry to be accepted as part of ordinary mathematics, by
their construction of convincing Euclidean and projective models which showed
that the new geometry was equally consistent with the ancient geometry. The
conformal models of non-Euclidean geometry exhibited in Poincar\'{e} [1882]
are, in fact, implicit among the models presented in Beltrami [1868b]. But by
applying them to his study of automorphic functions Poincar\'{e} also
contributed largely to the uprise of hyperbolic geometry to a respectable
mathematical discipline, which by 1890 was finally taught as a course at major universities.

\section{The classical geometries: Euclidean, spherical, and hyperbolic}

\subsection{A unified view of the three classical geometries}

In the 1860's the geometry of Bolyai and Lobachevsky, which Gauss had referred
to as non-Euclidean geometry, became also known as Lobachevskian geometry. In
fact, still more geometries with non-Euclidean properties were expounded by
Riemann, Klein, and others, and in fact, \emph{spherical} geometry was also
called Riemann's non-Euclidean geometry. Namely, the latter is the geometry of
figures drawn on a round spherical surface in 3-dimensional Euclidean space.

Being so closely related to Euclidean geometry, spherical geometry is really
the first example in history of a geometry other than the Euclidean one. The
ancients called it \emph{Sphaerica}, and they used it to describe the heavenly
bodies moving around on the celestial sphere. For example, Menelaus of
Alexandria and the Arabs (around 1000 AD) studied this geometry. The French
Albert Girard (1595-1632), whose treatise on trigonometry in 1626 contained
the first use of the abbreviations \emph{sin, cos, tan}, also gave the formula
(\ref{area1}) for the area of a spherical triangle, a formula which was later
generalized by Gauss to geodesic triangles on more general surfaces in
3-space. We refer to Rosenfeld [1988] for detailed information on the early
history of spherical geometry.

Euclidean, spherical, and Lobachevskian geometry constitute the three
\emph{classical geometries}, and for many good reasons. They have many
properties in common, deeply rooted in the human conception of space, and the
most basic ones are of a non-metric nature. \ But they are also, somehow,
elaborated within another type of "classical" geometry, namely the grand
unifying theory called projective geometry. However, in this section we shall
rather focus on the metric properties of the classical geometries, some of
which they have in common, but certainly there are also important differences.

One of Klein's achievements in [1871g] was his unification of these
geometries, subsumed by a common ambient projective space, by utilizing the
projective measure construction of Cayley [1859]. Klein proposed the
suggestive terms \emph{parabolic, elliptic}, and \emph{hyperbolic }geometry,
respectively, in accordance with their geometric properties and limiting
behavior in resemblance with conic sections. We shall henceforth follow Klein
and use the modern term hyperbolic geometry instead of Lobachevskian geometry.
However, the equivalence of the terms "parabolic" and "Euclidean" space ceased
in the early 20th century, when more general parabolic spaces were defined as
specific homogeneous spaces constructed in terms of Lie group theory.

With the unifying analytic approach inspired by Riemann, around 1870 the focus
on classical geometries had shifted to their space forms, which we shall
denote by $E^{n}$, $S^{n}$, and $H^{n}$. Namely, their differential geometric
properties, and to some extent also topological properties, were the subject
of study in the lowest dimensions $n=2,3$. The higher dimensional versions
with $n>3$ were gradually accepted since Riemann had extended the classical
notion of \emph{curvature} for surfaces to spaces of higher dimensions, and
the above three types of geometries were found to have constant curvature $K$,
namely
\end{subequations}
\begin{equation}
E^{n}:K=0\text{, \ \ \ }S^{n}:K>0\text{, \ \ }H^{n}:K<0 \label{spaces}%
\end{equation}
\qquad

It was Beltrami who discovered\ the link between hyperbolic geometry and
spaces of constant negative curvature. He came across the idea when he
compared Lobachevsky's 1837 paper with Minding's 1840 paper on the geometry of
the \emph{pseudosphere}, both papers printed in Crelle's Journal. The latter
surface has constant negative curvature, and Beltrami observed in 1868 that
the two papers had, in fact, the same trigonometric formulae. With the formula
(\ref{metric4}) below, Beltrami had expressed in polar coordinates
$(\rho,\theta)$ the Riemannian metric of $H^{2}$ for $K=-1$, so he certainly
knew the analogous expressions in all cases (\ref{spaces}) for $n=2$ and any
$K$, namely%

\begin{equation}
ds^{2}=d\rho^{2}+f(\rho)^{2}d\theta^{2}\text{, \ \ \ \ \ \ }f(\rho)=\left\{
\begin{array}
[c]{cc}%
\rho & K=0\\
\frac{1}{\sqrt{K}}\sin(\sqrt{K}\rho) & K>0\\
\frac{1}{\sqrt{-K}}\sinh(\sqrt{-K}\rho) & K<0
\end{array}
\right.  \label{metric4}%
\end{equation}
The distance $\rho$ from the chosen center $O$ of the coordinate system ranges
over $[0,\infty)$ when $K\leq0$, whereas $0\leq\rho\leq\pi/\sqrt{K}$ in the
sphere case. Another remarkable property of the three geometries is the
unifying Bolyai's sine law (\ref{sine law}), where by (\ref{metric4})%

\begin{equation}
\bigcirc_{a}=2\pi f(a) \label{circle1}%
\end{equation}
is the length of the circle of radius $a$, say the circle with the equation
$\rho=a$, which also yields a nice geometric interpretation of the size
function $f(\rho)$.

We also remark that the trichotomy (\ref{Saccheri}) of the angle $\alpha$ in
Saccheri's quadrilateral reflects the trichotomy of the classical geometries.
Thus case (i) is the Euclidean geometry, based upon the ancient Euclidean
postulates $E1,E2,..,E5$, and hyperbolic geometry was discovered in case (iii)
by renouncing $E5$. Finally, Riemann gave spherical geometry his full
recognition as a kind of non-Euclidean geometry, satisfying Saccheri's
hypothesis (ii) of the obtuse angle, with the great circles on the sphere
$S^{2}$ interpreted as the straight lines. In this case $E5$ is violated since
there are no parallel lines at all, but $E2$ is also violated because the
circles are of finite lengtht. In addition, $E1$ postulates that two points
determine a unique line, so $E1$ is violated as well. For an axiomatic
approach, Riemann therefore proposed a modification of these three postulates
so that, for example, (i) two points would determine at least one line, and
(ii) a line is unbounded. However, such an approach to spherical geometry has
never been found useful, since the geometry is best understood as a
subgeometry of Euclidean geometry.

However, the incompatibility of $E1$ with spherical geometry remained an
unsatisfactory issue for many years. Beltrami [1868a] refers to $E1$ as the
\emph{postulate of the straight line}, and several mathematicians, including
Beltrami and Weierstrass, believed that this failure was a characteristic
property of geometries of constant positive curvature. It also seemed that
Riemann [1867] had identified, although vaguely, such geometries with
spherical geometry. But Klein [1871g: 604], [1874c] removed this misconception
by describing a truly \emph{elliptic geomety }where $E1$ holds. Let us briefly
recall the basic ideas.

Klein's elliptic model in dimension n was taken to be the space obtained from
the n-sphere $S^{n}$ by identifying antipodal points, whereby each great
circle is reduced to a closed geodesic of half the original length. In fact,
the construction identifies the elliptic n-space with the projective n-space
$P^{n}$. This is clear from by the vector space model of $P^{n}$ (see Section
3.4.1), since a central line in $\mathbb{R}^{n+1}$ cuts the surrounding sphere
$S^{n}$ in two antipodal points. For example, the elliptic plane identifies
with the real projective plane $P^{2}$; in fact, by leaving out metrical
concepts and congruence, elliptic geometry becomes real projective geometry.
From the Riemannian viewpoint, $P^{n}$ and $S^{n}$ cannot be distinguished
locally, that is, they are locally isometric. In particular, they have the
same constant curvature. Spherical geometry on $S^{n}$ became also known as
\emph{doubly elliptic} geometry since the mapping $S^{n}\rightarrow P^{n}$ is
two to one.

Before basic topological concepts and constructions had been developed, the
distinction between local and global properties of a geometry was poorly
understood. This explains why people generally believed, in the 1870's and
maybe even later, that the three classical spaces (\ref{spaces}) and Klein's
elliptic space $P^{n}$ were the only (proper) space forms of constant
curvature. But, in fact, in 1873 the English mathematician Clifford described
the construction of a flat torus, which can be embedded in the sphere $S^{3}$
as a closed surface of zero curvature, cf. note to the letter of 4.11.73.

Hawkins [2000] points out that Riemann's discussion of manifolds of constant
curvature impressed the late 19th century geometers. Already in his paper
[1872c] Klein comments upon the work of Riemann and Helmholtz, and as late as
in his lectures on hyperbolic geometry in 1889 Klein maintains that the
concept of an n-dimensional manifold of constant curvature is the most
essential result of Riemann's approach. With modern eyes, however, its real
importance is the generality of his approach, allowing metrics (\ref{metric})
far more general than the classical ones. But for many years these
"speculations " were largely ignored or dismissed as useless. A notable
exception was Clifford, who with his paper \emph{On the space-theory of
matter} (1870) identified energy and matter with two types of curvature of
space. Many decades later his ideas were found to be important for the
development of Einstein's general relativity theory.

\subsection{The conception of higher-dimensional geometry}

Let us also briefly consider how the conception of dimension has influenced
the developments of geometric ideas in the past. Obviously, our geometric
intuition is largely based upon our experiences with the physical space we
live in, so it can be hard to imagine geometry in more than three dimensions.
Before 1870, the term geometry was in fact largely synonymous with the
classical geometries (\ref{spaces}) and projective geometry in dimensions $2$
and $3$. Geometry was still considered as "descriptive" and supposed to
describe physical space, or at least some idealization or mental conception of
it. L\"{u}tzen [1995] argues that the only notable exceptions were the
non-Euclidean geometers and certainly also Riemann.

Descartes and Euler knew how to identify the \emph{Cartesian n-space}
$\mathbb{R}^{n}$ in dimension $n=2$ or $3$ with the Euclidean plane or space,
respectively. Therefore they could apply analytic geometry in two or three
variables to study geometric problems in three dimensions or less. Conversely,
geometry was used to illustrate analytic problems in at most three variables.
In those days, however, geometry in higher dimensions did not make sense,
although the prominent 17th century scholar Pascal proposed, in fact, that a
4th dimension was allowed in geometry. Making a leap to M\"{o}bius in 1827, we
find that he introduced a 4-dimensional space in his paper \emph{Der
barycentrische Calcul}, but still with the reservation that "it cannot be imagined".

However, the \emph{Elements} of $\mathbb{R}^{n}$ are just n-tuples
$x=(x_{1},..,x_{n})$ of numbers, so gradually it became natural to attempt
generalize the ideas of analytic geometry to $n>3$ variables, for example as
generalized coordinates of mechanical systems. Depending on the problem at
hand, $\mathbb{R}^{n}$ would be referred to as the coordinate space or, say,
the number space in n dimensions. So, it is not surprising that early
geometric ideas in higher dimensions appeared in the study of special subsets
or "submanifolds" of $\mathbb{R}^{n}$. For example, Cauchy used geometric
notions to describe submanifolds of $\mathbb{R}^{n}$ of type $f(x_{1}%
,x_{2},..,x_{n})=0$.

Foremost scholars in the development of multidimensional geometry are Cayley,
Grassmann, Schl\"{a}fli, and Riemann. Grassmann was the first who attempted,
with his \emph{Ausdehnungslehre }(1844 and 1862), a systematic study of
n-dimensional vector spaces from a geometric viewpoint. But his work was hard
to understand and also too philosophical, so it was not appreciated until the
end of the century, when modern tensor calculus was developed. The Swiss
mathematician Ludwig Schl\"{a}fli (1814-1895) studied the differential
geometry of non-linear submanifolds of $\mathbb{R}^{n}$, such as ellipsoids
and their geodesics, and polytopes in higher dimensions, for example the
description of all regular solids in four dimensions. But being ahead of his
time, he too had problems when he tried to publish his major work
\emph{Theorie der vielfachen Kontinuit\"{a}t }(1852). Its importance was fully
appreciated when it was finally published at the turn of the century.

In 1843 Hamilton discovered the quaternions, which initiated a 4-dimensional
geometry, and Cayley announced his interest in multidimensional geometry, at
least in the title of his paper \emph{Chapters in analytic geometry of (n)
dimensions}. Like Grassmann he arrived independently at the notion of an
n-dimensional space, and joined by Sylvester they became the leading British
advocates of n-dimensional geometry. But they were widely opposed, until the
breakthrough at the end of the 1860's, mainly due to the gradual acceptance of
non-Euclidean geometry in France, Germany and Italy, and moreover, the
publication in 1868 of Riemann's celebrated lecture from 1854.\ Here Riemann's
space models are referred to as n-fold extended quantities or manifolds, whose
local coordinate systems make them look locally like $\mathbb{R}^{n}$. In
fact, Riemann had anticipated the idea of an n-dimensional differentiable
manifold, a concept belonging to the new discipline of the 20th century called topology.

Schl\"{a}fli and Riemann both extended the geometry of the Euclidean 2-sphere
$S^{2}$ to higher dimensions, thus providing the following model for Riemann's
n-dimensional spherical space imbedded in Euclidean (n+1)-space
\begin{equation}
S^{n}(r):x_{1}^{2}+x_{2}^{2}+...+x_{n+1}^{2}=r^{2} \label{sphere}%
\end{equation}
Riemann suggested in 1854 that the 3-sphere $S^{3}(r)$, with its intrinsic
geometry of constant positive curvature $K=r^{-2}$, provides a possible model
for the universe which is finite but still unbounded. In contrast to this, a
hyperbolic space model would yield an infinitely large universe. However,
since the curvature could be arbitrary small, it would be almost impossible to
distinguish a large 3-dimensional sphere from Euclidean 3-space. The 20th
century physicist Max Born has described this sphere model as "one of the
greatest ideas about the nature of the world which has ever been conceived".
In fact, as his first attempt at a cosmological model based on the general
theory of relativity, Einstein chose the 3-sphere as his model of the universe
(cf. Osserman [2005]).

\subsection{Foundations of geometry and the geometry of Space}

During the three decades prior to 1860 or so, the strictly synthetic approach
to geometry with Steiner in the forefront had a strong position in German
geometry. But in the late 1860's differential geometry came more to the
foreground, and the philosophically oriented essays of Riemann and Helmholtz
were published in 1868. Here we merely point out that Riemann, with his
mathematical formulation of the concept of space, paved the way for applying
geometry to physical reality, whereas Helmholtz, as an exponent for physical
geometry, became identified throughout the century with the philosophical
foundations of physical space. On the other hand, from a synthesist's
viewpoint, the foundations of geometry, projective or classical, still needed
a thorough revision of its basic concepts and postulates. Indeed, there were
also geometers who opposed the idea of reducing geometric thinking to analytic
geometry and perhaps relying too much on physical intuition or experience.

Moritz Pasch and David Hilbert were the major figures who resurrected the
classical synthetic or "elementary" geometry, to become a rigorous axiomatic
system with emphasis on the purely formal character of geometry. First of all,
the focus was on real projective geometry, considered to be the most basic
geometry, amply demonstrated in works by the French pioneers and by Steiner,
von Staudt, Cayley, and Klein. To ensure the desired properties of a line, for
example, Pasch was the first who recognized the importance of the relations of
"order" and "betweenness". With his book [1882]\emph{\ }he published the first
rigorous deductive system of geometry in history. But despite his insistence
on pure logical reasoning, Pasch still viewed geometry as the science of
physical space, and he wanted to justify the geometric axioms from experience.

In 1895 David Hilbert left K\"{o}nigsberg and joined Klein at G\"{o}ttingen
University. Influenced by Pasch's ideas he became heavily involved with the
foundations of elementary geometry, a subject which had also attracted him
while he was in K\"{o}nigsberg. Hilbert's seminal book, \emph{Grundlagen der
Geometrie}, first published in 1899 and was to appear in many editions, is a
purely logical approach to geometry. In comparison with Euclid's five
postulates, Hilbert formulated about 20 axioms, grouped into five types and
introduced stepwise, in order to clarify what type of theorems can be proved
at each step. These five types are nowadays referred to as the axioms of
incidence, order, congruence, continuity, and parallels.

In addition to Pasch, Hilbert also benefitted from many others who had
initiated techniques which he found useful for his program, as for example
Gergonne and his use of implicit definitions, von Staudt who created a
calculus of line segments, H. Wiener (1857-1939) and his lectures in 1891 on
the role of incidence theorems exemplified by the theorems of Pappus and
Desargues, and finally the Italian geometers G. Peano (1858-1932), M. Pieri
(1860-1913), and A. Padoa (1868-1937), who focused on the strictly logical
approach and replaced the processes of reasoning by symbols and formulas. As a
consequence, the meaning of the fundamental concepts, or the sense in which
the axioms are true, should be excluded altogether from geometry. The story
goes that Hilbert, on a certain occation (see Bos [1993]), expressed this
wisdom by saying "instead of points, lines, planes, one should always be able
to say tables, chairs, beer mugs".

The wave caused by Hilbert's "Grundlagen der Geometrie" pervaded the research
and teaching of geometry during the first half of the 20th century, but after
the initial excitements it also seemed that Hilbert had "killed" the subject,
in a tradition which had little to offer with regard to philosophical or
foundational reflection. On the other hand, in the tradition of Riemann,
Helmholtz and Poincar\'{e}, the term "Grundlagen der Geometrie"\ had in fact
been used several times before, for example in works of Lie and Killing.

Poincar\'{e} and Hilbert were the leading figures in mathematics at the turn
of the century. Poincar\'{e} wrote a favorable review [1903] of Hilbert's work
on the foundations of geometry, and as a great philosopher and scientist he
developed his geometrical conventionalism, closely linked to his own
mathematical studies. In his book [1902], Poincar\'{e} gives his answer to an
aged problem by declaring that the question of whether Euclidean geometry is
true has no meaning, and moreover, one geometry cannot be more true than
another, it can only be more convenient (cf. Torretti [1978], Bos [1993]).

\subsection{Riemann--Helmholtz--Lie space problem}

The geometric interests of Helmholtz arose from his work on physiological
optics in Heidelberg during the early 1860's. He had gained the reputation as
a leading world scientist, whose interests embraced all the sciences, as well
as philosophy and the fine arts. Around 1866, when he moved more towards
physics, he questioned the foundations of geometry, viewed as a science of
physical space. His first brief report \emph{On the factual foundations of
geometry} (1868) appeared in Heidelberg, but his more detailed essay [1868]
was printed in the G\"{o}ttinger Nachrichten shortly afterwards. In the
meantime, Helmholtz had obtained a copy of Riemann's essay from Ernst
Schering, who was Riemann's successor in G\"{o}ttingen. We refer to Nowak
[1989: 43] for a plausible reason why the titles of the two essays which
appeared in G\"{o}ttingen in 1868 differ literally by only one word - "facts"
versus "hypotheses"- as if Helmholtz aimed at a philosophical reproof of
Riemann's view, see Freudenthal [1964].\ 

Helmholtz renounced the conventional attitude to the basic geometric
questions, and like an empiricist he wanted to investigate the nature of space
on the basis of experimental facts (Thatsachen). Since he regarded
differential (or integral) quantities as derived concepts, Helmholtz could not
start from a hypothesis on\ the line element $ds$ as Riemann did. But he
wanted to support Riemann's assumption concerning the nature of $ds$. For him
the most basic observed fact is the free mobility of rigid bodies, and nobody
before him seems to have used mathematics to analyze the logical consequences
of this.

Helmholtz formulated four axioms $H1$ -$H4$, the first of which is similar to
Riemann's postulate that space is (in modern language) an n-dimensional
manifold with differentiability properties. For simplicity, we can say the
second axiom $H2$ expresses the idea of a metric (distance function), and the
motions of the space are the transformations which preserve both the metric
and the orientation. $H3$ explains the meaning of "free movability" (see
(\ref{freemobile})). In fact, a body moves due to the motion of the space, so
the whole space itself is like a rigid body. Helmholtz assumed that by fixing
n-1 general points of a body in n-dimensional space, the remaining mobility is
restricted to a 1-parameter family of motions regarded as "rotations". As a
consequence, the position of a rigid body, generally speaking, depends on
n(n+1)/2 quantities. Finally, $H4$ is the \emph{monodromy axiom}, according to
which the above 1-parameter family of motions is periodic, so that the body
eventually returns to its initial position as the parameter increases to a
certain value.

From his axioms Helmholtz deduced, in fact, Riemann's postulate about the
squared line element $ds^{2}$ (\ref{metric}). On the other hand, like Riemann
he readily accepted that free movability implies constant curvature, from
which he concluded that physical space is either Euclidean or spherical.
Furthermore, assuming the space is infinitely large, as was generally
believed, Helmholtz wrongfully concluded that its curvature must be zero and
is therefore Euclidean. But Beltrami pointed out to him the omission of the
hyperbolic geometry, whose curvature is negative, so Helmholtz issued in 1869
a correction to his paper. Now, from the simple observation that rigid bodies
exist in the space we live in, Helmholtz arrived at the final solution of his
space problem, which we may state as follows (see also Torretti [1978], Chap.
3): \textit{Physical space is a 3-dimensional Riemannian manifold with
constant curvature.}

In 1870 Helmholtz accepted a chair in physics at the university in Berlin,
where he became a highly respected colleague of Weierstrass. During the 1870's
he still gave\ lectures on the origin\ and meaning of geometric axioms, some
of which were published later. He was not the first scholar who argued that
the choice between Euclidean and non-Euclidean geometry cannot be resolved by
pure geometry. But as a physicist he pointed out that a change in the geometry
would impose a change in the laws of mechanics, so a given geometry can be
confirmed or refuted by experience. Like Gauss and Lobachevsky many decades
before him, Helmholtz claimed in 1877 that empirical measurements of triangles
will decide this question.

The papers Riemann [1867] and Helmholtz [1868] appealed to many scholars
interested in the foundations of geometry, as they raised deep mathematical
and philosophical questions about the relationship of geometry and the space
we live in. The motions in Helmholtz's approach are distance preserving
transformations, and they form a group, which is a concept Helmholtz never
mentioned and presumably did not know about. Still, he used notions and
methods with a group theoretic flavour, and in this respect he anticipated
Klein and Lie, who were the first to stress the importance of the notion of
groups in geometry. It was in the fall 1872 that Klein presented his Erlanger
Programm on this topic, whereas Lie started to develop his theory of
continuous groups in the fall 1873.

Lie was first informed by Klein in the 1870's about the geometrical works of
Helmholtz, but he was not attracted by philosophical speculations concerning
the foundation of geometry or the nature of space. On the other hand, it
seemed that noboby had so far really questioned the validity of Helmholtz's
mathematical reasoning and the role of his axioms. So Lie became interested in
the geometric problem behind the idea of "free movability", which he referred
to as the Riemann-Helmholtz space problem. Encouraged by Klein he realized the
problem was well suited for a demonstration of the power of his group theory.

Lie came to Leipzig in the spring 1886, having accepted the chair of geometry
after Klein, who had moved to G\"{o}ttingen. Now Lie was invited to Berlin
where the great Meeting of the German Natural Sciences would take place in
late September. It seems, however, that he accepted the invitation primarily
because of the opportunity to meet Klein in Berlin. Lie's lecture on September
21 was entitled "Tatsachen, welche der Geometrie zu Grunde liegen", or "Facts
lying at the foundation of geometry", and in this address he publicly
criticized Helmholtz's famous 1868 paper, even with sharp comments. Lie
claimed the monodromy axiom $H4$ was superfluous, at least if the axiom $H3$
on "free movability" is interpreted in the proper way.

Lie argued that his theory of continuous groups would give a more
comprehensive and better solution\ of the basic geometric problems discussed
by Helmholtz. His lecture appeared later in 1886 as a note of 5 pages in the
Leipziger Berichte, where also Lie's more detailed elaboration \emph{\"{U}ber
die Grundlagen der Geometrie} appeared as two articles in 1890 and one in
1892. Lie's five papers on the topic are collected in GA II, pp. 374- 479.
However, in Lie [1893] a complete recasting of his solution of the space
problem is presented. It should be mentioned that Poincar\'{e} had, in fact,
solved the special case of dimension $n=2$ in 1887, also by the use of
continuous groups.

Lie worked with infinitesimal groups (Lie algebras) $\mathcal{G}$ rather than
the corresponding continuous groups (Lie groups) $G$, so first of all he
introduced an infinitesimal version of \ "free movability", which enabled him
to determine the possible local structure of the group $G$ of motions
(isometries). This amounts to the determination of infinitesimal generators
$A_{i}$, that is, a basis of $\mathcal{G}$ viewed as a vector space. In Lie's
group theory, the \emph{Elements} of $\mathcal{G}$ can be interpreted as
vector fields acting on the underlying space. Lie first studied the case of a
3-dimensional geometry, as Helmholtz did, and he showed the condition of free
(or maximal) movability implies $\mathcal{G}$ has dimension 6, and moreover,
in terms of suitable coordinates there are just the three cases corresponding
to the classical geometries $E^{3},S^{3},H^{3}$, see (\ref{spaces}).

\ \ \ \ \ \ \ \ \ \ \ \ \ \ 

{\Large Remarks on homogeneous spaces}

Let us also describe the above result globally from a modern viewpoint, as a
pair $(G,M)$ where $G$ is the connected isometry group of a 3-dimensional
Riemannian manifold $M$. "Free movability" means, first of all, that $G$ is
transitive, namely each point $p$ can be moved to any other point $p^{\prime}%
$. Therefore, $M\simeq G/H$ is a homogeneous space, where $H$ $\subset G$ is
the subgroup keeping $p$ fixed. Moreover, $H$ identifies with a subgroup of
the group $\simeq$ $O(3)$ of local isometries around $p$. Now, maximal
movability around $p$ means $H$ should be transitive on the 2-sphere of points
at a fixed small distance from $p$, consequently $H\simeq SO(3)$ or $O(3)$
and, in particular, $G$ has dimension $6$. But the possibilities for such a
pair $H\subset G$ of Lie groups is very limited. With a fairly standard
notation for the groups involved, the solutions representing the classical
geometries (\ref{spaces}) and elliptic geometry, may be stated for n = 3 as
the following homogeneous spaces
\begin{equation}
E^{3}=\frac{\mathbb{R}^{3}\hat{\times}SO(3)}{SO(3)}\text{, \ }S^{3}%
=\frac{SO(4)}{SO(3)}\text{, \ }P^{3}=\frac{SO(4)}{O(3)}\text{, \ }H^{3}%
=\frac{SO(3,1)^{+}}{SO(3)}\text{, } \label{spaces3}%
\end{equation}
and these are, in fact, the only possibilities. They represent Euclidean,
spherical, elliptic, and hyperbolic geometry, respectively.

\ \ \ \ \ \ \ \ \ \ \ \ \ 

The above spaces, however, are not the only onces having constant curvature,
even among "well-behaved" manifolds with no boundary (cf. e.g. Section 5.1).
But ideas related to local and global topological properties of a space were
poorly understood at the time of Lie. In fact, topology and the modern theory
of Lie groups and their homogeneous spaces became theories first in the 20th
century. On the other hand, although Lie's solution of the space problem is
close to a satisfactory modern solution, he did not remove the
differentiability assumptions which sound rather artificial to modern taste,
and this remained a major obstacle at the turn of the century.

In the following years various mathematicians, including David Hilbert and
Hermann Weyl (1855-1955), contributed to the new formulation and final
solution of the problem, which was not found until 1953, by the Belgian
mathematician Jacques Tits (1930-). With rather weak topological assumptions,
the solution asserts that if a triple of points can be carried by a motion
into any other triple having the same mutual distances, then the space is one
of the classical geometries (\ref{spaces}) in some dimension n, more
precisely, either Euclidean, spherical or elliptic, or hyperbolic n-space. In
these spaces the above property is, in fact, the $SSS$-congruence property for
triangles, which is now seen to characterize the classical geometries uniquely.

\bigskip

\end{document}